\documentclass[11pt]{article}
\usepackage{amssymb,amsmath}
\usepackage{graphicx}

\newcommand{\cl}{\mathop{\rm cl}}
\newcommand{\inter}{\mathop{\rm Int}}
\newcommand{\conv}{\mathop{\rm conv}}

\newcommand{\dist}{\mathop{\rm dist}}

\newcommand{\Id}{\mathop{\rm Id}}
\newcommand{\rank}{\mathop{\rm rank}}

\renewcommand{\Im}{\mathop{\rm Im}}
\newcommand{\Proj}{\mathop{\rm Proj}}
\newcommand{\rot}{\mathop{\rm rot}}

\numberwithin{equation}{section}

\newtheorem{lemma}{Lemma}[section]
\newtheorem{theorem}{Theorem}[section]
\newtheorem{corollary}{Corollary}[section]

\title{Problems of unique determination of domains by the
relative metrics on their boundaries
\footnote{
(primary); Key words: intrinsic metric,
relative metric of boundary, local isometry of boundaries,
strict convexity}}
\author{Anatoly~P.~Kopylov\footnote{Sobolev Institute of Mathematics,
Acad. Koptyuga pr. 4, and Novosibirsk State University, Pirogova
str., 2, 630090 Novosibirsk, Russia; apkopylov@yahoo.com}}

\begin{document}

\maketitle

\begin{abstract}

This survey is devoted to discussing the problems of the unique
determination of surfaces that are the boundaries of
(generally speaking) nonconvex domains.

First (in Sec. 2) we examine some results on the problem of the unique determination
of domains by the relative metrics of the boundaries. Then, in Sec. 3, we study
rigidity conditions for the boundaries of submanifolds in a  Riemannian manifold.
The final part (Sec. 4) is concerned with the unique determination of domains
by the condition of the local  isometry of boundaries in the relative metrics.

\end{abstract}

\tableofcontents

\section{Introduction}\label{s1}

The topic of the article, although relatively new, has a
straightforward connection to the classical problems of a
bicentennial history. As the starting point we may view the
celebrated Cauchy theorem which claims that a convex polyhedron
is uniquely determined from its unfolding. Later, the problems
of unique determination of convex surfaces were studied by
Minkowski, Hilbert, Weyl, Blashke, Cohn-Vossen, and other
prominent mathematicians. The greatest progress in this
direction was achieved by A.~D.~Alexandrov and his students.
Mention Pogorelov's celebrated theorem on unique determination
of a closed convex surface from its intrinsic metric (for
example, see~\cite{Po}).

In~\cite{Ko3}, we proposed a new approach to the problem of
unique determination of surfaces, which allowed to substantially
enlarge the framework of the problem. The following model
situation illustrates the essence of this approach fairly well.

Let
$U_1$
and
$U_2$
be two domains (i.e., open connected sets) in the real
$n$-dimensional
Euclidean space
$\mathbb R^n$
whose closures
$\cl U_j$,
where
$j = 1,2$,
are Lipschitz manifolds (such that
$\partial (\cl U_j) = \partial U_j \ne\varnothing$
where
$\partial E$
is the boundary of
$E$
in
$\mathbb R^n$).
Assume also that the boundaries
$\partial U_1$
and
$\partial U_2$
of these domains, which coincide with the boundaries of the
manifolds
$\cl U_1$
and
$\cl U_2$,
are isometric with respect to their relative metrics
$\rho_{\partial U_j,U_j}$
($j = 1,2$),
i.e., with respect to the metrics that are the restrictions to
the boundaries
$\partial U_j$
of the extensions (by continuity) of the intrinsic metrics of the
domains
$U_j$
to
$\cl U_j$.
The following natural question arises: \textit{Under which
additional conditions are the domains
$U_1$
and
$U_2$
themselves isometric {\rm(}in the Euclidean metric{\rm)}}?
In particular, the natural character of this problem is
determined by the circumstance that the problem of unique
determination of closed convex surfaces mentioned in the
beginning of the article is its most important particular case.
Indeed, assume that
$S_1$
and
$S_2$
are two closed convex surfaces in
$\mathbb R^3$,
i.e., they are the boundaries of two bounded convex domains
$G_1 \subset \mathbb R^3$
and
$G_2 \subset \mathbb R^3$.
Let
$U_j = \mathbb R^3 \setminus \cl G_j$
be the complement of the closure
$\cl G_j$
of the domain
$G_j$,
$j = 1,2$.
Then the intrinsic metrics on the surfaces
$S_1 = \partial U_1$
and
$S_2 = \partial U_2$
coincide with the relative metrics
$\rho_{\partial U_1,U_1}$
and
$\rho_{\partial U_2,U_2}$
on the boundaries of the domains
$U_1$
and
$U_2$,
and thus, the problem of unique determination of closed convex
surfaces by their intrinsic metrics is indeed a particular case
of the problem of unique determination of domains by the
relative metrics on their boundaries.

The generalization of the problem of the unique determination of
surfaces ensuing from a new approach suggested in~\cite{Ko3}
manifests itself in the fact that the unique determination of
domains by the relative metrics on their boundaries holds not
only when their complements are bounded convex sets but, for
example, also in the following cases.

\textit{The domain
$U_1$
is bounded and convex, while the domain
$U_2$
is arbitrary} (A.~P.~Kopylov (see~\cite{Ko3})).

\textit{The domain
$U_1$
is strictly convex and the domain
$U_2$
is arbitrary} (A.~D.~Aleksandrov (see~\cite{Al})).

\textit{The domains
$U_1$,
$U_2$
are bounded and their boundaries are smooth} (V.~A.~Aleksandrov
(see~\cite{Al})).

\textit{The domains
$U_1$,
$U_2$
have nonempty bounded complements, while their boundaries are
($n - 1$)-dimensional
connected
$C^1$-manifolds
without boundary{\rm,}
$n > 2$} (V.~A.~Aleksandrov (see~\cite{Al1})).

Recall some definitions of unique determination that are needed
below (for example, see~\cite{Ko1}).

\textbf{Definition~1.1.}
Let
$\mathcal U$
be a class of domains (i.e., open connected sets) in real
Euclidean
$n$-dimensional
space
$\mathbb R^n$,
where
$n \ge 2$.
We say (see, e.g.,~\cite{Ko1}) that a domain
$U \in \mathcal U$
is uniquely determined in the class
$\mathcal U$
by the relative metric of its (Hausdorff) boundary if each
domain
$V \in \mathcal U$
whose Hausdorff boundary is isometric to the Hausdorff boundary
of the domain
$U$
with respect to the relative metrics is itself isometric to
$U$
(with respect to Euclidean metric).

\textbf{Remark~1.1.}
Suppose that
$U$
is a domain in
$\mathbb R^n$
($n \ge 2$)
and
$\rho_U$
is its intrinsic metric.\footnote{Recall that the distance
between points
$x,y \in U$
in the intrinsic metric of
$U$
equals the infimum of the lengths of the curves joining
$x$
and
$y$
within
$U$.}
Consider the Hausdorff completion
$(U_H,\rho_{U_H})$
of
the metric space
$(U,\rho_U)$,
i.e., the completion of this space in intrinsic metric~
$\rho_U$.
Identifying the points of this completion that correspond to
points of the domain
$U$
with these points themselves and removing them from the
completion, we obtain a metric space
$(\partial_H U,\rho_{\partial_H U,U})$;
the set
$\partial_H U$
of its elements is called the Hausdorff boundary of the domain
$U$,
and
$\rho_{\partial_H U,U}$
is the relative metric on this Hausdorff boundary. The isometry
of the Hausdorff boundaries of domains
$U$
and
$V$
with respect to their relative metrics means the existence of a
surjective isometry
$f: (\partial_H U,\rho_{\partial_H U,U}) \to
(\partial_H V,\rho_{\partial_H V,V})$
between these boundaries.

\textbf{Remark~1.2.}
If the intrinsic metric
$\rho_U$
of a domain
$U$
extends by continuity to the closure
$\cl U$
then
$\partial_H U$
is naturally identified with the Euclidean boundary
$\partial U$.

It is also worth noting that M.~V.~Korobkov obtained a complete
description of domains that are uniquely determined in the class
of all domains by the condition of the (global) isometry of their
boundaries in the relative metrics
(see~\cite{Kor1},~\cite{Kor2},~\cite{Kor4},~\cite{Kor3}\footnote{In
these works, the isometry
of the Hausdorff boundaries of domains
$U$
and
$V$
with respect to their relative metrics means the existence of a
bijective isometry
$f: (\partial_H U,\rho_{\partial_H U,U}) \to
(\partial_H V,\rho_{\partial_H V,V})$
between these boundaries.}).
In~\cite{Kor1},~\cite{Kor2}, he
obtained the following results.

Henceforth
$\inter E$ is the interior of a set
$E$,
while
$\cl E$
is the closure of
$E$,
$\partial E$
is the boundary of
$E$,
and
$\conv E$
is the convex hull of
$E$.

Given a domain
$U$,
we put
$F_U = U \setminus \cl (\conv \partial U)$
and denote by
$U_i$
the connected components of the open set
$U \cap \inter (\conv \partial U)$.

\begin{theorem}\label{t1.1}
Let
$U$
be a domain in
$\mathbb R^2$.
Then{\rm:}

$\rm{(I)}$
If
$\partial U$
lies on some straight line then
$U$
cannot be uniquely determined {\rm(}in the class of all domains in
$\mathbb R^2)$
if and only if
$\partial U$
is a connected set containing more than one point;

$\rm{(II)}$
If there is no straight line containing
$\partial U$
then
$U$
cannot be uniquely determined if and only if there exist a domain
$V$
nonisometric to
$U,$
a family of isometric mappings
$Q_i: \mathbb R^2 \to \mathbb R^2,$
and a homeomorphism
$\theta: \partial F_U \to \partial F_V$
satisfying the following conditions{\rm:}

$\rm{(IIa)}$
$Q_i(U_i) = V_i$
for each component
$U_i$.
The same equality is also valid for each component~
$V_i$.

$\rm{(IIb)}$
$x \in U \cap \partial U_i$
if and only if
$Q_i(x) \in V \cap \partial V_i;$
the validity of these containments implies that
$\theta(x) = Q_i(x)$.

$\rm{(IIc)}$
The homeomorphism
$\theta$
preserves the arc length {\rm(}i.e.{\rm,} for arbitrary two
points
$x$
and
$y$
the length of the arc joining
$x$
and
$y$
within
$\partial F_U$
coincides with the length of the image of this arc under the
mapping
$\theta)$.
\end{theorem}

As a straightforward consequence of Theorem~\ref{t1.1} we have

\begin{corollary}\label{c1.1}
Suppose that a domain
$U \subset \mathbb R^2$
satisfies one of the following two conditions{\rm:}

${\rm(1)}$
$F_U = \varnothing;$

${\rm(2)}$
the set
$F_U$
is nonempty and disconnected.

Then
$U$
is uniquely determined from the relative metric of its Hausdorff
boundary.
\end{corollary}

In turn, Corollary~\ref{c1.1} contains the following particular
case.

\begin{corollary}\label{c1.2}
Suppose that a domain
$U \subset \mathbb R^2$
is bounded. Then
$U$
is uniquely determined from the relative metric  of its
Hausdorff boundary.
\end{corollary}

The following theorem generalized A.~D.~Aleksandrov's theorem
about the unique determination of a strictly convex domain in
the class of domains whose intrinsic metrics extend by
continuity to their closure (see~\cite{Al}) to the case of
convex domains. In the plane case, this theorem ensues from
Corollary~\ref{c1.1}.

\begin{theorem}\label{t1.2}
Each convex domain
$U \subset \mathbb R^n,$
$n \ge 2,$
different from an open half-space in
$\mathbb R^n$
is uniquely determined from the relative metric of its Hausdorff
boundary.
\end{theorem}

In~\cite{Kor4},~\cite{Kor3}, Korobkov obtained a result which is
analogous to Theorem~\ref{t1.1} and contains a complete
description domains in
$\mathbb R^n$,
$n \ge 3$,
that are uniquely determined in the class of all
$n$-dimensional
domains by the condition of the (global) isometry of their
Hausdorff boundaries in the relative metrics.

In Section~2, we consider the problems of unique determination
of domains by the condition of the (global) isometry of their
boundaries in the relative metrics. We state of the proofs of
Theorems~\ref{t1.1},~\ref{t1.2}.

In this connection, there appears the following question: Is it
possible to construct an analog of the theory of rigidity of
surfaces in Euclidean spaces in the general case of the
boundaries of submanifolds in Riemannian manifolds?

Section~3 of our article is devoted to a detailed discussion of
this question. In it, we in particular obtain new results
concerning rigidity problems for the boundaries of
$n$-dimensional
connected submanifolds with boundary in
$n$-dimensional
smooth connected Riemannian manifolds without boundary
($n \ge 2$).

Results of~\cite{Kor1},~\cite{Kor2},~\cite{Kor4},~\cite{Kor3} imply, in
particular, that any bounded domain in
$\mathbb R^n$
is uniquely determined by the condition of isometry of boundaries
in the relative metrics. At the same time, according to results
of~\cite{Bor}, a bounded polygonal plane domain
$U$
is uniquely determined by the condition of local isometry of
boundaries in the relative metrics in the class of all such
domains if and only if the domain
$U$
is convex.

\textbf{Remark~1.3.}
Let
$\mathcal M$
be a class of domains in space
$\mathbb R^n$
with
$n \ge 2$.
Following~\cite{Ko1}, we say that a domain
$U \in \mathcal M$
is uniquely determined in the class
$\mathcal M$
by the condition of local isometry of the (Hausdorff) boundaries
of domains in the relative metrics if, for any domain
$V$
belonging to the class
$\mathcal M$,
the local isometry of its Hausdorff boundary to the Hausdorff
boundary of the domain
$U$
with respect to the relative metrics implies the isometry of the
domains
$U$
and
$V$
(with respect to the Euclidean metric). The local isometry in
the relative metrics between the Hausdorff boundaries
$\partial_H U$
and
$\partial_H V$
of the domains
$U$
and
$V$
means the existence of a bijective mapping
$f: \partial_H U \to \partial_H V$
of these boundaries which is a local isometry with respect to
their relative metrics, i.e., a mapping such that, for any
element
$y \in \partial_H U$,
there exists a number
$\varepsilon > 0$
satisfying the following condition: for any two elements
$a$
and
$b$
from the
$\varepsilon$-neighborhood
$Z(y) = \{z \in \partial_H U: \rho_{\partial_H U,U}(z,y) < \varepsilon\}$
of
$y$,
$\rho_{\partial_H U,U}(a,b) = \rho_{\partial_H V,V}(f(a),f(b))$.
It is clear that
$f^{-1}$
is also a local isometry with respect to relative metrics of
boundaries.

In Section~4 of this paper, we continue the study of the unique
determination of domains by the condition of local isometry of
their boundaries in the relative metrics.

It can be divided into two parts.

The first of them is mainly devoted to finding a complete description
of conditions that are necessary and sufficient for a plane
domain with smooth boundary to be uniquely determined by the
condition of local isometry of their boundaries in the class of
all domains with smooth boundaries (in the case of a bounded domain,
in the class of all bounded plane domains with smooth boundaries).

In the second part of Section~4, we obtain some new assertions on
the unique determination of space domains with smooth boundaries
by the considering in the Section condition. All of these results
emphasize the specific character of our approach to the problems
of rigidity of domains in $\mathbb R^n$.

Note that below
$[a,b] = \{bt + (1-t)a \in \mathbb R^n: 0 \le t \le 1 \}$,
$[a,b[ = \{bt + (1-t)a \in \mathbb R^n: 0 \le t < 1 \}$
($]a,b] = \{bt + (1-t)a \in \mathbb R^n: 0 < t \le 1 \}$)
and
$]a,b[ = \{bt + (1-t)a \in \mathbb R^n: 0 < t <1 \}$
are the segment (closed interval), the half-open interval and
the (open) interval in
$\mathbb R^n$
with endpoints
$a,b \in \mathbb R^n$,
$a \ne b$.
$\inter I$
is the interior of the segment (of the half-open interval)
$I$,
$\inter ]a,b[ = ]a,b[$.
$B(x_0,r) = \{x \in \mathbb R^n: |x - x_0| < r\}$
is the open ball in
$\mathbb R^n$
of radius
$r$
($0 <r < \infty$)
centered at
$x_0 \in \mathbb R^n$.
$\Id_E$
is the identity mapping of a set
$E$:
$\Id_E(x) = x$
for all
$x \in E$.

In what follows, all paths (curves)
$\gamma: [\alpha,\beta] \to \mathbb R^n$,
where
$\alpha,\beta \in \mathbb R$,
are assumed continuous and non-constant, and
$l(\gamma)$
means the length of a path
$\gamma$.
If
$\gamma : [\alpha,\beta] \to \mathbb R^n$
is continuous and injective then
$\gamma$
is also called an arc.

\section{On Unique Determination of Domains by Relative Metric
of Boundaries}\label{s2}

Below by connectedness we mean connectedness in the sense of
general topology.

The support of an element
$\widetilde x$
of the Hausdorff boundary
$\partial_H U$
of a domain
$U \subset \mathbb R^n$,
$n \ge 2$,
is a point
$x$
of the Euclidean boundary
$\partial U$
which is the limit of a Cauchy sequence (in the intrinsic metric
of
$U$)
of points
$x \in U$
representing
$\widetilde x$.
We denote the support of
$\widetilde x \in \partial_H U$
by
$p\widetilde x\,\, (= p_H\widetilde x)$.

It is clear that each element
$\widetilde x \in \partial_H U$
possesses the unique support
$x = p\widetilde x$.
At the same time, there can be points
$x \in \partial U$
that are not the support of any element of
$\partial_H U$;
on the other hand, there can be points
$x \in \partial U$
that are the supports of several (even uncountable many) elements
of
$\partial_H U$.
The following simple facts are valid:

\begin{lemma}\label{l2.1}
Let
$U$
be a domain in
$\mathbb R^n$.
Then
$\{p\widetilde x : \widetilde x \in \partial_H U\}$
is everywhere dense in
$\partial U$
{\rm(}with respect to the usual Euclidean metric{\rm)}.
\end{lemma}

\begin{lemma}\label{l2.2}
If
$\widetilde x \in \partial_H U$
then the support
$x = p\widetilde x$
is attainable from
$U$
along a rectifiable curve
$\gamma : [0,1] \to \cl U$
such that
$\gamma(1) = x,$
$\gamma(t) \in U$
for
$0 \le t < 1,$
and{\rm,} for each sequence
$\{t_j\}_{j \in \mathbb N}$
of points
$t_j \in [0,1[$
satisfying the condition
$\lim_{j \to \infty} t_j = 1,$
the sequence
$\{\gamma(t_j)\}_{j \in \mathbb N}$
is Cauchy in the intrinsic metric of
$U$
and represents the boundary element
$\widetilde x$.
On the other hand{\rm,} if a point
$x \in \partial U$
is attainable from
$U$
along a rectifiable curve
$\gamma$
then every sequence
$\{\gamma(t_j)\}_{j \in \mathbb N},$
where
$t_j \in [0,1[,$
$j = 1,2,\dots,$
and
$\lim_{j \to \infty} t_j = 1,$
is Cauchy in the intrinsic metric of
$U$
and determines the only element
$\widetilde x \in \partial_H U$
whose support is
$x$.
\end{lemma}

We need some additional notions.
Granted an ordered triple
$x,y,z$
of points, denote by
$\angle (x,y,z)$
the plane angle
$\angle (x,y,z) = \{y + t(x - y) + s(z - y) : t \ge 0, s \ge 0\}$;
and by
$\Delta(x,y,z)$,
the (closed) triangle with vertices
$x$,
$y$,
and
$z$.
We denote by
$|\angle(x,y,z)|$
the value of the angle
$\angle(x,y,z)$
(in radians).

Let
$U$
be a domain in
$\mathbb R^n$
different from the whole
$\mathbb R^n$
(we further consider only the domains of this kind). An
interval
$]x,y[$
is called a boundary interval for
$U$
if
$\varnothing \ne\,\, ]x,y[\,\, \subset U$
and
$x,y \in \partial U$.
We say that an ordered triple
$x,y,z$
of points determines a boundary angle for the domain
$U$
if these points do not lie on one straight line,
$x,y,z \in \partial U$,
$]x,y[\,\, \subset U$,
$]y,z[\,\, \subset U$,
and
$\exists r > 0$
$B(y,r) \cap \angle(x,y,z) \setminus \{y\} \subset U$.
Henceforth we say that a triangle
$\Delta(x,y,z)$
is a boundary triangle for
$U$
if
$x$,
$y$,
and
$z$
do not lie on one straight line,
$x,y,z \in \partial U$,
and
$\Delta(x,y,z) \setminus \{x,y,z\} \subset U$.
It is obvious that each side of a boundary triangle is a boundary
interval and the angles of this triangle are boundary angles.

Also, it is easily seen that each boundary interval
$]x,y[$
naturally generates a pair
$\widetilde x,\widetilde y \in \partial_H U$
of elements of the Hausdorff boundary such that
$p\widetilde x = x$,
$p\widetilde y = y$,
and
$\rho_{\partial_H U,U}(\widetilde x,\widetilde y) = |x - y|$
(the elements
$\widetilde x$
and
$\widetilde y$
are generated by Cauchy sequences of points in
$]x,y[$
converging to the respective points
$x$
and
$y$).
Similarly, each boundary angle
$\angle(x,y,z)$
naturally generates a triple
$\widetilde x,\widetilde y,\widetilde z$
of elements of the Hausdorff boundary
$\partial_H U$
such that
$p\widetilde x = x$,
$p\widetilde y = y$,
$p\widetilde z = z$,
$\rho_{\partial_H U,U}(\widetilde x,\widetilde y) = |x - y|$,
and
$\rho_{\partial_H U,U}(\widetilde z,\widetilde y) = |z - y|$.
The same can be said about the boundary triangle
$\Delta(x,y,z)$
for which we also have
$\rho_{\partial_H U,U}(\widetilde x,\widetilde z) = |x - z|$.

Assume that
$f : \partial_H U \to \partial_H V$
is a mapping,
$x \in \partial U$,
$\widetilde x \in \partial_H U$,
and
$p\widetilde x = x$
(usually, by default the point
$\widetilde x$
is determined from the context, for example, as the endpoint of
a boundary interval; see above). Now put
$x' = p f(\widetilde x)$.
Generally, given a point
$x \in \cl U$, we denote by
$x'$
the corresponding point
$x' \in \cl V$
when the correspondence is clear from the context.

Clearly, the intrinsic metric of
$U$
can be extended by continuity to
$U \cup \partial_H U$.
Thus, the value
$\rho_{\partial_H U,U}(\widetilde x,\widetilde y)$
is defined for every pair
$x,y$
in
$U \cup \partial_H U$.

We need a series of lemmas; moreover, the main tools here are
Lemmas~\ref{l2.3} and~\ref{l2.9}.

\begin{lemma}\label{l2.3} {\rm(on invariance of a boundary
interval).}
Suppose that
$U,V \subset \mathbb R^n$
are domains and
$f : \partial_H U \to \partial_H V$
is an isometry {\rm(}in the relative metrics{\rm)} of the
Hausdorff boundaries of these domains. Suppose that
$]x,y[$
is a boundary interval for the domain
$U$.
Then
$]x',y'[$
is a boundary interval for
$V;$
moreover{\rm,}
$|x' - y'| = |x - y)$.
\end{lemma}

Some particular cases of this assertion were known earlier: for
example, in the case when
$U$
is a bounded convex domain a similar assertion is contained in
Lemma~4 of~\cite{Ko3}.

\textit{Proof of Lemma~{\rm\ref{l2.3}}.}
Since
$f$
is an isometry, we have
\begin{equation}\label{eq2.1}
\rho_{\partial_H V,V}(\widetilde x,\widetilde y) = |x - y|.
\end{equation}
Take some Cauchy sequences
$v_\nu \in V \to x'$
and
$w_\nu \in V \to y'$
(with respect to the intrinsic metric of
$V$)
generating elements
$f(\widetilde x)$
and
$f(\widetilde y)$
of the Hausdorff boundary. (Existence of these sequences follows
from the definition of Hausdorff boundary.) Now, take a sequence
of curves
$\gamma_\nu : [0,1] \to V$
such that
$\gamma_\nu(0) = v_\nu$,
$\gamma_\nu(1) = w_\nu$,
and
$l(\gamma_\nu) \to \rho_{\partial_H V,VU}(f(\widetilde x),f(\widetilde y)) =
|x - y|$.
Without loss of generality we can also assume that the
parameterizations of
$\gamma_\nu$
coincide to within factor with the natural parameterizations and
the mappings
$\gamma_\nu$
converge uniformly to the mapping
$\gamma : [0,1] \to \cl V$
so that
$\gamma(0) = x'$,
$\gamma(1) = y'$,
$\gamma$
is a rectifiable curve, and
\begin{equation}\label{eq2.2}
l(\gamma) \le \rho_{\partial_H V,V}(f(\widetilde x),f(\widetilde y)).
\end{equation}
To complete the proof of Lemma~\ref{l2.3}, it suffices now to
show that
\begin{equation}\label{eq2.3}
\forall t \in\,\, ]0,1[ \quad \gamma(t) \in V.
\end{equation}
Indeed, if~(\ref{eq2.3}) is valid then, by the definition of the
relative metric on the Hausdorff boundary, we have
$l(\gamma) \ge \rho_{\partial_H V,V}(f(\widetilde x),f(\widetilde y))$.
Combining this with~(\ref{eq2.2}), we obtain
\begin{equation}\label{eq2.4}
l(\gamma) = \rho_{\partial_H V,V}(f(\widetilde x),f(\widetilde y)).
\end{equation}
From~(\ref{eq2.3}) and~(\ref{eq2.4}) we easily find that the
curve
$\gamma$
is a line segment with endpoints
$x'$
and
$y'$;
moreover,
$]x',y'[\,\, \subset V$
and
$|x' - y'| = \rho_{\partial_H V,V}(f(\widetilde x),f(\widetilde y))$.
This together with~(\ref{eq2.1}) gives the claim of
Lemma~\ref{l2.3}.

Thus, we are left with proving~(\ref{eq2.3}). Assume~(\ref{eq2.3})
false. Then there is
$t_0 \in\,\, ]0,1[$
such that
$\gamma(t_0) \in \partial V$.
Put
$r_\nu = l(\gamma_\nu([0,t_0]))$,
$s_\nu = l(\gamma_\nu([t_0,1]))$,
and
$R_\nu = \dist(\gamma_\nu(t_0),\partial V)$.
By the above assumptions,
$R_\nu \to 0$
as
$\nu \to \infty$,
$\forall \nu = 1,2,\dots\,\,\, (B(\gamma_\nu(t_0),R_\nu) \subset V
\&\,\, \exists A_\nu \in \cl B(\gamma_\nu(t_0),R_\nu) \cap
\partial V)$,
and
\begin{equation}\label{eq2.5}
r_\nu +s_\nu \to |x - y| \quad {\rm as}\,\, \nu \to \infty,
\end{equation}
\begin{equation}\label{eq2.6}
\exists \varepsilon_0 > 0\,\,\, \forall \nu\,\, \in \mathbb N\,\,\,
(r_\nu > \varepsilon_0\,\, \&\,\, s_\nu > \varepsilon_0).
\end{equation}

By~(\ref{eq2.5}) and~(\ref{eq2.6}), we can assume without loss
of generality that
$r_\nu \to r > 0$,
$s_\nu \to s > 0$,
and
\begin{equation}\label{eq2.7}
r + s = |x - y|.
\end{equation}
If we take a sequence of points in the ball
$B(\gamma_\nu(t_0),R_\nu)$
converging to
$A_\nu$
then it generates some element of the Hausdorff boundary
$\widetilde A_\nu \in \partial_H V$.
It is easy to see that
$\limsup_{\nu \to \infty} \rho_{\partial_H V,V}(f(\widetilde x),
\widetilde A_\nu)
\le r$
and~\linebreak
$\limsup_{\nu \to \infty}
\rho_{\partial_H V,V}(f(\widetilde y),\widetilde A_\nu)
\le s$.
Denote
$\widetilde B_\nu = f^{-1}(\widetilde A_\nu)$.
Since the mapping
$f$
is an isometry, we find
\begin{equation}\label{eq2.8}
\limsup_{\nu \to \infty}
\rho_{\partial_H V,V}(\widetilde x,\widetilde B_\nu)
\le r,
\end{equation}
\begin{equation}\label{eq2.9}
\limsup_{\nu \to \infty}
\rho_{\partial_H V,V}(\widetilde y,\widetilde B_\nu)
\le s.
\end{equation}
Put
$B_\nu = p\widetilde B_\nu$.
Then, by the definition of
$\rho_{\partial_H U,U}$,
the following are valid:
\begin{equation}\label{eq2.10}
|x - B_\nu| \le \rho_{\partial_H V,V}(\widetilde x,\widetilde B_\nu),
\end{equation}
\begin{equation}\label{eq2.11}
|y - B_\nu| \le \rho_{\partial_H V,V}(\widetilde y,\widetilde B_\nu),
\end{equation}
It follows from~(\ref{eq2.7})-(\ref{eq2.11}) and the triangle
inequality that
$B_\nu \to B = \frac{s}{|x - y|}x + \frac{r}{|x - y|}y$.
Since
$B_\nu \in \partial U$,
we have
$B \in \partial U$.
However, the last fact contradicts the condition
$]x,y[ \subset U$
of Lemma~\ref{l2.3}.
The contradiction completes the proof of Lemma~\ref{l2.3}.

In the lemmas below, we suppose unless the contrary is specified
that
$U$
and
$V$
are domains in
$\mathbb R^2$
and
$f : \partial_H U \to \partial_H V$
is an isometry (in the relative metrics) of the Hausdorff
boundaries of these domains.

\textbf{Remark~2.1.}
With these definitions, the mapping
$f : \partial_H U \to \partial_H V$
is an isometry  of the Hausdorff boundaries if and only if the
inverse mapping
$f^{-1} : \partial_H V \to \partial_H U$
is an isometry of the Hausdorff boundaries. Therefore, the
following lemmas remain valid with
$U$
and
$V$
interchanged.

\begin{lemma}\label{l2.4}
Let
$\angle (x,y,z)$
be a boundary angle for
$U$.
Then
$|x' - y'| = |x - y|,$
$|z' - y'| = |z - y|,$
and
$|x' - z'| \le |x - z|$.
\end{lemma}

We can express the gist of this lemma as follows: The mapping
$f$
takes boundary angles into the angles equal or less than the
original (so far we do not claim that these angles are boundary
angles).

Below we need one more denotation. Let
$\angle (x,y,z)$
be a boundary angle for some domain
$\Omega$.
Put
$E = E(\angle(x,y,z),\Omega) = \conv(\angle(x,y,z) \cap
(\partial \Omega) \setminus \{y\})$.
It follows from the definitions that
$[x,z] \subset E$
and
\begin{equation}\label{eq2.12}
\exists r_1 > 0 \quad B(y,r_1) \cap E = \varnothing.
\end{equation}
Denote
\begin{equation}\label{eq2.13}
\Gamma^{\Omega}_{\angle(x,y,z)} = \{w \in (\Omega) \cap
(\partial E) : (y,w) \subset \Omega \setminus E\}.
\end{equation}
Observe, in particular, that by construction
\begin{equation}\label{eq2.14}
x,z \in \Gamma^{\Omega}_{\angle(x,y,z)},
\end{equation}
\begin{equation}\label{eq2.15}
\forall w \in \Gamma^{\Omega}_{\angle(x,y,z)} \quad\,\,
]y,w[\,\, {\rm is\,\, a\,\, boundary\,\, interval\,\, for}\,\,
\Omega.
\end{equation}

\textit{Proof of Lemma~{\rm\ref{l2.4}}.}
Making parallel translations, if necessary, we can assume without
loss of generality that
\begin{equation}\label{eq2.16}
y = y' = 0 \in \mathbb R^2.
\end{equation}

Let
$\Phi = |\angle(x,0,z)|$
and
$E = E(\angle(x,0,z),U)$
(see above). Construct an arc
$\gamma : [0,\Phi] \to \partial E$,
defining the point
$\gamma(\varphi)$
by the following properties:
$\gamma(\varphi) \in \angle(x,0,z)$,
$\gamma(\varphi) \in E$,
$]0,\gamma(\varphi)[\,\, \cap E = \varnothing$,
and
$|\angle(x,0,\gamma(\varphi))| = \varphi$.
It is clear that
$\gamma(0) = x$,
$\gamma(\Phi) = z$,
and
$\gamma$
is a convex arc in
$E \cap \partial E$
joining
$x$
and
$z$.

Denote
$\Gamma = \gamma([0,1])$,
$\Gamma_1 = \Gamma \cap \partial U$,
and
$\Gamma_2 = \Gamma \setminus \Gamma_1$.
From the above
$\Gamma_1 = \Gamma^U_{\angle(x,0,z)}$,
whence, by~(\ref{eq2.15}),
\begin{equation}\label{eq2.17}
\forall w \in \Gamma_1, \quad\,\, ]0,w[\,\, {\rm is\,\, a\,\,
boundary\,\, interval}.
\end{equation}
Moreover, by construction,
$\forall w \in \Gamma_2, \quad \exists w_1,w_2 \in \Gamma_1 \quad
w \in\,\, ]w_1,w_2[\,\, \subset U \cap \Gamma$.

Define
$g : \Gamma_1 \to \partial V$
by the rule
$g(w) = w'$
for
$w \in \Gamma_1$.
By Lemma~\ref{l2.3} (see~(\ref{eq2.16}) and~(\ref{eq2.17})),
\begin{equation}\label{eq2.18}
\forall w \in \Gamma_1, \quad |g(w)| = |w|.
\end{equation}
From the above constructions we also derive the equalities
\begin{equation}\label{eq2.19}
x' = g(x) , \quad z' = g(z).
\end{equation}
From here and~(\ref{eq2.18}) we obtain in particular that
\begin{equation}\label{eq2.20}
|x'| = |x|, \quad |z'| = |z|.
\end{equation}

Since
$\Gamma$
is a convex arc and all its extreme\footnote{An extreme point of
a set
$F$
is a point
$A \in F$
such that there is no pair of points
$B,C \in F$
different from
$A$
and a number
$p \in\,\, ]0,1[$
for which
$A = pB + (1 - p)C$.}
points lie in
$\Gamma_1$,
it is geometrically obvious that, for each
$\varepsilon > 0$,
there is a partition
$0 < \varphi_0 < \varphi_1 < \dots < \varphi_N = \Phi$,
$N = N(\varepsilon)$,
such that
$\forall \nu = 0,1,\dots,N\,\,\, \gamma(\varphi_\nu) \in
\Gamma_1$,
\begin{equation}\label{eq2.21}
\forall \nu = 0,1,\dots,N - 1 \quad
\frac{l(\gamma([\varphi_\nu,\varphi_{\nu + 1}]))}
{|\gamma(\varphi_\nu) - \gamma(\varphi_{\nu + 1})|} \le
1 + \varepsilon.
\end{equation}
Using the definition of Hausdorff metric, we find that
\begin{equation}\label{eq2.22}
\forall \nu = 0,1,\dots,N - 1 \quad
\rho_{U_H}(\widetilde {\gamma(\varphi_\nu)},
\widetilde {\gamma(\varphi_{\nu + 1})}) \le
l(\gamma([\varphi_\nu,\varphi_{\nu + 1}])).
\end{equation}
From~(\ref{eq2.21}) and~(\ref{eq2.22}) and the definition of
Hausdorff metric we conclude that
\begin{multline}\label{eq2.23}
\forall \nu = 0,1,\dots,N - 1 \quad |g(\gamma(\varphi_\nu)) -
g(\gamma(\varphi_{\nu + 1}))| \le
\rho_{V_H}(f(\widetilde {\gamma(\varphi_\nu)}),
f(\widetilde {\gamma(\varphi_{\nu + 1})})) =
\\
\rho_{U_H}(\widetilde {\gamma(\varphi_\nu)},
\widetilde {\gamma(\varphi_{\nu + 1})}) \le
l(\gamma([\varphi,\varphi_{\nu + 1}])) \le
(1 + \varepsilon) |\gamma(\varphi_\nu) - \gamma(\varphi_{\nu + 1})|.
\end{multline}
From the arbitrariness of
$\varepsilon > 0$,~(\ref{eq2.18}),~(\ref{eq2.23}), and obvious
geometric arguments, we infer that
$|g(z) - g(x)| \le |x - z|$
which together with~(\ref{eq2.19}) and~(\ref{eq2.20}) gives the
claim of Lemma~\ref{l2.4}.

Introduce one more definition. Take
$\widetilde x,\widetilde y \in \partial_H U$.
We say that a rectifiable curve
$\gamma : [0,L] \to \cl U$
is an
$H$-shortest curve between
$\widetilde x$
and
$\widetilde y$
if there is a sequence of rectifiable naturally parameterized curves
$\gamma_\nu : [0,L_\nu] \to U$
such that
$L_\nu \to L = |\widetilde x - \widetilde y|$,
$\sup_{t \in [0,L]}|\gamma_\nu(t\frac{L_\nu}{L}) - \gamma(t)| \to 0$,
$\rho_{\partial_H U,U}(\gamma_\nu(0),\widetilde x) \to 0$,
and
$\rho_{\partial_H U,U}(\gamma_\nu(L_\nu),\widetilde y) \to 0$
as
$\nu \to \infty$.

\begin{lemma}\label{l2.5}
For every pair
$\widetilde x,\widetilde y \in \partial_H U,$
there is an
$H$-shortest
curve. Moreover{\rm,} every
$H$-shortest
curve
$\gamma : [0,L] \to \cl U$
possesses the following properties{\rm:}

{\rm(H-i)}
Given a pair of numbers
$\varphi_1,\varphi_2 \in [0,L],$
if
\begin{equation}\label{eq2.24}
\gamma(\varphi_1),\gamma(\varphi_2) \in \partial U\,\, {\rm and}\,\,
\gamma(\varphi) \in U \,\, {\rm for\,\, all}\,\, \varphi \in\,\,
]\varphi_1,\varphi_2[
\end{equation}
then
$\gamma(\varphi) \in\,\, ]\gamma(\varphi_1),\gamma(\varphi_2)[$
for all
$\varphi \in\,\, ]\varphi_1,\varphi_2[$.

{\rm(H-ii)}
If
$\gamma(\varphi) \in U$
then
$$
l(\gamma([0,\varphi])) \le \rho_{U_H}(\widetilde x,
\gamma(\varphi)) = \varphi, \quad l(\gamma([\varphi,L])) \le
\rho_{U_H}(\widetilde y,\gamma(\varphi)) = L - \varphi.
$$
\end{lemma}

The assertion (H-i) means that the intersection of an
$H$-shortest
curve with
$U$
consists of boundary intervals.

\textit{Proof of Lemma~{\rm\ref{l2.5}}} is simple and therefore
omitted.

\begin{lemma}\label{l2.6}
Take
$\widetilde x,\widetilde y \in \partial_H U$
and let
$\gamma : [0,L] \to \cl V$
be an
$H$-shortest
curve between
$\widetilde x$
and
$\widetilde y$.
Then there is an
$H$-shortest
curve
$\gamma_V : [0,L] \to \cl V$
between
$f(\widetilde x)$
and
$f(\widetilde y)$
with the following properties{\rm:}

{\rm(H-iii)}
$\forall \varphi \in\,\, ]0,L[ \quad (\gamma(\varphi) \in U
\Leftrightarrow \gamma_V(\varphi) \in V)$.

{\rm(H-iv)}
For every pair of numbers
$\varphi_1,\varphi_2 \in [0,L]$
satisfying~{\rm(\ref{eq2.24})} we have{\rm:}
$f(\widetilde {\gamma(\varphi_1)}) =
\widetilde {\gamma_V(\varphi_1)}$
and
$f(\widetilde {\gamma(\varphi_2)}) =
\widetilde {\gamma_V(\varphi_2)}$.
\end{lemma}

The proof relies on the corresponding definitions and
Lemmas~\ref{l2.2},~\ref{l2.3},~\ref{l2.5} and is carried out by
the standard arguments of calculus; therefore, we omit it.

\begin{lemma}\label{l2.7}
Let
$\Delta(x,y,z)$
be a boundary triangle for
$U$.
Then
$\Delta(x',y',z')$
is a boundary triangle for
$V;$
moreover{\rm,}
$$
|x' - y'| = |x - y|,\,\,\, |z' - y'| = |z - y|,\,\,\, |x' - z'| = |x - z|.
$$
\end{lemma}

We can express the essence of this lemma as follows: The mapping
$f$
takes boundary triangles into equal boundary triangles.

\textit{Proof.}
From the definition of a boundary triangle we obtain in particular
\begin{equation}\label{eq2.25}
\inter \Delta(x,y,z) \subset U.
\end{equation}
It follows from Lemma~\ref{l2.3} that
$x',y',z' \in \partial V$ and
\begin{equation}\label{eq2.26}
]x',y'[\,\, \cup\,\, ]y',z'[\,\, \cup\,\, ]x',z'[\,\, \subset V,
\end{equation}
\begin{equation}\label{eq2.27}
|x' - y'| = |x - y|,\,\,\, |z' - y'| = |z - y|,\,\,\, |x' - z'| = |x - z|.
\end{equation}
It remains to prove that
\begin{equation}\label{eq2.28}
\inter \Delta(x',y',z') \subset V.
\end{equation}
Assume that some point
$w$
satisfies the conditions
$w \in \partial U$,
$]y,w[\,\, \subset U$,
and
$]y,w[\,\, \cap\,\, ]x,z[\,\, \ne\varnothing$.
Then it is geometrically obvious that
$\angle(x,y,w)$
and
$\angle(w,y,z)$
are boundary angles. Then, by Lemmas~\ref{l2.3} and~\ref{l2.4},
$]y',w'[\,\, \in V$,
$|w' - x'| \le |w - x|$,
$|w' - z'| \le |w - z|$,
and
\begin{equation}\label{eq2.29}
|w' - y'| = |w - y|.
\end{equation}
Hence, using~(\ref{eq2.27}), we conclude that
\begin{equation}\label{eq2.30}
|w' - x'| = |w - x|,\,\,\, |w' - z'| = |w - z|.
\end{equation}
Equalities~(\ref{eq2.27}),~(\ref{eq2.29}), and~(\ref{eq2.30})
demonstrate that the figure consisting of the four points
$x,y,z$,
and
$w$
is Euclideanly isometric to the figure consisting of
$x',y',z'$,
and
$w'$.

We have thus proven the implication
\begin{multline}\label{eq2.31}
(w \in \partial U,\,\, ]y,w[\,\, \subset U,\,\, ]y,w[\,\, \cap\,\,
]x,z[\,\, \ne\varnothing) \Rightarrow (]y',w'[\,\, \subset V,
\\
|w' - y'| = |w - y|,\,\,\, |w' - x'| = |w - x|,\,\,\, |w' - z'| = |w - z|).
\end{multline}
Similarly, we can prove the implications
\begin{multline*}
(w \in \partial U,\,\, ]x,w[\,\, \subset U,\,\, ]x,w[\,\, \cap\,\,
]y,z[\,\, \ne\varnothing) \Rightarrow (]x',w'[\,\, \subset V,
\\
|w' - y'| = |w - y|,\,\,\, |w' - x'| = |w - x|,\,\,\, |w' - z'| = |w - z|),
\end{multline*}
\begin{multline*}
(w \in \partial U,\,\, ]z,w[\,\, \subset U,\,\, ]z,w[\,\, \cap\,\,
]x,y[\,\, \ne\varnothing) \Rightarrow (]z',w'[\,\, \subset V,
\\
|w' - y'| = |w - y|,\,\,\, |w' - x'| = |w - x|,\,\,\, |w' - z'| = |w - z|).
\end{multline*}
Now, consider the case when
\begin{equation}\label{eq2.32}
\exists r > 0 \quad B(y',r) \cap (\inter \Delta(x',y',z')) \cap
\partial V = \varnothing.
\end{equation}
Then the desired equality~(\ref{eq2.28}) is very easy to prove.
Indeed, if~(\ref{eq2.32}) is valid but~(\ref{eq2.28}) is violated
then there is a point
$v \in (\inter \Delta(x',y',z')) \cap \partial V$
such that
$\angle(x',y',v)$
and
$\angle(v,y',z')$
are boundary angles. Take the corresponding element
$\widetilde v \in \partial_H V$,
$p\widetilde v = v$,
and consider
$u = pf^{-1}(\widetilde v) \in \partial U$.
Using the same method as in the proof of~(\ref{eq2.31}), we can
show that the figure consisting of the four points
$x'$,
$y'$,
$z'$,
and
$v$
is Euclideanly isometric to the figure consisting of the points
$x$,
$y$,
$z$,
and
$u$.
But then
$u \in (\inter \Delta(x,y,z)) \cap \partial U$
which contradicts~(\ref{eq2.25}).

Thus, it remains to consider the case when~(\ref{eq2.32})
is false, i.e., when
\begin{equation}\label{eq2.33}
\forall r > 0 \quad B(y',r) \cap (\inter \Delta(x',y',z')) \cap
\partial V \ne\varnothing.
\end{equation}
For the same reasons, without loss of generality we can additionally
assume that
\begin{equation}\label{eq2.34}
\forall r > 0 \quad B(x',r) \cap (\inter \Delta(x',y',z')) \cap
\partial V \ne\varnothing,
\end{equation}
\begin{equation}\label{eq2.35}
\forall r > 0 \quad B(z',r) \cap (\inter \Delta(x',y',z')) \cap
\partial V \ne\varnothing.
\end{equation}
Suppose now that at least one of the following two assertions is
valid:
\begin{equation}\label{eq2.36}
\exists w_0 \in\,\, ]x,z[\,\,\, \forall w_1 \in\,\, ]w_0,z[\,\,\,
\exists w \in \partial U \quad w_1 \in\,\, ]y,w[\,\, \subset U
\end{equation}
or
\begin{equation}\label{eq2.37}
\exists w_0 \in\,\, ]x,z[\,\,\, \forall w_1 \in\,\, ]x,w_0[\,\,\,
\exists w \in \partial U \quad w_1 \in\,\, ]y,w[\,\, \subset U.
\end{equation}
It follows from~(\ref{eq2.31}) that
\begin{equation}\label{eq2.38}
\exists w'_0 \in\,\, ]x',z'[\,\,\, \forall w'_1 \in\,\, ]w'_0,z'[\,\,\,
\exists w' \in \partial V \quad
w'_1 \in\,\, ]y',w'[\,\, \subset V
\end{equation}
or
\begin{equation}\label{eq2.39}
\exists w'_0 \in\,\, ]x',z'[\,\,\, \forall w'_1 \in\,\, ]x',w'_0[\,\,\,
\exists w' \in \partial V \quad
w'_1 \in\,\, ]y',w'[\,\, \subset V.
\end{equation}
If~(\ref{eq2.38}) holds, then in particular
$\exists w'_0 \in\,\, ]x',z'[\,\,\, \forall w'_1 \in\,\, ]w'_0,z'[
\quad ]y',w'_1[\,\, \subset V$
which contradicts~(\ref{eq2.35}). If~(\ref{eq2.39}) holds
then we similarly obtain a contradiction with~(\ref{eq2.34}).

We have thus proven that~(\ref{eq2.34}) and~(\ref{eq2.35}) imply
that none of the assertions~(\ref{eq2.36}) and~(\ref{eq2.37}) is
valid, i.e.,
\begin{equation*}
\forall w_0 \in\,\, ]x,z[\,\,\, \exists w_1 \in\,\, ]w_0,z[ \quad
\{w \in \mathbb R^2 : w - y = \tau(w_1 -y),\tau > 0\} \subset U;
\end{equation*}
\begin{equation*}
\forall w_0 \in\,\, ]x,z[\,\,\, \exists w_1 \in\,\, ]x,w_0[ \quad
\{w \in \mathbb R^2 : w - y = \tau(w_1 -y),\tau > 0\} \subset U.
\end{equation*}
It is geometrically obvious from here,~(\ref{eq2.25}), and
property~(H-i) that for every point
$\widetilde u \in \partial_H U$
such that
$u = p\widetilde u \in \inter \angle(x,y,z)$
and for each
$H$-shortest
curve
$\gamma : [0,L] \to \mathbb R^2$
joining
$\widetilde y$
and
$\widetilde u$,
there exist
$w \in \angle(x,y,z) \cap \partial U \setminus \Delta(x,y,z)$
and
$\alpha_0 \in\,\, ]0,L]$
such that
$\gamma(\alpha_0) =w$
and
$\gamma(\tau) \in\,\, ]y,w[\,\, \subset U$
for
$\tau \in\,\, ]0,\alpha_0[$.
On the other hand, it is geometrically obvious that if an
$H$-shortest
curve joins the vertex
$f(\widetilde y)$
with a point of
$\inter \Delta(x'y'z')$
then the
$H$-shortest
curve cannot leave the triangle
$\Delta(x'y'z')$.
Therefore, from the above, Lemma~\ref{l2.6}, and~(\ref{eq2.31})
we obtain the implications
\begin{equation}\label{eq2.40}
\forall \widetilde v \in \partial_H V \quad (v = p\widetilde v \in
\inter \Delta(x',y',z') \Rightarrow u = pf^{-1}(\widetilde v) \not\in
\inter \angle(x,y,z).
\end{equation}
Thus, from~(\ref{eq2.34}) and~(\ref{eq2.35}) we arrive
at~(\ref{eq2.40}). Similarly, from~(\ref{eq2.33})
and~(\ref{eq2.35}) we can also obtain two more implications
\begin{equation}\label{eq2.41}
\forall \widetilde v \in \partial_H V \quad (v = p\widetilde v \in
\inter \Delta(x',y',z') \Rightarrow
u = pf^{-1}(\widetilde v) \not\in \inter \angle(z,x,y));
\end{equation}
\begin{equation}\label{eq2.42}
\forall \widetilde v \in \partial_H V \quad (v = p\widetilde v \in
\inter \Delta(x',y',z') \Rightarrow
u = pf^{-1}(\widetilde v) \not\in \inter \angle(y,z,x)).
\end{equation}
We can assume without loss of generality that
\begin{equation}\label{eq2.43}
|y' - z'|,|x' - y'| \le |x' - z'|.
\end{equation}
From~(\ref{eq2.26}),~(\ref{eq2.34}), and~(\ref{eq2.35}) we easily
find now that there exist sequences
\begin{equation}\label{eq2.44}
\inter \Delta(x',y',z') \ni v^1_\nu \to x', \quad
\inter \Delta(x',y',z') \ni v^2_\nu \to z'
\end{equation}
such that
$v^1_\nu,v^2_\nu \in \partial V$
and
$]v^1_\nu,v^2_\nu[\,\, \subset V$
(i.e.,
$]v^1_\nu,v^2_\nu[$
is a boundary interval for
$V$)
and the following hold:
$$
\rho_{\partial_H V,V}(\widetilde v^1_\nu,f(\widetilde x)) \to 0,\quad
\rho_{\partial_H V,V}(\widetilde v^2_\nu,f(\widetilde z)) \to 0.
$$
From~(\ref{eq2.43}) and~(\ref{eq2.44}) we obtain
$|v^1_\nu - v^2_\nu| < |x' -z'|$.
Put
$u^1_\nu = pf^{-1}(\widetilde v^1_\nu)$
and
$u^2_\nu = pf^{-1}(\widetilde v^2_\nu)$.
Then from Lemma~\ref{l2.3} and the above assumptions we find that
$u^1_\nu,u^2_\nu \in \partial U$,
$]u^1_\nu,u^2_\nu[\,\, \subset U$,
and
\begin{equation}\label{eq2.45}
u^1_\nu \to x, \quad u^2_\nu \to z, \quad
|u^1_\nu - u^2_\nu| = |v^1_\nu -v^2_\nu| < |x' -z'| = |x - z|.
\end{equation}
Moreover,~(\ref{eq2.40})-(\ref{eq2.42}) and the previous
computations imply that
\begin{equation}\label{eq2.46}
u^1_\nu, u^2_\nu \not\in (\inter \Delta(x,y,z)) \cup
(\inter \Delta(z,x,y))  \cup (\inter \Delta(y,z,x)) \cup
\Delta(x,y,z).
\end{equation}
But~(\ref{eq2.45}) and~(\ref{eq2.46}) contradict each other.
The contradiction completes the proof of Lemma~\ref{l2.7}.

\begin{lemma}\label{l2.8}
Let
$\angle(x,y,z)$
be a boundary angle for
$U$.
Then
$\angle(x',y',z')$
is a boundary angle for
V.
\end{lemma}

We can express the essence of this lemma as follows: The mapping
$f$
takes boundary angles into boundary angles.

\textit{Proof.}
It follows from Lemmas~\ref{l2.3} and~\ref{l2.4} that
$x',y',z' \in \partial V$,
$]x',y'[\,\, \cup\,\, ]y',z'[\,\, \subset V$,
$|x' - y'| = |x - y|$,
$|z' - y'| = |z - y|$,
and
$|x' - z'| \le |x - z|$.
These formulas and the fact that
$x,y$,
and
$z$
do not lie on one straight line imply that
$x',y'$,
and
$z'$
do not lie on one straight line, too.

It remains to prove that
\begin{equation}\label{eq2.47}
\exists r > 0 \quad B(y',r) \cap \angle(x',y',z') \setminus \{y'\}
\subset V.
\end{equation}

Making parallel translations, if necessary, we can assume without
loss of generality that
$y = y' = 0 \in \mathbb R^2$.
We proceed firstly as in the proof of Lemma~\ref{l2.4}: repeat
all arguments of the proof of Lemma~\ref{l2.4} up
to~(\ref{eq2.20}). Consider the above-constructed mapping
$g :\Gamma_1 \to \mathbb R^2$.
By construction,
\begin{equation}\label{eq2.48}
\forall w \in \Gamma_1, \quad ]0,g(w)[\,\, \subset V.
\end{equation}
Extend
$g$
to the whole set
$\Gamma$
as follows: Let
$w \in \Gamma \setminus \Gamma_1\,\,(= \Gamma_2)$.
It is geometrically obvious that
$\exists w_1,w_2 \in \Gamma_1 \quad w \in\,\, ]w_1,w_2[\,\, \subset U$.
It is also geometrically obvious that the points
$w_1$
and
$w_2$
from the previous formula are determined uniquely; moreover,
\begin{equation}\label{eq2.49}
]w_1,w_2[\,\, \subset U \cap \Gamma,
\end{equation}
\begin{equation}\label{eq2.50}
\Delta(0,w_1,w_2)\,\, {\rm is\, a\, boundary\, triangle\, for}\,\, U.
\end{equation}
In this case put
$g(w) =sg(w_1) + (1 - s)g(w_2)$,
where
$s \in\,\, ]0,1[$
is determined from the equality
$w = sw_1 + (1 - s)w_2$.

From the definition of
$g$
and~(\ref{eq2.18}) we immediately find that
\begin{equation}\label{eq2.51}
\forall w \in \Gamma, \quad |g(w)| = |w|.
\end{equation}
Define the curve
$\gamma_V : [0,\Phi] \to \mathbb R^2$
by the rule
$\gamma_V(\varphi) = g(\gamma(\varphi))$
for
$\varphi \in [0,\Phi]$,
where
$\gamma(\varphi)$
is defined in the
proof of Lemma~\ref{l2.4}. By construction,
$\gamma_V$
is continuous, and from~(\ref{eq2.19}) and~(\ref{eq2.51}) we
obtain
\begin{equation}\label{eq2.52}
x' = \gamma_V(0), \quad z' = \gamma_V(\Phi)
\end{equation}
and
$\forall \varphi \in [0,\Phi] \quad |\gamma_V(\varphi)| =
|\gamma(\varphi)|$.
Using~(\ref{eq2.12}) and the construction of
$\gamma(\varphi)$,
from the last identity we conclude that there exists
$r_1 >0$
such that
\begin{equation}\label{eq2.53}
\forall \varphi \in [0,\Phi], \quad
|\gamma_V(\varphi)| \ge r_1.
\end{equation}
By~(\ref{eq2.48}),
\begin{equation}\label{eq2.54}
\forall \varphi \in [0,\Phi], \quad (\gamma(\varphi) \in \Gamma_1
\Rightarrow\,\, ]0,\gamma_V(\varphi)[\,\,
\subset V).
\end{equation}
Now, by construction (in particular, see~(\ref{eq2.49})
and~(\ref{eq2.50})) and Lemma~\ref{l2.7},
\begin{equation}\label{eq2.55}
\forall \varphi \in [0,\Phi], \quad (\gamma(\varphi) \in \Gamma_2
\Rightarrow\,\, ]0,\gamma_V(\varphi)[\,\,
\subset V).
\end{equation}
Combining~(\ref{eq2.54}) and~(\ref{eq2.55}), we finally obtain
the desired inclusions
\begin{equation}\label{eq2.56}
\forall \varphi \in [0,\Phi], \quad ]0,\gamma_V(\varphi)[\,\,
\subset V.
\end{equation}
Now,~(\ref{eq2.56})),~(\ref{eq2.53})),~(\ref{eq2.52})), and
continuity of
$\gamma_V$
yield the desired relation~(\ref{eq2.47})). Lemma~\ref{l2.8} is
proven.

\begin{lemma}\label{l2.9}
Let
$\angle(x,y,z)$
be a boundary angle for
$U$.
Then
$\angle(x',y',z')$
is a boundary angle for
$V;$
moreover{\rm,}
$|x' -y'| = |x - y|,$
$|z' -y'| = |z - y|,$
and
$|x' -z'| = |x - z|$.
\end{lemma}

We can express the essence of this lemma as follows: The mapping
$f$
takes boundary angles into equal boundary angles.

\textit{Proof.}
We have to apply Lemmas~\ref{l2.4} and~\ref{l2.8} twice.

In view of~(\ref{eq2.14}), the technical assertion below
contains Lemma~\ref{l2.9} as a particular case:

\begin{lemma}\label{l2.10}
Let
$\angle(x,y,z)$
be a boundary angle for
$U$.
Then
\begin{equation}\label{eq2.57}
\Gamma^V_{\angle(x',y',z')} = \{w' : w \in
\Gamma^U_{\angle(x,y,z)}\};
\end{equation}
moreover{\rm,}
\begin{equation}\label{eq2.58}
|w'_1 - w'_2| = |w_1 - w_2|, \quad \forall w_1,w_2 \in
\Gamma^U_{\angle(x,y,z)} \cup \{y\}.
\end{equation}
\end{lemma}

\textit{Proof.}
By Lemma~\ref{l2.9},
$\angle(x',y',z')$
is a boundary angle for
$V$
which is Euclideanly isometric to the angle
$\angle(x,y,z)$.
It is geometrically obvious now that, for every
$w \in \Gamma^U_{\angle(x,y,z)} \setminus \{x,z\}$,
$\angle(x,y,w)$
and
$\angle(w,y,z)$
are boundary angles for
$U$.
Similarly, for every
$v \in \Gamma^V_{\angle(x',y',z')} \setminus \{x',z'\}$,
$\angle(x',y',v)$
and
$\angle(v,y',z')$
are boundary angles for
$V$.
Now, the desired assertions~(\ref{eq2.57}) and~(\ref{eq2.58})
are obtained from the last two propositions, Lemma~\ref{l2.9},
and the following fact:

($P_1$) \textit{The position of a point on a plane is determined
uniquely from its distances to three points not lying on one
straight line.}

We drop the remaining technical details which are plain.

\begin{lemma}\label{l2.11}
Suppose that
$x_i \in \partial U,\,\, i = 1,\dots,m,$
are such that
$\angle(x_i,x_{i + 1},x_{i + 2})$
constitute boundary angles for
$U$
for all
$i = 1,\dots,m - 2$.
Then
$|x'_i - x'_j| = |x_i - x_j|$
for all
$i,j = 1,\dots,m$.
\end{lemma}

\textit{Proof.}
For
$m = 3$
the assertion of Lemma~\ref{l2.11} coincides with the assertion
of Lemma~\ref{l2.9}. Now, it clearly suffices to prove
Lemma~\ref{l2.11} in the case
$m = 4$,
since for
$m > 4$
the assertion of Lemma~\ref{l2.11} is derived by induction with
use of the simple fact
($P_1$)
(see above). Thus, we assume below that
$m = 4$.

It follows from Lemma~\ref{l2.9} and the definition of boundary
angles that
\begin{equation}\label{eq2.59}
|x'_i - x'_j| = |x_i - x_j|, \quad  i,j = 1,2,3,
\end{equation}
\begin{equation}\label{eq2.60}
|x'_i - x'_j| = |x_i - x_j|, \quad  i,j = 2,3,4,
\end{equation}
\begin{equation}\label{eq2.61}
x'_1,x'_2,\,\, {\rm and}\,\, x'_3\,\, {\rm are\,\, not\,\,
collinear},
\end{equation}
\begin{equation}\label{eq2.62}
x'_2,x'_3,\,\, {\rm and}\,\, x'_4\,\, {\rm are\,\, not\,\,
collinear}.
\end{equation}

It suffices to prove the only equality
\begin{equation}\label{eq2.63}
|x'_1 - x'_4| = |x_1 - x_4|.
\end{equation}

Making isometric transformations, if necessary, and
using~(\ref{eq2.59}), we can assume without loss of generality
that
\begin{equation}\label{eq2.64}
x'_i = x_i, \quad i = 1,2,3.
\end{equation}
Then, by~(\ref{eq2.59})-(\ref{eq2.62}), the desired
equality~(\ref{eq2.63}) is equivalent to
$x'_4 = x_4$.
Assume that the latter is false. Then,
by~(\ref{eq2.59})-(\ref{eq2.62}),
\begin{equation}\label{eq2.65}
x'_4\,\, {\rm is\,\, the\,\, reflection\,\, of}\,\, x_4\,\,
{\rm in\,\, the\,\, straight\,\, line}\,\, x_2x_3.
\end{equation}
It follows from~(\ref{eq2.64}) and~(\ref{eq2.58}) that
\begin{equation}\label{eq2.66}
w' = w, \quad {\rm for\,\, all}\,\, w \in \Gamma^U_{\angle(x_1,x_2,x_3)}.
\end{equation}
From~(\ref{eq2.65}),~(\ref{eq2.64}), and~(\ref{eq2.58}) we find
that
\begin{equation}\label{eq2.67}
w'\,\, {\rm is\,\, the\,\, reflection\,\, of}\,\, w\,\, {\rm in\,\,
the\,\, straight\,\, line}\,\, x_2x_3\,\, {\rm for\,\, all}\,\, w \in
\Gamma^U_{\angle(x_2,x_3,x_4)}.
\end{equation}
Then, by interchanging the domains
$U$
and
$V$,
if necessary, we can assume without loss of generality that
$x_1$
and
$x_4$
lie on one side of the straight line
$x_2x_3$
and
$x'_1$
and
$x'_4$
lie on the different sides of the straight line
$x'_2x'_3$.
It is geometrically obvious now that in our situation we deal
with at least one of the following three possibilities:\footnote{
Assertion (iii) is valid if
$x'_2 \in \cl\bigl(\Gamma^V_{\angle(x'_2,x'_3,x'_4)} \setminus
\{x'_2\}\bigr)$
and
$x'_3 \in \cl\bigl(\Gamma^V_{\angle(x'_1,x'_2,x'_3)} \setminus
\{x'_3\}\bigr)$,
and if at least one of these memberships is violated then one of
the assertions (i) and (ii) is valid.}

(i) there is a boundary angle
$\angle(x_1,x_2,u)$
for
$U$,
where
$u \in \Gamma^U_{\angle(x_2,x_3,x_4)} \setminus \{x_2\}$;

(ii) there is a boundary angle
$\angle(x_4,x_3,u)$
for
$U$,
where
$u \in \Gamma^U_{\angle(x_1,x_2,x_3)} \setminus \{x_3\}$;

(iii) there is a boundary interval
$]v_1,v_2[$
for
$V$,
joining
$v_1 \in \Gamma^V_{\angle(x'_1,x'_2,x'_3)} \setminus \{x'_3\}$
and
$v_2 \in \Gamma^V_{\angle(x'_2,x'_3,x'_4)} \setminus \{x'_2\}$.

By Lemmas~\ref{l2.3} and~\ref{l2.10}, (iii) is equivalent to the
following:

(iv)
there is a boundary interval
$]u_1,u_2[$
(for
$U$)
joining
$u_1 \in \Gamma^U_{\angle(x_1,x_2,x_3)} \setminus \{x_3\}$
and
$u_2 \in \Gamma^U_{\angle(x_2,x_3,x_4)} \setminus \{x_2\}$.

Therefore, we can assume that at least one of the three
assertions is always valid: (i), (ii), or (iv). But, by
Lemmas~\ref{l2.3} and~\ref{l2.9}, each of these assertions
implies the following:

(v) there exist\footnote{The points
$u_1$
and
$u_2$
from (v) are determined as follows: if (iv) is valid then they
coincide with the points in (iv) with the same names; if (i) is valid
then
$u_1 = x_1$
and
$u_2 = u$,
where
$u$
is defined in (i); if (ii) is valid then
$u_1 = u$
and
$u_2 = x_4$,
where
$u$
is defined in (ii).} points
$u_1 \in \Gamma^U_{\angle(x_1,x_2,x_3)} \setminus \{x_3\}$
and
$u_2 \in \Gamma^U_{\angle(x_2,x_3,x_4)} \setminus \{x_2\}$
such that
\begin{equation}\label{eq2.68}
|u'_2 - u'_1| = |u_2 - u_1|.
\end{equation}
It is geometrically obvious that
\begin{equation}\label{eq2.69}
u_1\,\, {\rm does\,\, not\,\, lie\,\, on\,\, the\,\, straight\,\,
line}\,\, x_2x_3,
\end{equation}
\begin{equation}\label{eq2.70}
u_2\,\, {\rm does\,\, not\,\, lie\,\, on\,\, the\,\, straight\,\,
line}\,\, x_2x_3.
\end{equation}
From Lemma~\ref{l2.10} (see~(\ref{eq2.14}) and~(\ref{eq2.58}))
and the membership
$u_2 \in \Gamma^U_{\angle(x_2,x_3,x_4)}$
we also obtain the equalities
\begin{equation}\label{eq2.71}
|u'_2 - x'_2| = |u_2 - x_2|,
\end{equation}
\begin{equation}\label{eq2.72}
|u'_2 - x'_3| = |u_2 - x_3|.
\end{equation}
Using~(\ref{eq2.64}) and~(\ref{eq2.66}), we can
rewrite~(\ref{eq2.68}),~(\ref{eq2.71}), and~(\ref{eq2.72}) as
$|u'_2 - u_1| = |u_2 - u_1|$,
$|u'_2 - x_2| = |u_2 - x_2|$,
and
$|u'_2 - x_3| = |u_2 - x_3|$.
Then from~(\ref{eq2.69}) and property
($P_1$)
we find that
$u'_2 = u_2$.
But the last equality contradicts~(\ref{eq2.67})
and~(\ref{eq2.70}). The contradiction completes the proof of
Lemma~\ref{l2.11}.

\begin{lemma}\label{l2.12}
Suppose that
$x_1,x_2,x_3 \in \partial U$
and
$z \in\,\, ]x_1,x_2[$
are such that
$]x_1,x_2[$
is a boundary interval and
$]x_3,z] \subset U$,
$z = \tau x_1 + (1 - \tau)x_2 \in\,\, ]x_1,x_2[$,
$\tau \in\,\, ]0,1[$.
Put
$z' = \tau x'_1 + (1 - \tau)x'_2$
and
$x'_3 = pf(\widetilde x_3)$,
where
$\widetilde x_3$
denotes the element of the Hausdorff boundary
$\partial_H U$
generated by a Cauchy sequence of points of the interval
$]x_3,z[$
converging to
$x_3$.
Then
$|x'_i - x'_j| = |x_i - x_j|$
for all
$i,j = 1,\dots,3$,
$]x'_3,z'] \subset V$,
and
$f(\widetilde x_3)$
coincides with the element of the Hausdorff boundary
$\partial_H V$
generated by a Cauchy sequence of points of the interval
$]x'_3,z'[$
converging to
$x'_3$.
\end{lemma}

\textit{Proof.}
If
$x_1$,
$x_2$,
and
$x_3$
are collinear, then
$x_3$
coincides with
$x_1$
or
$x_2$
and there  is nothing to prove. Henceforth we assume that
$x_1$,
$x_2$,
and
$x_3$
are not collinear. Put
$\Gamma_1 = \Gamma^U_{\angle(x_1,z,x_3)}$
and
$\Gamma_2 = \Gamma^U_{\angle(x_2,z,x_3)}$,
where these sets are defined as in the case of boundary angles
(see~(\ref{eq2.13})). It is geometrically obvious that there is
a finite or infinite sequence of points
$u_\nu$,
$\nu \in \mathbb Z \cap\,\, ]\alpha,\omega[$,
$-\infty \le \alpha < \omega \le +\infty$,
such that
$u_\nu \in \Gamma_1$
for odd
$\nu$
and
$u_\nu \in \Gamma_2$
for even
$\nu$;
moreover,
$$
\angle(u_\nu,u_{\nu + 1},u_{\nu + 2})\,\, {\rm is\,\, a\,\,
boundary\,\, angle\,\,for}\,\, \nu \in \mathbb Z \cap\,\,
]\alpha,\omega - 2[,
$$
$$
{\rm if}\,\, \omega = +\infty,\,\, {\rm then}\,\, u_\nu \to x_3\,\,
{\rm as}\,\, \nu \to +\infty,
$$
$$
{\rm if}\,\, \omega < +\infty,\,\, {\rm then}\,\, u_{\omega - 1}
= x_3,
$$
$$
{\rm if}\,\, \alpha = -\infty,\,\, {\rm then}\,\, u_{2k + 1} \to
x_1\,\, {\rm and}\,\, u_{2k} \to x_2\,\, {\rm as}\,\, k \to -\infty,
$$
$$
{\rm if}\,\, \alpha > -\infty,\,\, {\rm then}\,\,
\{u_{\alpha + 1},u_{\alpha + 2}\} = \{x_1,x_2\}.
$$
Now, the assertion of Lemma~\ref{l2.12} easily follows from
Lemmas~\ref{l2.9}-\ref{l2.11} and obvious geometric arguments.
We skip the details of computations in view of their simplicity.

\textbf{Remark~2.2.}
Actually, Lemma~\ref{l2.10} and the method of the proof of
Lemma~\ref{l2.12} yield the following stronger assertion:
$$
\Gamma^V_{\angle(x'_1,z',x'_3)} = \bigl\{w' : w \in
\Gamma^U_{\angle(x_1,z,x_3)}\bigr\}, \quad
\Gamma^V_{\angle(x'_2,z',x'_3)} = \bigl\{w' : w \in
\Gamma^U_{\angle(x_2,z,x_3)}\bigr\};
$$
moreover,
$$
|w'_1 - w'_2| = |w_1 - w_2| \quad {\rm for\,\, all}\,\, w_1,w_2 \in
\Gamma^U_{\angle(x_1,z,x_3)} \cup \Gamma^U_{\angle(x_2,z,x_3)}
\cup \{z\}.
$$

\textit{Proof of Item {\rm(I)} of Theorem~{\rm\ref{t1.1}}.}
Necessity easily follows from Lemmas~\ref{l2.12} and~\ref{l2.3},
while sufficiency is obvious.

\textbf{Remark~2.3.}
In view of Item (I) of Theorem~\ref{t1.1}, we will suppose unless
the contrary is specified that there is no straight line
containing
$\partial U$
or
$\partial V$.

\begin{lemma}\label{l2.13}
Suppose that
$x_1,x_2,y_1,y_2 \in \partial U$
and
$z_1,z_2 \in U$
are such that
$]x_1,y_1[$
and
$]x_2,y_2[$
are boundary intervals{\rm,}
$x_1 = \tau_1 x_1 + (1 - \tau_1)y_1 \in\,\, ]x_1,y_1[,$
and
$x_2 = \tau_2 x_2 + (1 - \tau_2)y_2 \in\,\, ]x_2,y_2[,$
$\tau_1,\tau_2 \in\,\, ]0,1[;$
moreover{\rm,}
$[x_1,x_2] \in U$.
Put
$x'_1 : = \tau_1 x'_1 + (1 - \tau_1)y'_1$
and
$x'_2 : = \tau_2 x'_2 + (1 - \tau_2)y'_2$.
Then
$[z'_1,z'_2] \in V$
and
the quadrangle
$x_1y_1y_2x_2$
is Euclideanly isometric to the quadrangle
$x'_1y'_1y'_2x'_2$.
\end{lemma}

\textit{Proof.}
If
$]x_1,y_1[\,\, \cap\,\, ]x_2,y_2[\,\, \ne\varnothing$,
then the assertion of Lemma~\ref{l2.13} is easy from
Lemma~\ref{l2.12} and Remark 2.2. Now, consider the situation when
$]x_1,y_1[\,\, \cap\,\, ]x_2,y_2[\,\, = \varnothing$.
We assume without loss of generality that
$x_1$
and
$x_2$
lie on one side of the straight line
$x_1x_2$.
Consider the quadrangle
$G_x$
with vertices
$z_1x_1x_2z_2$.
Put
$E_x = \{w :\,\, {\rm there\,\, is\,\, a\,\, collection\,\,
of\,\, points}\,\, w_i \in \partial U,\, i = 1,\dots, k,\,\,
{\rm such\,\, that}\,\, w \in \conv(w_1,\dots,w_k) \subset G_x\}$.
By construction,
$[x_1,x_2] \subset E_x$.
Put
$$
\Gamma^U_x(x_1,y_1,x_2,y_2,z_1,z_2) = \{w \in (\partial U) \cap
(\partial E_x) : \exists z \in [z_1,z_2]\,\,\, ]z,w[\,\, \subset U
\setminus E_x\}.
$$
By construction, we find also that
$x_1,x_2 \in \Gamma^U_x(x_1,y_1,x_2,y_2,z_1,z_2)$.
Similarly, consider the quadrangle
$G_y$
with vertices
$z_1y_1y_2z_2$.
Put
\begin{multline*}
E_y = \{w :\text{ there is a collection
of points } w_i \in \partial U,\, i = 1,\dots, k,
\\
\text{ such that } 
w \in \conv(w_1,\dots,w_k) \subset G_y\}.
\end{multline*}
Denote
$$
\Gamma^U_y(x_1,y_1,x_2,y_2,z_1,z_2) = \{w \in (\partial U) \cap
(\partial E_y) : \exists z \in [z_1,z_2]\}\,\,\, ]z,w[\,\,
\subset U \setminus E_y\}.
$$
By construction,
$y_1,y_2 \in \Gamma^U_y(x_1,y_1,x_2,y_2,z_1,z_2)$.
It is geometrically obvious that there is a finite or infinite
sequence of points
$u_\nu$,
$\nu \in \mathbb Z \cap\,\,]\alpha,\omega[$,
$-\infty \le \alpha < \omega \le +\infty$,
such that
$u_\nu \in \Gamma^U_x(x_1,y_1,x_2,y_2,z_1,z_2)$
for odd
$\nu$
and
$u_\nu \in \Gamma^U_y(x_1,y_1,x_2,y_2,z_1,z_2)$
for even
$\nu$;
moreover,
$$
\angle(u_\nu,u_{\nu + 1},u_{\nu + 2})\,\, {\rm is\,\, a\,\,
boundary\,\, angle\,\,for}\,\, \nu \in \mathbb Z\,\, \cap\,\,
]\alpha,\omega - 2[,
$$
$$
{\rm if}\,\, \omega = +\infty,\,\, {\rm then}\,\, u_{2k + 1} \to
x_2\,\, {\rm and}\,\, u_{2k} \to y_2\,\, {\rm as}\,\,
k \to +\infty,
$$
$$
{\rm if}\,\, \omega < +\infty,\,\, {\rm then}\,\,
\{u_{\omega - 1},u_{\omega - 2}\} = \{x_1,y_1\},
$$
$$
{\rm if}\,\, \alpha = -\infty,\,\, {\rm then}\,\, u_{2k + 1} \to
x_1\,\, {\rm and}\,\, u_{2k} \to y_1\,\, {\rm as}\,\,
k \to -\infty,
$$
$$
{\rm if}\,\, \alpha > -\infty,\,\, {\rm then}\,\,
\{u_{\alpha + 1},u_{\alpha + 2}\} = \{x_2,y_2\}.
$$
Now, the assertion of Lemma~\ref{l2.13} is easy from
Lemmas~\ref{l2.9}-\ref{l2.11} and obvious geometric arguments.
We skip the details of computations in view of their simplicity.

\textbf{Remark~2.4.}
From Lemmas~\ref{l2.10} and~\ref{l2.11} and the method of the
proof of Lemma~\ref{l2.13} we derive the equalities
$$
\Gamma^V_x(x'_1,y'_1,x'_2,y'_2,z'_1,z'_2) = \{w' : w \in
\Gamma^U_x(x_1,y_1,x_2,y_2,z_1,z_2)\},
$$
$$
\Gamma^V_y(x'_1,y'_1,x'_2,y'_2,z'_1,z'_2) = \{w' : w \in
\Gamma^U_y(x_1,y_1,x_2,y_2,z_1,z_2)\},
$$
$$
\forall w_1,w_2 \in \Gamma^U_x(x_1,y_1,x_2,y_2,z_1,z_2) \cup
\Gamma^U_y(x_1,y_1,x_2,y_2,z_1,z_2)\} \quad |w'_1 -w'_2| =
|w_1 - w_2|.
$$

We say that boundary intervals
$]x_1,y_1[$
and
$]x_2,y_2[$
(with respect to
$U$)
are
$U$-joinable
if there exist points
$x_1,x_2 \in U$
such that the conditions of Lemma~\ref{l2.13} are satisfied.

\begin{lemma}\label{l2.14}
Let
$]x_i,y_i[,$
$i = 1,\dots,m,$
be boundary intervals {\rm(}with respect to
$U)$.
Suppose that
the intervals
$]x_i,y_i[$
and
$]x_{i + 1},y_{i + 1}[$
are
$U$-joinable
for all
$i = 1,\dots, m - 1$.
Then the boundary intervals
$]x'_i,y'_i[$
and
$]x'_{i + 1},y'_{i + 1}[$
are
$V$-joinable
for all
$i = 1,\dots,m - 1;$
moreover{\rm,}
\begin{equation}\label{eq2.73}
|x'_i - x'_j| = |x_i - x_j|,\,\,\, |x'_i - y'_j| = |x_i - y_j|,\,\,\,
|y'_i - y'_j| = |y_i - y_j|,\,\,\, i,j = 1,\dots,m.
\end{equation}
\end{lemma}

\textit{Proof.}
For
$m = 2$
the assertion of Lemma~\ref{l2.14} coincides with that of
Lemma~\ref{l2.13}. Now, it suffices clearly to prove
Lemma~\ref{l2.14} in the case
$m = 3$,
since for
$m > 3$
the assertion of Lemma~\ref{l2.14} is derived by induction on
using the simple fact
($P_1$)
(see above). Thus, we assume below that
$m = 3$.

By the definition of
$U$-joinable
intervals, there exist
$z^+_1,z^-_2,z^+_2,z^-_3 \in U$
such that
$z^+_i = \tau^+_i x_i + (1 - \tau^+_i)y_i \in\,\, ]x_i,y_i[$,
$z^-_i = \tau^-_i x_i + (1 - \tau^-_i)y_i \in\,\, ]x_i,y_i[$,
and
$[z^+_i,z^-_{i + 1}] \subset U$.
Denote
$z'^+_i = \tau^+_i x'_i + (1 - \tau^+_i)y'_i$
and
$z'^-_i = \tau^-_i x'_i + (1 - \tau^-_i)y'_i$.
Then Lemmas~\ref{l2.13} and~\ref{l2.3} immediately imply that
$[z'^+_i,z'^-_{i + 1}] \subset V$,
$i = 1,2$,
\begin{equation}\label{eq2.74}
|x'_i - x'_j| = |x_i - x_j|,\,\,\, |x'_i - y'_j| = |x_i - y_j|,\,\,\,
|y'_i - y'_j| = |y_i - y_j|,\,\,\, i,j = 1,2,
\end{equation}
\begin{equation}\label{eq2.75}
|x'_i - x'_j| = |x_i - x_j|,\,\,\, |x'_i - y'_j| = |x_i - y_j|,\,\,\,
|y'_i - y'_j| = |y_i - y_j|,\,\,\, i,j = 2,3.
\end{equation}

Making isometric transformations, if necessary, and
using~(\ref{eq2.74}), we can assume without loss of generality
that
\begin{equation}\label{eq2.76}
x'_i = x_i ,\,\,\, y'_i = y_i,\,\,\, i = 1,2.
\end{equation}
Then, by~(\ref{eq2.74}) and~(\ref{eq2.75}), the desired
equality~(\ref{eq2.73}) is equivalent to
$x'_3 = x_3$
and
$y'_3 = y_3$.
Assume that the latter are false. Then, by~(\ref{eq2.75})
and~(\ref{eq2.76}),
\begin{equation}\label{eq2.77}
x'_3\,\, {\rm is\,\, the\,\, reflection\,\, of}\,\, x_3\,\,
{\rm in\,\, the\,\, straight\,\, line}\,\, x_2y_2,
\end{equation}
\begin{equation}\label{eq2.78}
y'_3\,\, {\rm is\,\, the\,\, reflection\,\, of}\,\, y_3\,\,
{\rm in\,\, the\,\, straight\,\, line}\,\, x_2y_2.
\end{equation}
Suppose first that
\begin{equation}\label{eq2.79}
]x_1,y_1[\,\, \cap\,\, ]x_2,y_2[\,\, =\,\,
]x_2,y_2[\,\, \cap\,\, ]x_3,y_3[\,\, = \varnothing.
\end{equation}
We can assume without loss of generality that
$x_1$
and
$z^+_1$
lie on one side of the straight line
$x_2x_3$,
$x_1$
and
$x_2$
lie on one side of the straight line
$z^+_1z^-_2$,
and
$x_2$
and
$x_3$
lie on one side of the straight line
$z^+_2z^-_3$.
It follows from Remark~2.4,~(\ref{eq2.77}), and~(\ref{eq2.78})
that
\begin{equation}\label{eq2.80}
w' = w\,\, {\rm for\,\, all}\,\, w \in
\Gamma^U_x(x_1,y_1,x_2,y_2,z^+_1,z^-_2) \cup
\Gamma^U_y(x_1,y_1,x_2,y_2,z^+_1,z^-_2),
\end{equation}
\begin{multline}\label{eq2.81}
w'\,\, {\rm is\,\, the\,\, reflection\,\, of}\,\, w\,\,
{\rm in\,\, the\,\ straight\,\, line}\,\, x_2y_2\,\, {\rm for\,\, all}
\\
w \in \Gamma^U_x(x_3,y_3,x_2,y_2,z^-_3,z^+_2) \cup
\Gamma^U_y(x_3,y_3,x_2,y_2,z^-_3,z^+_2).
\end{multline}
Interchanging the domains
$U$
and
$V$,
if necessary, we can assume without loss of generality that
$x_1$,
$z^+_1$,
and
$z^-_3$
lie on one side of the straight line
$x_2y_2$.
It is geometrically obvious now that in our situation we always
deal with at least one of the following three possibilities:

(i) there is a point
$u_0 \in \Gamma^U_x(x_1,y_1,x_2,y_2,z^+_1,z^-_2) \cap
\Gamma^U_x(x_3,y_3,x_2,y_2,z^-_3,z^+_2) \setminus \{x_2\}$;

(ii) there is a point
$u_0 \in \Gamma^U_x(x_1,y_1,x_2,y_2,z^+_1,z^-_2)$
such that
$]x_2,u_0[$
is a boundary interval intersecting
$]z^+_2,z^-_3[\,\, \cup\,\, ]x_3,y_3[$
and, thereby, the interval
$]x_2,u_0[$
is
$U$-joinable
with the interval
$]x_3,y_3[$;

(iii) there is a point
$u_0 \in \Gamma^U_x(x_3,y_3,x_2,y_2,z^-_3,z^+_2)$
such that
$]x_2,u_0[$
is a boundary interval intersecting
$]z^+_1,z^-_2[\,\, \cup\,\, ]x_1,y_1[$
and, thereby, the interval
$]x_2,u_0[$
is
$U$-joinable
with the interval
$]x_1,y_1[$.

But, by Lemma~\ref{l2.13} and Remark~2.4, each of these
assertions implies the following:

(iv) there is a point
$u_0 \in \Gamma^U_x(x_1,y_1,x_2,y_2,z^+_1,z^-_2) \cup
\Gamma^U_x(x_3,y_3,x_2,y_2,z^-_3,z^+_2)$
such that
$u_0$
does not lie on the straight line
$x_2y_2$
and
$|u'_0 - x'_i| = |u_0 - x_i|$
and
$|u'_0 - y'_i| = |u_0 - y_i|$
for all
$i = 1,2,3$.

However, this obviously
contradicts~(\ref{eq2.77}),~(\ref{eq2.78}),~(\ref{eq2.80}),
and~(\ref{eq2.81}). This contradiction completes the proof of
Lemma~\ref{l2.14} in the case when~(\ref{eq2.79}) are valid. If
equalities~(\ref{eq2.79}) are false then the proof is carried
out in a similar way; moreover, it is even simpler: we only need
to use Lemma~\ref{l2.12} and Remark~2.2. We skip the corresponding
computation in view of their simplicity.

\textbf{Definition~2.1.}
We say that boundary intervals
$]x,y[$
and
$]w,z[$
(with respect to
$U$)
are
$U$-equivalent
if there is a sequence of boundary intervals
$]x_i,y_i[$,
$i = 1,\dots,m$,
such that
$x_1 = x$,
$y_1 = y$,
$x_m = w$,
$y_m = z$,
and
the intervals
$]x_i,y_i[$
and
$]x_{i + 1},y_{i + 1}[$
are
$U$-joinable
for all
$i = 1,\dots,m - 1$.

Now, we can restate Lemma~\ref{l2.14} as follows:

\begin{lemma}\label{l2.15}
Boundary intervals
$]x,y[$
and
$]w,z[$
{\rm(}with respect to
$U)$
are
$U$-equivalent
if and only if the intervals
$]x',y'[$
and
$]w',z'[$
are
$V$-equivalent.
If boundary intervals
$]x,y[$
and
$]w,z[$
are
$U$-equivalent
then the quadrangle
$xywz$
is Euclideanly isometric to the
quadrangle
$x'y'w'z'$.
\end{lemma}

\begin{lemma}\label{l2.16}
Boundary intervals
$]x,y[$
and
$]w,z[$
{\rm(}with respect to
$U)$
are
$U$-equivalent
if and only if there is a component
$U_i$
{\rm(}see the definition of
$U_i$
in Section~{\rm 1)} such that
$]x,y[\,\, \cup\,\, ]w,z[\,\, \subset \cl U_i$.
\end{lemma}

\textit{Proof} is simple and carried out by elementary means;
therefore, we omit it.

\begin{lemma}\label{l2.17}
Let
$\widetilde x \in (\partial_H U_i) \cap (\partial_H U)$
and
$x = p\widetilde x$.
Then there are sequences
$\{x_\nu\}_{\nu \in \mathbb N},\{y_\nu\}_{\nu \in \mathbb N}$
of points such that
$]x_\nu,y_\nu[$
are boundary intervals for
$U$
and
$]x_\nu,y_\nu[\,\, \subset \cl U_i;$
moreover{\rm,}
$\rho_{\partial_H U,U}(\widetilde x_\nu,\widetilde x) \le
\rho_{\partial_H U_i,U_i}(\widetilde x_\nu,\widetilde x) \to 0$
as
$\nu \to \infty$.
\end{lemma}

In connection with Lemma~\ref{l2.17}, observe that by
$(\partial_H U_i) \cap (\partial_H U)$
we naturally mean the set of elements
$\widetilde x \in \partial_H U$
such that there is a sequence of points
$w_\nu \in U_i$
which is Cauchy with respect to the intrinsic metric of
$U_i$
and such that
$\rho_{U_H}(w_\nu,\widetilde x) \to 0$.

\textit{Proof of Lemma~{\rm\ref{l2.17}}.}
Take a sequence
$w_\nu$.
Then from the corresponding definitions we immediately see that
$\rho_{U_H}(w_\nu,\widetilde x) \le \rho_{U_{i_H}}(w_\nu,
\widetilde x) \to 0$
as
$\nu \to \infty$.

Put
$R_\nu = \dist(w_\nu,\partial U)$.
By the choice of
$w_\nu$
and
$R_\nu$,
we find that
$R_\nu \to 0$
as
$\nu \to \infty$
and
\begin{equation}\label{eq2.82}
\forall \nu \quad (B(w_\nu,R_\nu) \subset U\,\, \&\,\, \exists
x_\nu \in \cl B(w_\nu,R_\nu) \cap \partial U).
\end{equation}

By construction,
$]x_\nu,w_\nu[\,\, \subset U_i$,
$x_\nu \in \partial U_i$.
If we take a sequence of points of the (open) interval
$]x_\nu,w_\nu[$
converging to
$x_\nu$
then it generates some element of the Hausdorff boundary
$\widetilde x_\nu \in (\partial_H U_i) \cap (\partial_H U)$.
Hence, we find in turn that
$\rho_{U_H}(\widetilde x_\nu,w_\nu) =
\rho_{U_{i_H}}(\widetilde x_\nu,w_\nu) = R_\nu \to 0$
as
$\nu \to \infty$.
It is geometrically obvious now that
\begin{equation}\label{eq2.83}
\forall \nu\,\,\,  \exists y_\nu \in \partial U \quad
(]x_\nu,y_\nu[\,\, \subset U\,\, \&\,\, ]x_\nu,y_\nu[\,\, \cap
B(w_\nu,R_\nu) \ne\varnothing).
\end{equation}
Indeed, if this were false then, by~(\ref{eq2.82}), the
half-plane passing through the point
$x_\nu$
and containing the ball
$B(w_\nu,R_\nu)$
would lie entirely in
$U$.
But this contradicts the membership
$w_\nu \in U_i$
and the definition of
$U_i$.

Thus,~(\ref{eq2.83}) is valid. Then from the definition of
$U_i$
and the previous constructions we immediately find that
$]x_\nu,y_\nu[\,\, \subset \cl U_i$.
The last inclusion together with the facts established above
yields the claim of Lemma~\ref{l2.17}.

\begin{lemma}\label{l2.18}
We can enumerate the components
$U_i$
and
$V_i$
so that this enumeration determines a bijective correspondence
between the components
$U_i$
and
$V_i;$
moreover{\rm,} there exist Euclideanly isometric mappings
$Q_i : \mathbb R^2 \to \mathbb R^2$
such that
$Q_i(U_i) = V_i$
for each component
$U_i$.
Now{\rm,}
$\widetilde x \in (\partial_H U_i) \cap (\partial_H U)$
if and only if
$f(\widetilde x) \in (\partial_H V_i) \cap (\partial_H V);$
moreover{\rm,} once these containments hold{\rm,} we also have
$x' = Q_i(x),$
where
$x = p\widetilde x$.
Finally{\rm,}
$y \in U \cap \partial U_i$
if and only if
$Q_i(y) \in V \cap \partial V_i$.
\end{lemma}

Lemma~\ref{l2.18} is a simple consequence of Lemmas~\ref{l2.1}
and~\ref{l2.15}-\ref{l2.17}; details are omitted. Henceforth we
assume that the components
$U_i$
and
$V_i$
are enumerated as in Lemma~\ref{l2.18}.

From Lemma~\ref{l2.18} we immediately obtain the following:

\begin{lemma}\label{l2.19}
If
$F_U = \varnothing$
{\rm(}see the definition of
$F_U$
in Section~{\rm 1)} then
$U$
is determined uniquely from the relative metric of its Hausdorff
boundary.
\end{lemma}

\textbf{Remark~2.5.}
Now, by Lemma~\ref{l2.19}, we suppose unless the contrary is
specified that
$F_U \ne\varnothing$.

\textbf{Remark~2.6.}
Let
$x \in (\partial F_U) \cap (\partial U)$.
Then, obviously, there is a unique\footnote{Here we use
Remark~2.3.} element
$\widetilde x \in (\partial_H F_U) \cap (\partial_H U)$
such that
$p\widetilde x = x$.
In this case we put
$x' =pf(\widetilde x)$
unless the contrary is specified.

From Lemma~\ref{l2.18} and Remark~2.6 we also easily obtain the
following:

\begin{lemma}\label{l2.20}
If
$x \in (\partial F_U) \cap (\partial U)$
then
$f(\widetilde x) \in (\partial_H F_V) \cap (\partial_H V)$
and{\rm,} in particular{\rm,}
$x' \in (\partial F_U) \cap (\partial U)$.
If
$x \in (\partial F_U) \setminus (\partial U)$
then there is a unique number
$i = i(x)$
such that
$x \in U \cap \partial U_i$
and
$Q_i(x) \in V \cap \partial V_i$.
\end{lemma}

The following geometric lemma is trivial:

\begin{lemma}\label{l2.21}
The boundary
$\partial F_U$
has at most two connected components. Moreover{\rm,} if the
boundary
$\partial _U$
has two connected components then these connected components are
two parallel straight lines and the set
$F_U$
itself represents the union of two half-planes by these straight
lines.
\end{lemma}

Lemma~\ref{l2.21} will help us to derive the two lemmas:

\begin{lemma}\label{l2.22}
The boundary
$\partial F_U$ contains two connected components if and only if
the boundary
$\partial F_V$
contains two connected components.
\end{lemma}

\textit{Proof} bases on using Lemmas~\ref{l2.18},~\ref{l2.20},
and~\ref{l2.21} together with the general facts on the
continuity of
$f$.
We skip the details of computations in view of their simplicity.

\begin{lemma}\label{l2.23}
If the boundary
$\partial F_U$
is disconnected then
$U$
is determined uniquely from the relative metric of its Hausdorff
boundary.
\end{lemma}

\textit{Proof} is carried out by simple application of
Lemmas~\ref{l2.18} and~\ref{l2.20}-\ref{l2.22}.

\begin{lemma}\label{l2.24}
Suppose that each of the boundaries
$\partial F_U$
and
$\partial F_V$
is a connected set {\rm(}i.e.{\rm,} a convex curve{\rm)}. Define
the mapping
$\theta : \partial F_U \to \mathbb R^2$
as
$$
\theta(w) =
\begin{cases}
w', & w \in (\partial F_U) \cap (\partial U);\\
Q_i(w), & w \in U \cap \partial U_i.
\end{cases}
$$
Then
$\theta$
is a homeomorphism between the convex curves
$\partial F_U$
and
$\partial F_V;$
moreover{\rm,} this homeomorphism preserves the arc length
{\rm(}i.e.{\rm,} for arbitrary two points
$x$
and
$y$
the length of the arc connecting them in
$\partial F_U$
coincides with the length of the image of this arc under the mapping
$\theta)$.
\end{lemma}

\textit{Proof.}
The fact that
$\theta$
is a homeomorphism of the convex curves
$\partial F_U$
and
$\partial F_V$
follows from construction and Lemmas~\ref{l2.18} and~\ref{l2.20}
(also see Remarks~2.6 and~2.1). Now, prove that this
homeomorphism preserves arc length. Take an arbitrary pair of
points
$x,y \in \partial F_U$.
Since all
$Q_i$
are isometries, by the construction of
$\theta$,
we can assume without loss of generality that
$x,y \in \partial U$.
Denote by
$\Gamma$
the convex arc joining
$x$
and
$y$
in
$\partial F_U$
and put
$\Gamma_1 = \Gamma \cap \partial U$.
Let
$\gamma : [0,\Phi] \to \Gamma$
be some parametrization of the arc
$\Gamma$.

Since
$\Gamma$
is a convex arc whose all extreme points lie in
$\Gamma_1$,
it is geometrically obvious that, for each
$\varepsilon > 0$,
there is a partition
$0 < \varphi_0 < \varphi_1 < \dots < \varphi_N = \Phi$,
$N = N(\varepsilon)$,
such that
$\forall \nu = 0,1,\dots,N \quad \gamma(\varphi_\nu) \in \Gamma_1$
and
\begin{equation}\label{eq2.84}
\forall \nu = 0,1,\dots,N - 1 \quad
\frac{l(\gamma([\varphi_\nu,\varphi_{\nu + 1}]))}
{|\gamma(\varphi_\nu) - \gamma(\varphi_{\nu + 1})|} \le 1 + \varepsilon.
\end{equation}
Using the definition of Hausdorff metric and the definition of
$F_U$,
we find that
\begin{equation}\label{eq2.85}
\forall \nu = 0,1,\dots,N - 1 \quad
\rho_{V_H}(f(\widetilde {\gamma(\varphi_\nu)}),
f(\widetilde {\gamma(\varphi_{\nu + 1})})) \le
l(\theta(\gamma([\varphi_\nu,\varphi_{\nu + 1}]))).
\end{equation}
From~(\ref{eq2.84}) and~(\ref{eq2.85}) we conclude that
\begin{multline}\label{eq2.86}
\forall \nu = 0,1,\dots,N - 1 \quad
\frac{1}{1 + \varepsilon}l(\gamma([\varphi_\nu,\varphi_{\nu + 1}])) \le
|\gamma(\varphi_\nu) - \gamma(\varphi_{\nu + 1})| \le
\\
\rho_{U_H}(\widetilde {\gamma(\varphi_\nu)},
\widetilde {\gamma(\varphi_{\nu + 1})}) =
\rho_{V_H}(f(\widetilde {\gamma(\varphi_\nu)}),
f(\widetilde {\gamma(\varphi_{\nu + 1})})) \le
l(\theta(\gamma([\varphi_\nu,\varphi_{\nu + 1}]))).
\end{multline}
It follows from the arbitrariness of
$\varepsilon > 0$
and~(\ref{eq2.86}) that
$l(\Gamma) \le l(\theta(\Gamma))$.
The reverse inequality follows, if we interchange the domains
$U$
and
$V$.
Lemma~\ref{l2.24} is proven.

\textit{Proof of Item {\rm(II)} of Theorem~{\rm\ref{t1.1}}.}
Necessity in Item (II) follows from Item (I) proven above and
Lemmas~\ref{l2.3},~\ref{l2.18},~\ref{l2.23}, and~\ref{l2.24}.
Sufficiency in Item (II) is obvious: in the presence of the
homeomorphism
$\theta : \partial F_U \to \partial F_V$
and the isometries
$Q_i : \mathbb R^2 \to \mathbb R^2$
with properties (IIa)-(IIc) we can naturally construct the
isometry
$f : \partial_H U \to \partial_H V$.

\textit{Proof of Theorem~{\rm\ref{t1.2}}.}
Let
$U$
be a convex domain in
$\mathbb R^n$,
$n \ge 2$,
different from a half-space. Under these conditions, we can
identify the elements of the Hausdorff boundary
$\partial_H U$
with the elements of the usual Euclidean boundary
$\partial U$.
Suppose that the domain
$V \subset \mathbb R^n$
is such that there is an isometry
$f : \partial_H U \to \partial_H V$.
Then the isometry
$f$
naturally generates the mapping
$g : \partial U \to \partial V$
by the formula
$g(x) = pf(x)$.
It is easily seen that, by convexity of
$U$,
the resulting mapping
$g$
is
$1$-Lipschitz.
We have the following elementary lemma:

\begin{lemma}\label{l2.25}
Let
$U$
be a convex domain in
$\mathbb R^n,$
$n \ge 2,$
different from a half-space and take
$[x,y] \subset U$.
Then there is a sequence of boundary intervals
$]x_\nu,y_\nu[$
{\rm(}with respect to
$U)$
such that
$\forall z \in [x,y] \quad \dist(z,[x_\nu,y_\nu)\,\,
(= \inf_{w \in [x_\nu,y_\nu]} |z - w|)\,\, \to 0$
as
$\nu \to \infty$.
\end{lemma}

It follows from this lemma, Lemma~\ref{l2.3}, and the
$1$-Lipschitz
continuity of
$g$
that, in fact,
$|g(x) - g(y)| = |x - y|$
for all
$x,y \in \partial U$.
Hence, there is a Euclidean isometry
$Q : \mathbb R^n \to \mathbb R^n$
such that
$Q(x) = g(x)$
for all
$x \in \partial U$.
It is easily seen that
$Q(U) = V$.
The proof of Theorem~\ref{t1.2} is complete.

\section{Rigidity conditions for the boundaries of submanifolds in a
Riemannian manifold}\label{s3}

{\sc Rigidity problems and intrinsic geometry of submanifolds
in Riemannian manifolds}

Let
$(X,g)$
be an
$n$-dimensional
smooth connected Riemannian manifold without boundary and let
$Y$
be its
$n$-dimensional
compact connected
$C^0$-submanifold
with nonempty boundary
$\partial Y$
($n \ge 2$).

A classical object of investigations (see, for example,~\cite{A})
is given by the intrinsic metric
$\rho_{\partial Y}$
on a hypersurface
$\partial Y$
defined for
$x,y \in \partial Y$
as the infimum of the lengths of curves
$\nu \subset \partial Y$
joining
$x$
and
$y$.
In the recent decades, an alternative approach arose in the
rigidity theory for submanifolds of Riemannian manifolds (see,
for instance, the recent articles~\cite{Kor3},~\cite{Ko5},
and~\cite{Ko6}, which also contain a historical survey of works
on the topic). In accordance with this approach, the metric on
$\partial Y$
is induced by the intrinsic metric of the interior
$\inter Y$
of the submanifold
$Y$.

Namely, suppose that
$Y$
satisfies the following condition:

${\rm(i)}$
if
$x,y \in Y$,
then
\begin{equation}\label{eq3.1}
\rho_Y(x,y) = \liminf_{x' \to x, y' \to y; x',y' \in \inter Y}
\{\inf[l(\gamma_{x',y',\inter Y})]\} < \infty,
\end{equation}
where
$\inf[l(\gamma_{x',y',\inter Y})]$
is the infimum of the lengths
$l(\gamma_{x',y',\inter Y})$
of smooth paths
$\gamma_{x',y',\inter Y} : [0,1] \to \inter Y$
joining
$x'$
and
$y'$
in the interior
$\inter Y$
of Y.

\textbf{Remark~3.1.}
Easy examples show that if
$X$
is an
$n$-dimensional
connected smooth Riemannian manifold without boundary then an
$n$-dimensional
compact connected
$C^0$-submanifold
in
$X$
with nonempty boundary may fail to satisfy condition
${\rm(i)}$.
For
$n = 2$,
we have the following counterexample:

Let
$(X,g)$
be the space
$\mathbb R^2$
equipped with Euclidean metric and let
$Y$
be a closed Jordan domain in
$\mathbb R^2$
whose boundary is the union of the singleton
$\{0\}$
consisting of the origin
$0$,
the segment
$\{(1 - t)(e_1 + 2e_2) + t(e_1 + e_2) : 0 \le t \le 1\}$,
and
of the segments of the following four types:
$$
\biggl\{\frac{(1 - t)(e_1 + e_2)}{n} + \frac{te_1}{n + 1} :
0 \le t \le 1\biggr\}
\quad (n = 1,2,\dots);
$$
$$
\biggl\{\frac{e_1 + (1 - t)e_2}{n} : 0 \le t \le 1\biggr\}
\quad (n = 2,3,\dots);
$$
$$
\biggl\{\frac{(1 - t)(e_1 + 2e_2)}{n} +
\frac{2t(2e_1 + e_2)}{4n + 3} : 0 \le t \le 1\biggr\} \quad
(n = 1,2,\dots);
$$
$$
\biggl\{\frac{(1 - t)(e_1 + 2e_2)}{n + 1} +
\frac{2t(2e_1 + e_2)}{4n + 3} : 0 \le t \le 1\biggr\} \quad
(n = 1,2,\dots).
$$
Here
$e_1,e_2$ is the canonical basis in
$\mathbb R^2$.
By the construction of
$Y$,
we have
$\rho_Y(0,E) = \infty$
for every
$E \in Y \setminus \{0\}$
(see figure~1).

\begin{center}     
\begin{figure}[ht]  
\includegraphics{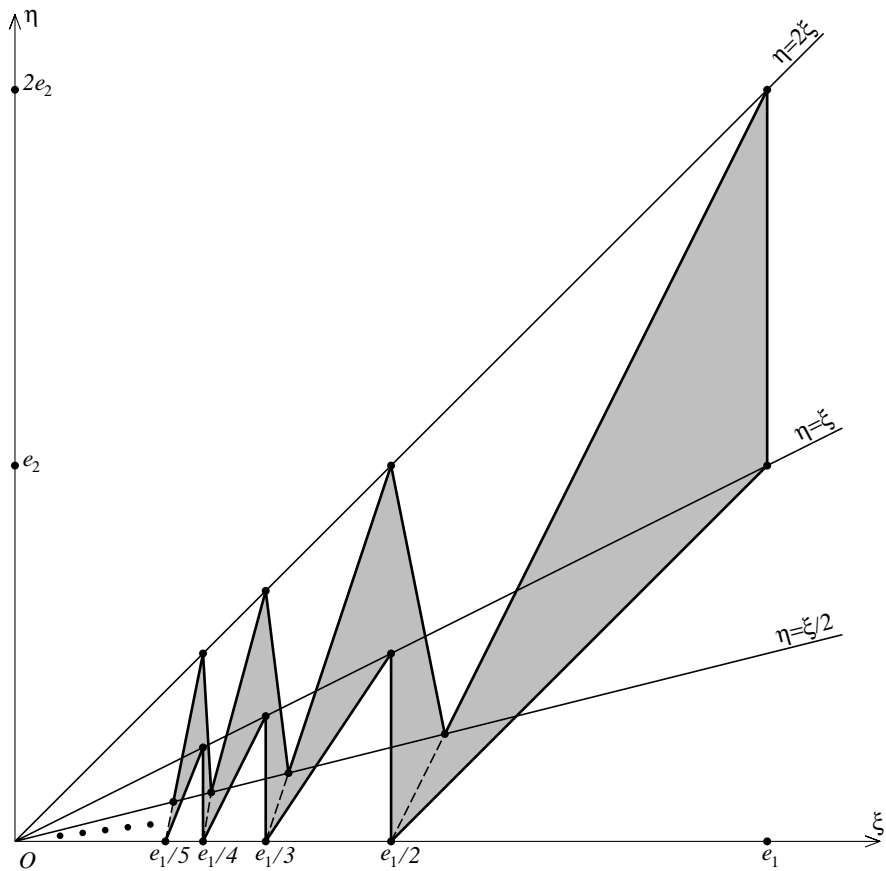}  
\caption{An example of
$2$-dimensional
compact connected
$C^0$-submanifold
with nonempty boundary which does not satisfy condition
${\rm(i)}$}  
\end{figure}       
\end{center}       

\textbf{Remark~3.2.}
Note that if
$X = \mathbb R^n$
and
$U$
is a domain in
$\mathbb R^n$
whose closure
$Y = \cl U$
is a Lipschitz manifold (such that
$\partial(\cl U) = \partial U \ne\varnothing$),
then
$\rho_{\partial U,U}(x,y) = \rho_Y(x,y)$
($x,y \in \partial U$)
and
$Y$
satisfies~
${\rm(i)}$.
Hence, this example is an important particular case of
submanifolds
$Y$
in a Riemannian manifold
$X$
satisfying
${\rm(i)}$.

To prove our rigidity results for boundaries of submanifolds in
a Riemannian manifold, we first need to study the
properties of the intrinsic geometry of these submanifolds.

One of the main results of this section is as follows:

\begin{theorem}\label{t3.1}
Let
$n = 2$.
Then, under condition
${\rm(i)},$
the function
$\rho_Y$
defined by~\ref{eq3.1} is a metric on
$Y$.
\end{theorem}

\textit{Proof.}
It suffices to prove that
$\rho_Y$
satisfies the triangle inequality. Let
$A$,
$O$,
and
$D$
be three points on the boundary of
$Y$
(note that this case is basic because the other cases are
simpler). Consider
$\varepsilon > 0$
and assume that
$\gamma_{A_\varepsilon O^1_\varepsilon} : [0,1] \to \inter Y$
and
$\gamma_{O^2_\varepsilon D_\varepsilon} : [2,3] \to \inter Y$
are smooth paths with the endpoints
$A_\varepsilon = \gamma_{A_\varepsilon O^1_\varepsilon}(0)$,
$O^1_\varepsilon = \gamma_{A_\varepsilon O^1_\varepsilon}(1)$
and
$D_\varepsilon = \gamma_{O^2_\varepsilon D_\varepsilon}(3)$,
$O^2_\varepsilon = \gamma_{O^2_\varepsilon D_\varepsilon}(2)$
satisfying the conditions
$\rho_X(A_\varepsilon,A) \le \varepsilon$,
$\rho_X(D_\varepsilon,D) \le \varepsilon$,
$\rho_X(O^j_\varepsilon,O) \le \varepsilon$
($j = 1,2$),
$|l(\gamma_{A_\varepsilon O^1_\varepsilon}) - \rho_Y(A,O)| \le
\varepsilon$,
and
$|l(\gamma_{O^2_\varepsilon D_\varepsilon}) - \rho_Y(O,D)| \le
\varepsilon$.
Let
$(U,h)$
be a chart of the manifold
$X$
such that
$U$
is an open neighborhood of the point
$O$
in
$X$,
$h(U)$
is the unit disk
$B(0,1) = \{(x_1,x_2) \in \mathbb R^2 : x^2_1 + x^2_2 < 1\}$
in
$\mathbb R^2$,
and
$h(O) = 0$
($0 = (0,0)$
is the origin in
$\mathbb R^2$);
moreover,
$h : U \to h(U)$
is a diffeomorphism having the following property: there exists
a chart
$(Z,\psi)$
of
$Y$
with
$\psi(O) = 0$,
$A,D \in U \setminus \cl_X Z$
($\cl_X Z$
is the closure of
$Z$
in the space
$(X,g)$)
and
$Z = \widetilde U \cap Y$
is the intersection of an open neighborhood
$\widetilde U$
($\subset U$)
of
$O$
in
$X$
and
$Y$
whose image
$\psi(Z)$
under
$\psi$
is the half-disk
$B_+(0,1) = \{(x_1,x_2) \in B(0,1) : x_1 \ge 0\}$.
Suppose that
$\sigma_r$
is an arc of the circle
$\partial B(0,r)$
which is a connected component of the set
$V \cap \partial B(0,r)$,
where
$V = h(Z)$
and
$0 < r < r^* = \min\{|h(\psi^{-1}(x_1,x_2))| : x^2_1 + x^2_2 = 1/4,
x_1 \ge 0\}$.
Among these components, there is at least one (preserve the notation
$\sigma_r$
for it) whose ends belong to the sets
$h(\psi^{-1}(\{-t e_2 : 0 < t < 1\}))$
and
$h(\psi^{-1}(\{t e_2 : 0 < t < 1\}))$
respectively. Otherwise, the closure of the connected component
of the set
$V \cap B(0,r)$
whose boundary contains the origin would contain a point
belonging to the arc
$\{e^{i\theta}/2 : |\theta| \le \pi/2\}$
(here we make use of the complex notation
$z = re^{i\theta}$
for points
$z \in \mathbb R^2$
($= \mathbb C$)).
But this is impossible. Therefore, the above-mentioned arc
$\sigma_r$
exists.

It is easy to check that if
$\varepsilon$
is sufficiently small then the images of the paths
$h \circ \gamma_{A_\varepsilon O^1_\varepsilon}$
and
$h \circ \gamma_{O^2_\varepsilon D_\varepsilon}$
also intersect the arc
$\sigma_r$,
i.e., there are
$t_1 \in\, ]0,1[$,
$t_2 \in\, ]2,3[$
such that
$\gamma_{A_\varepsilon O^1_\varepsilon}(t_1) = x^1 \in Z$,
$\gamma_{O^2_\varepsilon D_\varepsilon}(t_2) = x^2 \in Z$
and
$h(x^j) \in \sigma_r$,
$j = 1,2$.
Let
$\gamma_r : [t_1,t_2] \to \sigma_r$
be a smooth parametrization of the corresponding subarc of
$\sigma_r$,
i.e.,
$\gamma_r(t_j) = h(x^j)$,
$j = 1,2$.
Now we can define a mapping
$\widetilde \gamma_\varepsilon : [0,3] \to \inter Y$
by setting
$$
\widetilde \gamma_\varepsilon(t) =
\begin{cases}
\gamma_{A_\varepsilon O^1_\varepsilon}(t), & t \in [0,t_1]; \\
h^{-1}(\gamma_r(t)), & t \in ]t_1,t_2[;\\
\gamma_{O^2_\varepsilon D_\varepsilon}(t), & t \in [t_2,3].
\end{cases}
$$
By construction,
$\widetilde \gamma_\varepsilon$
is a piecewise smooth path joining the points
$A_\varepsilon = \widetilde \gamma_\varepsilon(0)$,
$D_\varepsilon = \widetilde \gamma_\varepsilon(3)$
in
$\inter Y$;
moreover,
$$
l(\widetilde \gamma_\varepsilon) \le
l(\gamma_{A_\varepsilon O^1_\varepsilon}) +
l(\gamma_{O^2_\varepsilon D_\varepsilon}) +
l(h^{-1}(\sigma_r)).
$$
By an appropriate choice of
$\varepsilon > 0$,
we can make
$r > 0$
arbitrarily small, and since a piecewise smooth path can be
approximated by smooth paths, we have
$\rho_Y(A,D) \le \rho_Y(A,O) + \rho_Y(O,D)$.

In connection with Theorem~\ref{t3.1}, there appears a natural
question: Are there analogs of this theorem for
$n \ge 3$?
The following Theorem~\ref{t3.2} answers this question in the
negative:

\begin{theorem}\label{t3.2}
If
$n \ge 3$
then there exists an
$n$-dimensional
compact connected
$C^0$-manifold
$Y \subset \mathbb R^n$
with nonempty boundary
$\partial Y$
such that the condition
${\rm(i)}$
{\rm(}where now
$X = \mathbb R^n)$
is fulfilled for
$Y$
but the function
$\rho_Y$
in this condition is not a metric on
$Y$.
\end{theorem}

\textit{Proof.}
It suffices to consider the case of
$n = 3$.
Suppose that
$A$,
$O$,
$D$
are points in
$\mathbb R^3$,
$O$
is the origin in
$\mathbb R^3$,
$|A| = |D| = 1$,
and the angle between the segments
$OA$
and
$OD$
is equal to
$\frac{\pi}{6}$.

The manifold
$Y$
will be constructed so that
$O \in \partial Y$,
and
$]O,A] \subset \inter Y$,
$]O,D] \subset \inter Y$.
Under these conditions,
$\rho_Y(O,A) = \rho_Y(O,D) = 1$.
However, the boundary of
$Y$
will create "obstacles" between
$A$
and
$D$
such that the length of any curve joining
$A$
and
$D$
in
$\inter Y$
will be greater than
$\frac{12}{5}$
(this means the violation of the triangle inequality for
$\rho_Y$).

Consider a countable collection of mutually disjoint segments\linebreak
$\{I^k_j\}_{j \in \mathbb N,k = 1,\dots,k_j}$
lying in the interior  of the triangle
$6\Delta AOD$
(which is obtained from the original triangle
$\Delta AOD$
by dilation with coefficient
$6$)
with the following properties:

$(*)$ every segment
$I^k_j = [x^k_j,y^k_j]$
lies on a ray starting at the origin,
$y^k_j = 11 x^k_j$,
and
$|x^k_j| = 2^{-j}$;

$(**)$ any curve
$\gamma$
with ends
$A$,
$D$
whose interior points lie in the interior of the triangle
$4\Delta AOD$
and belong to no segment
$I^k_j$,
satisfies the estimate
$l(\gamma) \ge 6$.

The existence of such a family of segments is certain: the
segments of the family must be situated chequer-wise so that any
curve disjoint from them be sawtooth, with the total length of
its "teeth" greater than
$6$
(it can clearly be made greater than any prescribed positive
number). However, below we exactly describe the construction.

It is easy to include the above-indicated family of segments in
the boundary
$\partial Y$
of
$Y$.
Thus, it creates a desired "obstacle" to joining
$A$
and
$D$
in the plane of
$\Delta AOD$.
But it makes no obstacle  to joining
$A$
and
$D$
in the space. The simplest way to create such a space obstacle
is as follows: Rotate each segment
$I^k_j$
along a spiral around the axis
$OA$.
Make the number of coils so large that the length of this spiral
be large its pitch (i.e., the distance between the origin and
the the end of a coil) be sufficiently small. Then the set
$S^k_j$
obtained as the result of the rotation of the segment
$I^k_j$
is diffeomorphic to a plane rectangle, and it lies in a small
neighborhood of the cone of revolution with axis
$AO$
containing the segment
$I^k_j$.
The last circumstance guarantees that the sets
$S^k_j$
are disjoint as before, and so (as above) it is easy to include
them in the boundary
$\partial Y$
but, due to the properties of the
$I^k_j$'s
and a large number of coils of the spirals
$S^k_j$,
any curve joining
$A$,
$D$
and disjoint from each
$S^k_j$
has length
$\ge \frac{12}{5}$.

We turn to an exact description of the constructions used. First
describe the construction of the family of segments
$I^k_j$.
They are chosen on the basis of the following observation:

Let
$\gamma : [0,1] \to 4\Delta AOD$
be any curve with ends
$\gamma(0) = A$,
$\gamma(1) = D$
whose interior points lie in the interior of the triangle
$4\Delta AOD$.
For
$j \in \mathbb N$,
put
$R_j = \{x \in 4\Delta AOD : |x| \in [8 \cdot 2^{-j},4 \cdot
3^{-j}]\}$.
It is clear that
$$
4\Delta AOD \setminus \{O\} = \cup_{j \in \mathbb N} R_j.
$$
Introduce the polar system of coordinates on the plane of the
triangle
$\Delta AOD$
with center
$O$
such that the coordinates of the points
$A,D$
are
$r = 1$,
$\varphi = 0$
and
$r = 1$,
$\varphi = \frac{\pi}{6}$,
respectively. Given a point
$x \in 6\Delta AOD$,
let
$\varphi_x$
be the angular coordinate of
$x$
in
$[0,\frac{\pi}{6}]$.
Let
$\Phi_j = \{\varphi_{\gamma(t)} : \gamma(t) \in R_j\}$.
Obviously, there is
$j_0 \in \mathbb N$
such that
\begin{equation}\label{eq3.2}
\mathcal H^1(\Phi_{j_0}) \ge 2^{-j_0}\frac{\pi}{6},
\end{equation}
where
$\mathcal H^1$
is the Hausdorff
$1$-measure.
This means that, while in the layer
$R_{j_0}$,
the curve
$\gamma$
covers the angular distance
$\ge 2^{-j_0}\frac{\pi}{6}$.
The segments
$I^k_j$
must be chosen such that~(\ref{eq3.2}) together with the condition
$$
\gamma(t) \cap I^k_j = \varnothing \quad \forall t \in [0,1],\,\,
\forall j \in \mathbb N,\,\, \forall k \in \{1,\dots, k_j\}
$$
give the desired estimate
$l(\gamma) \ge 6$.
To this end, it suffices to take
$k_j = [(2\pi)^j]$
(the integral part of
$(2\pi)^j$)
and
$$
I^k_j = \bigl\{x \in 6\Delta AOD : \varphi_x = k(2\pi)^{-j}\frac{\pi}{6},
|x| \in [11 \cdot 2^{-j},2^{-j}]\bigr\},
$$
$k = 1,\dots,k_j$.
Indeed, under this choice of the
$I^k_j$'s,
estimate~(\ref{eq3.2}) implies that
$\gamma$
must intersect at least
$(2\pi)^{j_0} 2^{-j_0} = \pi^{j_0} > 3^{j_0}$
of the figures
$$
U_k = \bigl\{x \in R_{j_0} : \varphi_x \in
\bigl]k(2\pi)^{-j_0}\frac{\pi}{6}, (k + 1)(2\pi)^{-j_0}\frac{\pi}{6}
\bigr[\bigr\}.
$$
Since these figures are separated by the segments
$I^k_{j_0}$
in the layer
$R_{j_0}$,
the curve
$\gamma$
must be disjoint from them each time in passing from one figure
to another. The number of these passages must be  at least
$3^{j_0} - 1$,
and a fragment of
$\gamma$
of length at least
$2 \cdot 3 \cdot 2^{-j_0}$
is required for each passage (because the ends of the segments
$I^k_{j_0}$
go beyond the boundary of the layer
$R_{j_0}$
containing the figures
$U_k$
at distance
$3 \cdot 2^{-j_0}$).
Thus, for all these passages, a section of
$\gamma$
is spent of length at least
$$
6 \cdot 2^{-j_0}(3^{j_0} - 1) \ge 6.
$$
Hence, the construction of the segments
$I^k_j$
satisfying
$(*)$ --- $(**)$
is finished.

Let us now describe the construction of the above-mentioned
space spirals.

For
$x \in \mathbb R^3$,
denote by
$\Pi_x$
the plane that passes through
$x$
and is perpendicular to the segment
$OA$.
On
$\Pi_{x^k_j}$,
consider the polar coordinates
$(\rho,\psi)$
with origin at the point of intersection of
$\Pi_{x^k_j}$
and
$[O,A]$
(in this system,
the point
$x^k_j$
has coordinates
$\rho = \rho^k_j$,
$\psi = 0$).
Suppose that a point
$x(\psi) \in \Pi_{x^k_j}$
moves along an Archimedean spiral, namely, the polar coordinates
of the point
$x(\psi)$
are
$\rho(\psi) = \rho^k_j - \varepsilon_j \psi$,
$\psi \in [0,2\pi M_j]$,
where
$\varepsilon_j$
is a small parameter to be specified below, and
$M_j \in \mathbb N$
is chosen so large that the length of any curve passing between
all coils of the spiral is at least
$10$.

Describe the choice of
$M_j$
more exactly. To this end, consider the points
$x(2\pi)$,
$x(2\pi(M_j - 1))$,
$x(2\pi M_j)$,
which are the ends of the first, penultimate, and last coils of
the spiral respectively (with
$x(0) = x^k_j$
taken as the starting point of the spiral). Then
$M_j$
is chosen so large that the following condition hold:

$(*_1)$
\textit{The length of any curve on the plane
$\Pi_{x^k_j,}$
joining the segments
$[x^k_j,x(2\pi)]$
and
$[x(2\pi(M_j - 1)),x(2\pi M_j)]$
and disjoint from the spiral
$\{x(\psi) : \psi \in [0,2\pi M_j]\},$
is at least
$10$}.

Figuratively speaking, the constructed spiral bounds a
"labyrinth", the mentioned segments are the entrance to and the
exit from this labyrinth, and thus any path through the has length
$\ge 10$.

Now, start rotating the entire segment
$I^k_j$
in space along the above-mentioned spiral, i.e., assume that
$I^k_j(\psi) = \{y = \lambda x(\psi) : \lambda \in [1,11]\}$.
Thus, the segment
$I^k_j(\psi)$
lies on the ray joining
$O$
with
$x(\psi)$
and has the same length as the original segment
$I^k_j = I^k_j(0)$.
Define the surface
$S^k_j = \cup_{\psi \in [0,2\pi M_j]}I^k_j(\psi)$.
This surface is diffeomorphic to a plane rectangle (strip).
Taking
$\varepsilon_j > 0$
sufficiently small, we may assume without loss of generality
that
$2\pi M_j \varepsilon_j$
is substantially less than
$\rho^k_j$;
moreover, that the surfaces
$S^k_j$
are mutually disjoint (obviously, the smallness of
$\varepsilon_j$
does not affect property
$(*_1)$
which in fact depends on
$M_j$).

Denote by
$y(\psi) = 11 x(\psi)$
the second end of the segment
$I^k_j(\psi)$.
Consider the trapezium
$P^k_j$
with vertices
$y^k_j$,
$x^k_j$,
$x(2\pi M_j)$,
$y(2\pi M_j)$
and sides
$I^k_j$,
$I^k_j(2\pi M_j)$,
$[x^k_j,x(2\pi M_j)]$,
and
$[y^k_j,y(2\pi M_j)]$
(the last two sides are parallel since they are perpendicular to
the segment
$AO$).
By construction,
$P^k_j$
lies on the plane
$AOD$;
moreover, taking
$\varepsilon_j$
sufficiently small, we can obtain the situation where the
trapeziums
$P^k_j$
are mutually disjoint (since
$P^k_j \to I^k_j$
under fixed
$M_j$
and
$\varepsilon_j \to 0$).
Take an arbitrary triangle whose vertices lie on
$P^k_j$
and such that one of these vertices is also a vertex at an acute
angle in
$P^k_j$.
By construction, this acute angle is at least
$\frac{\pi}{2} - \angle AOD  = \frac{\pi}{3}$.
Therefore, the ratio of the side of the triangle lying inside the
trapezium
$P^k_j$
to the sum of the other two sides (lying on the corresponding
sides of
$P^k_j$)
is at least
$\frac{1}{2}\sin\frac{\pi}{3} > \frac{2}{5}$.
If we consider the same ratio for the case of a triangle with a
vertex at an obtuse angle of
$P^k_j$
then it is greater than
$\frac{1}{2}$.
Thus we have the following property:

$(*_2)$
\textit{For arbitrary triangle whose vertices lie on the
trapezium
$P^k_j$
and one of these vertices is also a vertex in
$P^k_j,$
the sum of length of the sides situated on the corresponding
sides of
$P^k_j$
is less than
$\frac{5}{2}$
of the length of the third side {\rm(}lying inside
$P^k_j)$}.

Let a point
$x$
lie inside the cone
$K$
formed by the rotation of the angle
$\angle AOD$
around the ray
$OA$.
Denote by
$\Proj_{\rot} x$
the point of the angle
$\angle AOD$
which is the image of
$x$
under this rotation. Finally, let
$K_{4\Delta AOD}$
stand for the corresponding truncated cone obtained by the
rotation of the triangle
$4\Delta AOD$,
i.e.,
$K_{4\Delta AOD} = \{x \in K : \Proj_{\rot} x \in 4\Delta AOD\}$.

The key ingredient in the proof of our theorem is the following
assertion:

$(*_3)$
\textit{For arbitrary space curve
$\gamma$
of length less than
$10$
joining the points
$A$
and
$D,$
contained in the truncated cone
$K_{4\Delta AOD} \setminus \{O\},$
and disjoint from each strip
$S^k_j,$
there exists a plane curve
$\widetilde \gamma$
contained in the triangle
$4\Delta AOD \setminus \{O\},$
that joins
$A$
and
$D,$
is disjoint from all segments
$I^k_j$
and such that the length of
$\widetilde \gamma$
is less than
$\frac{5}{2}$
of the length of
$\Proj_{\rot} \gamma$.}

Prove
$(*_3)$.
Suppose that its hypotheses are fulfilled. In particular, assume
that the inclusion
$\Proj_{\rot} \gamma \subset 4\Delta AOD \setminus \{O\}$
holds. We need to modify
$\Proj_{\rot} \gamma$
so that the new curve be contained in the same set but be
disjoint from each of the
$I^k_j$'s.
The construction splits into several steps.

{\bf Step 1.}
If
$\Proj_{\rot} \gamma$
intersects a segment
$I^k_j$
then it necessarily intersects also at least one of the shorter
sides of
$P^k_j$.

Recall that, by construction,
$P^k_j = \Proj_{\rot} S^k_j$;
moreover,
$\gamma$
intersects no spiral strip
$S^k_j$.
If
$\Proj_{\rot} \gamma$
intersected
$P^k_j$
without intersecting its shorter sides then
$\gamma$
would pass through all coils of the corresponding spiral. Then,
by
$(*_1)$,
the length of the corresponding fragment of
$\gamma$
would be
$\ge 10$
in contradiction to our assumptions. Thus, the assertion of
step~1 is proved.

{\bf Step 2.}
Denote by
$\gamma_{P^k_j}$
the fragment of the plane curve
$\Proj_{\rot} \gamma$
beginning at the first point of its entrance into the trapezium
$P^k_j$
to the point of its exit from
$P^k_j$
(i.e., to its last intersection point with
$P^k_j$).
Then this fragment
$\gamma_{P^k_j}$
can be deformed changing the first and the last points so that
the corresponding fragment of the new curve lie entirely on the
union of the sides of
$P^k_j$;
moreover, its length is less than
$\frac{5}{2}$
of the length of
$\gamma_{P^k_j}$.

The assertion of step~2 immediately follows from the assertions
of step~1 and
$(*_2)$.

The assertion of step~2 in turn directly implies the desired
assertion
$(*_3)$.
The proof of
$(*_3)$
is finished.

Now, we are ready to pass to the final part of the proof of
Theorem~\ref{t3.2}.

$(*_4)$
\textit{The length of any space curve
$\gamma \subset \mathbb R^3 \setminus \{O\}$
joining
$A$
and
$D$
and disjoint from each strip
$S^k_j$
is at least
$\frac{12}{5}$.}

Prove the last assertion. Without loss of generality, we may
also assume that all interior points of
$\gamma$
are inside the cone
$K$
(otherwise the initial curve can be modified without any increase
of its length so that assumptions of
$(*_4)$
are still fulfilled and the modified curve lies in
$K$).
If
$\gamma$
is not included in the truncated cone
$K_{4\Delta AOD} \setminus \{O\}$
then
$\Proj_{\rot} \gamma$
intersects the segment
$[4A,4D]$;
consequently, the length of
$\gamma$
is at least
$2(4\sin \angle OAD - 1) = 2(4 \sin\frac{\pi}{3} - 1) =
2(2\sqrt 3 - 1) > 4$,
and the desired estimate is fulfilled. Similarly, if the length
of
$\gamma$
is at least
$10$
then the desired estimate is fulfilled automatically, and there
is nothing to prove. Hence, we may further assume without loss
of generality that
$\gamma$
is included in the truncated cone
$K_{4\Delta AOD} \setminus \{O\}$
and its length is less than
$10$.
Then, by
$(*_3)$,
there is a plane curve
$\widetilde \gamma$
contained in the triangle
$4\Delta AOD \setminus \{O\}$,
joining the points
$A$
and
$D$,
disjoint from each segment
$I^k_j$,
and such that the length of
$\widetilde \gamma$
is at most
$\frac{5}{2}$
of the length of
$\Proj_{\rot} \gamma$.
By property
$(**)$
of the family of segments
$I^k_j$,
the length of
$\widetilde \gamma$
is at least
$6$.
Consequently, the length of
$\Proj_{\rot} \gamma$
is at least
$\frac{12}{5}$,
which implies the desired estimate. Assertion
$(*_4)$
is proved.

The just-proven property
$(*_4)$
of the constructed objects implies Theorem~\ref{t3.2}. Indeed,
since the strips
$S^k_j$
are mutually disjoint and, outside every neighborhood of the
origin
$O$,
there are only finitely many of these strips, it is easy to
construct a
$C^0$-manifold
$Y \subset \mathbb R^3$
that is homeomorphic to a closed ball (i.e.,
$\partial Y$
is homeomorphic to a two-dimensional sphere) and has the
following properties:

(I)
$O \in \partial Y$,
$[A,O[ \cup [D,O[ \subset \inter Y$;

(II)
for every point
$y \in (\partial Y) \setminus \{O\}$,
there exists a neighborhood
$U(y)$
such that
$U(y) \cap \partial Y$
is
$C^1$-diffeomorphic
to the plane square
$[0,1]^2$;

(III)
$S^k_j \subset \partial Y$
for all
$j \in \mathbb N$,
$k = 1,\dots,k_j$.

The construction of
$Y$
with properties (I)---(III) can be carried out, for example, as
follows: As the surface of the zeroth step, take a sphere
containing
$O$
and such that
$A$
and
$D$
are inside the sphere. At the
$j$th
step, a small neighborhood of the point
$O$
of our surface is smoothly deformed so that the modified surface
is still smooth, homeomorphic to a sphere, and contains all
strips
$S^k_j$,
$k = 1,\dots,k_j$.
Besides, we make sure that, at the each step, the so-obtained
surface be disjoint from the half-intervals
$[A,O[$
and
$[D,O[$,
and, as above, contain all strips
$S^k_j$,
$i \le j$,
already included therein. Since the neighborhood we are deforming
contracts to the point
$O$
as
$j \to \infty$,
the so-constructed sequence of surfaces converges (for example,
in the Hausdorff metric) to a limit surface which is the
boundary of a
$C^0$-manifold
$Y$
with properties (I)---(III).

Property (I) guarantees that
$\rho_Y(A,O) = \rho_Y(A,D) = 1$
and
$\rho_Y(O,x) \le 1 + \rho_Y(A,x)$
for all
$x \in Y$.
Property (II) implies the estimate
$\rho_Y(x,y) < \infty$
for all
$x,y \in Y \setminus \{O\}$,
which, granted the previous estimate, yields
$\rho_Y(x,y) < \infty$
for all
$x,y \in Y$.
However, property (III) and the assertion
$(*_4)$
imply that
$\rho_Y(A,D) \ge \frac{12}{5} > 2 = \rho_Y(A,O) + \rho_Y(A,D)$.
Theorem~\ref{t3.2} is proved.

If
$\rho_Y$
is a metric (the dimension
$n$
($\ge 2$)
is arbitrary) then the question of the existence of geodesics is
solved in the following assertion, which implies that
$\rho_Y$
is an \textit{intrinsic metric} (see, for example, \S 6
from~\cite{Al}).

\begin{theorem}\label{t3.3}
Assume that
$\rho_Y$
is a finite function and is a metric on
$Y$.
Then any two points
$x,y \in Y$
can be joined in
$Y$
by a shortest curve
$\gamma : [0,L] \to Y$
in the metric
$\rho_Y;$
i.e.{\rm,}
$\gamma(0) = x,$
$\gamma(L) = y,$
and
\begin{equation}\label{eq3.3}
\rho_Y(\gamma(s),\gamma(t)) = t - s, \quad \forall s,t \in [0,L],
\quad s < t.
\end{equation}
\end{theorem}

\textit{Proof.}
Fix a pair of distinct points
$x,y \in Y$
and put
$L = \rho_Y(x,y)$.
Now, take a sequence of paths
$\gamma_j : [0,L] \to Y$
such that
$\gamma_j(0) = x_j$,
$\gamma_j(L) = y_j$,
$x_j \to x$,
$y_j \to y$,
and
$l(\gamma_j) \to L$
as
$j \to \infty$.
Without loss of generality, we may also assume that the
parameterizations of the curves
$\gamma_j$
are their natural parameterizations up to a factor (tending to
$1$)
and the mappings
$\gamma_j$
converge uniformly to a mapping
$\gamma : [0,L] \to Y$
with
$\gamma(0) = x$,
$\gamma(L) = y$.
By these assumptions,
\begin{equation}\label{eq3.4}
\lim_{j \to \infty} l(\gamma_j|_{[s,t]}) = t - s
\quad \forall s,t \in [0,L], \quad s < t.
\end{equation}

Take an arbitrary pair of numbers
$s,t \in [0,L]$,
$s < t$.
By construction, we have the convergence
$\gamma_j(s) \in \inter Y \to \gamma(s)$,
$\gamma_j(t) \in \inter Y \to \gamma(t)$
as
$j \to \infty$.
From here and the definition of the metric
$\rho_Y(\cdot,\cdot)$
it follows that
$\rho_Y(\gamma(s),\gamma(t)) \le \lim_{j \to \infty}
l(\gamma_j|_{[s,t]})$.
By~\ref{eq3.4},
\begin{equation}\label{eq3.5}
\rho_Y(\gamma(s),\gamma(t)) \le t - s \quad \forall s,t
\in [0,L],\,\, s < t.
\end{equation}
Prove that~(\ref{eq3.5}) is indeed an equality. Assume that
$\rho_Y(\gamma(s'),\gamma(t')) < t' - s'$
for some
$s',t' \in [0,L]$,
$s' < t'$.
Then applying the triangle inequality and then~(\ref{eq3.5}), we
infer
$$
\rho_Y(x,y) \le \rho_Y(x,\gamma(s')) +
\rho_Y(\gamma(s'),\gamma(t')) + \rho_Y(\gamma(t'),y) < s' +
(t' - s') + (L - t') = L,
$$
which contradicts the initial equality
$\rho_Y(x,y) = L$.
The so-obtained contradiction completes the proof of
identity~(\ref{eq3.3}). q.e.d.

\textbf{Remark~3.3.}
Identity~(\ref{eq3.3}) means that the curve of Theorem~\ref{t3.3}
is a geodesic in the metric
$\rho_Y$,
i.e., the length of its fragment between points
$\gamma(s)$,
$\gamma(t)$
calculated in
$\rho_Y$
is equal to
$\rho_Y(\gamma(s),\gamma(t)) = t - s$.
Nevertheless, if we compute the length of the above-mentioned
fragment of the curve in the initial Riemannian metric then this
length need not coincide with
$t - s$;
only the easily verifiable estimate
$l(\gamma|_{[s,t]}) \le t - s$
holds (see~(\ref{eq3.4})). In the general case, the equality
$l(\gamma|_{[s,t]}) = t - s$
can only be guaranteed if
$n = 2$
(if
$n \ge 3$
then the corresponding counterexample is constructed by analogy
with the counterexample in the proof of Theorem~\ref{t3.2}, see
above). In particular, though, by Theorem~\ref{t3.3}, the metric
$\rho_Y$
is always intrinsic in the sense of the definitions in \cite[\S 6]{Al},
the space
$(Y,\rho_Y)$
may fail to be \textit{a space with intrinsic metric} in the
sense of [ibid].

{\sc Rigidity theorems for the boundaries of submanifolds in
Riemannian manifolds}

As in the 1st part of Sec.~3, let
$(X,g)$
be an
$n$-dimensional
smooth connected Riemannian manifold without boundary and let
$\rho_X$
be its intrinsic metric (i.e., let
$\rho_X(x,y)$
be the infimum of the lengths
$l(\gamma_{x,y,X})$
of smooth paths
$\gamma_{x,y,X} : [0,1] \to X$
joining
points
$x$
and
$y$
in a manifold
$X$).

Assume that
$Y$
is an
$n$-dimensional
compact connected
$C^0$-submanifold
with nonempty boundary
$\partial Y$
in
$X$
satisfying condition
${\rm(i)}$
in the 1st part of Sec.~3, moreover,
$\rho_Y$
is a metric on
$Y$.
Then
$Y$
is called strictly convex in the metric
$\rho_Y$
if, for any
$\alpha,\beta \in Y$,
any shortest path
$\gamma = \gamma_{\alpha,\beta,Y} : [0,1] \to Y$
between
$\alpha$
and
$\beta$
(in the metric
$\rho_Y$)
satisfies
$\gamma(]0,1[) \subset \inter Y$.

\begin{theorem}\label{t3.4}
Let
$n = 2$.
Assume that condition~
${\rm(i)}$
holds for a
$2$-dimensional
compact connected
$C^0$-submanifold
$Y_1$
with nonempty boundary
$\partial Y_1$
of a
$2$-dimensional
smooth connected Riemannian manifold
$X$
without boundary which is strictly convex in the metric
$\rho_{Y_1}$.
Suppose that
$Y_2 \subset X$
is also a
$2$-dimensional
compact connected
$C^0$-submanifold
of
$X$
with
$\partial Y_2 \ne\varnothing$
satisfying
${\rm(i)};$
moreover{\rm,}
$\partial Y_1$
and
$\partial Y_2$
are isometric in the metrics
$\rho_{Y_j},$
for
$j = 1,2$.
Then{\rm,}
$Y_2$
is strictly convex with respect to the metric
$\rho_{Y_2}$.
\end{theorem}

\textit{Proof.}
Suppose that, for points
$x,y \in Y_2$,
there exists a shortest path
$\gamma_{x,y,Y_2} : [0,1] \to Y_2$
in the metric
$\rho_{Y_2}$
joining
$x$
and
$y$
and such that
$\{\gamma_{x,y,Y_2}(]0,1[)\} \cap \partial Y_2 \ne\varnothing$,
i.e.,
$x' = \gamma_{x,y,Y_2}(t') \in \{\gamma_{x,y,Y_2}(]0,1[) \cap
\partial Y_2\}$
for a point
$t' \in ]0,1[$.
By Theorem~\ref{t3.3} and the fact that
$Y_2$
is a
$2$-dimensional
compact connected
$C^0$-submanifold
in
$X$,
for a sufficiently small neighborhood of
$x'$
in
$Y_2$,
we can find points
$x_0,y_0 \in \partial Y_2$
and a shortest path
$\gamma_{x_0,y_0,Y_2} : [0,1] \to Y_2$
between
$x_0$
and
$y_0$
in the same metric satisfying the condition
$x' \in \{\gamma_{x_0,y_0,Y_2}(]0,1[) \cap \partial Y_2\}$.
Further, we will suppose that
$x = x_0$
and
$y = y_0$.

Now, assume that
$f : \partial Y_1 \to \partial Y_2$
is an isometry of
$\partial Y_1$
and
$\partial Y_2$
in the metrics
$\rho_{Y_1}$
and
$\rho_{Y_2}$
of the boundaries
$\partial Y_1$
and
$\partial Y_2$
of the submanifolds
$Y_1$
and
$Y_2$
of
$X$.
From Theorem~\ref{t3.3}, we have
$$
\rho_{Y_2}(x,x') + \rho_{Y_2}(x',y) = l_1 + l_2 = l =
\rho_{Y_2}(x,y).
$$
Since
$f$
is an isometry,
$$
\rho_{Y_1}(f^{-1}(x),f^{-1}(x')) + \rho_{Y_1}(f^{-1}(x'),f^{-1}(y))
= \rho_{Y_2}(x,x') + \rho_{Y_2}(x',y).
$$
Next, consider shortest paths
$\gamma_{f^{-1}(x),f^{-1}(x'),Y_1} : [0,1/2] \to Y_1$
and
$\gamma_{f^{-1}(x'),f^{-1}(y),Y_1} : [1/2,1] \to Y_1$
in
$\rho_{Y_1}$
between (respectively)
$f^{-1}(x)$
and
$f^{-1}(x')$
and
$f^{-1}(x')$
and
$f^{-1}(y)$,
and then construct a path
$\gamma : [0,1] \to Y_1$
by setting
$\gamma(t) = \gamma_{f^{-1}(x),f^{-1}(x'),Y_1}(t)$
if
$0 \le t < 1/2$
and
$= \gamma_{f^{-1}(x'),f^{-1}(y),Y_1}(t)$
for
$1/2 \le t \le 1$.
Let
$l_{Y_1}(\delta)$
be the length of a path
$\delta : [0,1] \to Y_1$
in the metric
$\rho_{Y_1}$.
Since
$\rho_{Y_1}$
is a metric on
$Y_1$,
it is not difficult to show that
$$
l_{Y_1}(\gamma) \le l_{Y_1}(\gamma_{f^{-1}(x),f^{-1}(x'),Y_1}) +
l_{Y_1}(\gamma_{f^{-1}(x'),f^{-1}(y),Y_1}) = l_1 + l_2.
$$
Hence
$\gamma$
is a shortest path in
$\rho_{Y_1}$
joining
$f^{-1}(x)$
and
$f^{-1}(y)$
in
$Y_1$.
This contradicts the strict convexity of
$Y_1$.
The theorem is proved.

\begin{corollary}\label{c3.1}
Suppose that the conditions of Theorem~{\rm\ref{t3.4}} hold and
the manifold
$X$
has the following property{\rm:}
$\rho_X(x,y) = \rho_Y(x,y)$
for any two points
$x$
and
$y$
from every
$2$-dimensional
compact connected
$C^0$-submanifold
$Y \subset X$
with
$\partial Y \ne\varnothing$
satisfying condition
${\rm(i)}$
and strictly convex with respect to the metric
$\rho_Y$.
Then,
$\partial Y_1$
and
$\partial Y_2$
are isometric in the metric
$\rho_X$
on the ambient manifold
$X$.
\end{corollary}

\textbf{Remark~3.4.}
The condition imposed on the manifold
$X$
in Corollary~\ref{c3.1} can be reformulated as follows: in this
manifold, every
$2$-dimensional
compact connected
$C^0$-submanifold
$Y$
with boundary satisfying condition
${\rm(i)}$
and strictly convex with respect to its intrinsic metric
$\rho_Y$
is a convex set in the ambient space
$X$
with respect to the metric
$\rho_X$
(for the notion of a convex set in a metric space the reader is
referred, for example, to~\cite{A}).

We have the following analog of Theorem~\ref{t3.4} and
Corollary~\ref{c3.1} (combined together) for
$n \ge 3$.

\begin{theorem}\label{t3.5}
Let
$n \ge 3$.
Suppose that
$(X,g)$
is an
$n$-dimensional
smooth connected Riemannian manifold without boundary and
$Y_1$
and
$Y_2$
are
$n$-dimensional
compact connected
$C^0$-submanifolds
with nonempty boundaries
$\partial Y_1$
and
$\partial Y_2$
satisfying conditions
${\rm(i)}$,

${\rm(ii)}$
$\rho_{Y_j}$
is a metric on
$Y_j$
$(j = 1,2)$,

and

${\rm(iii)}$
for any two points
$a,b \in Y_j$,
there exist points
$c,d \in \partial Y_j$
which can be joined in
$Y_j$
by a shortest path
$\gamma : [0,1] \to Y_j$
in the metric
$\rho_{Y_j}$
so that
$a,b \in \gamma([0,1])$.

Furthermore{\rm,} assume that
$Y_1$
is strictly convex in the metric
$\rho_{Y_1}$,
$X$
has the additional property that
$\rho_X(x,y) = \rho_Y(x,y)$
for any two points
$x$
and
$y$
in every
$n$-dimensional
compact connected
$C^0$-submanifold
$Y \subset X$
with
$\partial Y \ne\varnothing$
satisfying conditions~
${\rm(i)}$-${\rm(iii)}$
and strictly convex with respect to
$\rho_Y$
and the boundaries
$\partial Y_1$
and
$\partial Y_2$
of the submanifolds
$Y_1$
and
$Y_2$
are isometric with respect to the metrics
$\rho_{Y_j},$
where
$j = 1,2$.
Then,
$\partial Y_1$
and
$\partial Y_2$
are isometric with respect to
$\rho_X$.
\end{theorem}

\textbf{Remark~3.5.}
For a submanifold
$Y$
in
$X$,
${\rm(i)}$
and
${\rm(ii)}$
can be considered as conditions of generalized regularity near
its boundary.

\textbf{Remark~3.6.}
Theorem~\ref{t3.4}, Corollary~\ref{c3.1}, and Theorem~\ref{t3.5}
are closely related to a theorem of A.~D.~Aleksandrov about the
rigidity of the boundary
$\partial U$
of a strictly convex domain
$U$
in Euclidean
$n$-space
$\mathbb R^n$
($n \ge 2$)
by the relative metric
$\rho_{\partial U,U}$
on the boundary (A.~D.~Aleksandrov's theorem was first published
(with his consent) by V.~A.~Aleksandrov in~\cite{Al}). The following
is an important particular case of this theorem.

\begin{theorem}\label{t3.6}
Let
$U_1$
be a strictly convex domain in
$\mathbb R^n$
{\rm(}i.e.{\rm,} for any
$\alpha,\beta \in \cl U_1$
every shortest path
$\gamma = \gamma_{\alpha,\beta,\cl U_1} : [0,1] \to \cl U_1$
between
$\alpha$
and
$\beta$
{\rm(}in the metric
$\rho_{\cl U_1})$
satisfies
$\gamma(]0,1[) \subset U_1)$.
Assume that
$U_2 \subset \mathbb R^n$
is any domain whose closure is a Lipschitz manifold {\rm(}such
that
$\partial (\cl U_2) = \partial U_2 \ne\varnothing);$
moreover{\rm,}
$\partial U_1$
and
$\partial U_2$
are isometric {\rm(}globally{\rm)} in their relative metrics
$\rho_{\partial U_1,U_1}$
and
$\rho_{\partial U_2,U_2}$.
Then
$\partial U_1$
and
$\partial U_2$
are isometric in Euclidean metric.
\end{theorem}

We say that an
$n$-dimensional
compact (closed) connected
$C^0$-submanifold
$Y$
with boundary
$\partial Y \ne\varnothing$
of an
$n$-dimensional
smooth connected (respectively,
$n$-dimensional
smooth complete connected) Riemannian manifold
$X$
without boundary has property
${\rm(\circ)}$
if
$\gamma_{x,y,Y}(]0,1[) \subset \inter Y$
for any two points
$x,y \in \partial Y$
and for every shortest path
$\gamma_{x,y,Y} : [0,1] \to Y$
in the metric
$\rho_Y$
joining these points.

\begin{theorem}\label{t3.7}
Let
$n = 2$.
Suppose that
${\rm(i)}$
holds for a
$2$-dimensional
compact connected
$C^0$-submanifold
$Y_1$
with boundary
$\partial Y_1 \ne\varnothing$
in a
$2$-dimensional
smooth connected Riemannian manifold
$X$
without boundary{\rm;} moreover{\rm,}
$Y_1$
has property
${\rm(\circ)}$.
Assume that
$Y_2 \subset X$
is a
$2$-dimensional
compact connected
$C^0$-submanifold
with
$\partial Y_2 \ne\varnothing$
in
$X$
satisfying
${\rm(\circ)};$
moreover{\rm,}
$\partial Y_1$
has property
$(\circ)$.
Assume that
$Y_2 \subset X$
is a
$2$-dimensional
compact connected
$C^0$-submanifold
with
$\partial Y_2 \ne\varnothing$
in
$X$
and
$\partial Y_1$
and
$\partial Y_2$
are isometric in the metrics
$\rho_{Y_j}$
$(j = 1,2)$.
Then
$\partial Y_2$
also has the property
${\rm(\circ)}$.
\end{theorem}

This theorem has the following generalization.

\begin{theorem}\label{t3.8}
Let
$n = 2$.
Suppose that
${\rm(i)}$
holds for a
$2$-dimensional
closed connected
$C^0$-submanifold
$Y_1$
with boundary
$\partial Y_1$
$(\ne\varnothing)$
in a
$2$-dimensional
smooth complete connected Riemannian manifold
$X$
without boundary satisfying
${\rm(\circ)}$.
Assume that
$Y_2 \subset X$
is a
$2$-dimensional
closed connected
$C^0$-submanifold
with
$\partial Y_2 \ne\varnothing$
in
$X;$
moreover{\rm,}
$\partial Y_1$
and
$\partial Y_2$
are isometric in the metrics
$\rho_{Y_j}$
$(j = 1,2)$.
Then
$Y_2$
has the property
${\rm(\circ)}$
as well.
\end{theorem}

\begin{corollary}\label{c3.2}
{\rm(of Theorem~\ref{t3.7}).}
Assume that the hypothesis of Theorem{\rm~\ref{t3.7}} hold and
that the manifold
$X$
has the following property{\rm:}
$\rho_X(x,y) = \rho_Y(x,y)$
for any two points
$x$
and
$y$
on the boundary
$\partial Y$
of every
$2$-dimensional
compact connected
$C^0$-submanifold
$Y \subset X$
with
$\partial Y \ne\varnothing$
satisfying
${\rm(i)}$
and
${\rm(\circ)}$.
Then
$\partial Y_1$
and
$\partial Y_2$
are isometric in the metric
$\rho_X$
of the ambient manifold~
$X$.
\end{corollary}

\begin{corollary}\label{c3.3} {\rm(of Theorem~\ref{t3.8})}{\bf.}
Assume that the hypothesis of Theorem{\rm~\ref{t3.8}} hold and
that the manifold
$X$
has the following property{\rm:}
$\rho_X(x,y) = \rho_Y(x,y)$
for any two points
$x$
and
$y$
on the boundary
$\partial Y$
of every
$2$-dimensional
closed connected
$C^0$-submanifold
$Y \subset X$
with
$\partial Y \ne\varnothing$
satisfying
${\rm(i)}$
and
${\rm(\circ)}$.
Then
$\partial Y_1$
and
$\partial Y_2$
are isometric with respect to the metric
$\rho_X$.
\end{corollary}

\begin{theorem}\label{t3.9}
Let
$n \ge 3$.
Suppose that
$(X,g)$
is an
$n$-dimensional
smooth connected Riemannian manifold without boundary and
$Y_1$
and
$Y_2$
are
$n$-dimensional
compact connected
$C^0$-submanifolds
with nonempty boundaries
$\partial Y_1$
and
$\partial Y_2$
in
$X$
satisfying conditions
${\rm(i)}$
and
${\rm(ii)}$
{\rm(}in Theorem{\rm~\ref{t3.5}}{\rm)}. Assume that
$Y_1$
has property
${\rm(\circ)}$
and
$X$
satisfies the following condition{\rm:}
$\rho_X(x,y) = \rho_Y(x,y)$
for any two points
$x$
and
$y$
on the boundary
$\partial Y$
of every
$n$-dimensional
compact connected
$C^0$-submanifold
$Y \subset X$
with
$\partial Y \ne\varnothing$
satisfying
${\rm(i)},$
${\rm(ii)},$
and
${\rm(\circ)}$.
Suppose also that
$\partial Y_1$
and
$\partial Y_2$
are isometric in the metrics
$\rho_{Y_j},$
where
$j = 1,2$.
Then
$\partial Y_1$
and
$\partial Y_2$
are isometric in
$\rho_X$.
\end{theorem}

\begin{theorem}\label{t3.10}
Let
$n \ge 3$.
Suppose that
$(X,g)$
is an
$n$-dimensional
smooth complete connected Riemannian manifold without boundary
and
$Y_1$
and
$Y_2$
are
$n$-dimensional
closed connected
$C^0$-submanifolds
with nonempty boundaries
$\partial Y_1$
and
$\partial Y_2$
in
$X$
satisfying
${\rm(i)}$
and
${\rm(ii)}$.
Assume that
$\partial Y_1$
has property
${\rm(\circ)}$
and
$X$
satisfies the following condition{\rm:}
$\rho_X(x,y) = \rho_Y(x,y)$
for any two points
$x$
and
$y$
on the boundary
$\partial Y$
of every
$n$-dimensional
closed connected
$C^0$-submanifold
$Y$
with
$\partial Y \ne\varnothing$
in
$X$
satisfying
${\rm(i)},$
${\rm(ii)},$
and
${\rm(\circ)}$.
Suppose also that
$\partial Y_1$
and
$\partial Y_2$
are isometric in the metrics
$\rho_{Y_j}$
$(j = 1,2)$.
Then
$\partial Y_1$
and
$\partial Y_2$
are isometric in
$\rho_X$.
\end{theorem}

\textit{Proofs of Theorems{\rm~\ref{t3.5}}
and{\rm~\ref{t3.7}}-{\rm\ref{t3.10}}}
are similar to the proof of Theorem~\ref{t3.4} (Theorems~\ref{t3.5}
and~\ref{t3.7}-\ref{t3.10} can be proved using the corresponding
analogs of Theorems~\ref{t3.1} and~\ref{t3.3}).

\section{On unique determination of domains by the condition of
local isometry of boundaries in the relative metrics}\label{s4}

{\sc The case of plane domains}

The first main result of this Section is the following
theorem.

\begin{theorem}\label{t4.1}
Let
$U$
be a domain in
$\mathbb R^2$
with smooth boundary. Then

$(i)$ if
$U$
is bounded{\rm,} then it is uniquely determined in the class of
all bounded plane domains with smooth boundaries by the condition
of local isometry of boundaries in the relative metrics if and
only if this domain is convex;

$(ii)$ if
$U$
is unbounded{\rm,} then the unique determination of
$U$
in the class of all plane domains with smooth boundaries by the
condition of local isometry of boundaries in the relative metrics
is equivalent to the strict convexity of this domain.
\end{theorem}

\textbf{Remark~4.1.}
Let
$U$
be a domain in
$\mathbb R^n$.
As in~\cite{Ko1}, we say that
$U$
has smooth boundary, Lipschitz boundary if the Euclidean
boundary
$\partial U$
of this domain is an
$(n - 1)$-submanifold
of class
$C^1$
(a Lipschitz submanifold) without boundary in
$\mathbb R^n$.
In the case of domain
$U$
with Lipschitz boundary, its Hausdorff boundary
$\partial_H U$
is in natural way identified with Euclidean boundary and metric
$\rho_{\partial U,U}$,
corresponding to Hausdorff metric can be defined in the following
manner:
$$
\rho_{\partial U,U}(x,y) = \liminf_{x' \to x,y' \to y;\, x',y' \in U}
\{\inf[l(\gamma_{x',y',U})]\},
$$
where
$x,y \in \partial U$
and
$\inf[l(\gamma_{x',y',U})]$
is the infimum of lengths
$l(\gamma_{x',y',U})$
of smooth paths
$\gamma_{x',y',U}: [0,1] \to U$
joining
$x'$
and
$y'$
in
$U$.
Recall also that a domain
$U$
is said to be strictly convex if it is convex and the interior
of the segment joining any two points in its closure
$\cl U$
is contained in
$U$
(cf. with the hypothesis of Theorem~\ref{t3.6}).

\begin{lemma}\label{l4.1}
Let
$U$
and
$V$
be two plane domains with smooth boundaries and
$f: \partial U \to \partial V$
be a bijective mapping which is a local isometry of boundaries
of these domains in the relative metrics. Then
$f$
is a {\rm(}global{\rm)} isometry of boundaries
$\partial U$
and
$\partial V$
in their intrinsic metrics.
\end{lemma}

\begin{lemma}\label{l4.2}
Suppose that domains
$U$
and
$V$
and mapping
$f: \partial U \to \partial V$
satisfy to the hypothesis of Lemma~{\rm\ref{l4.1},} moreover{\rm,}
$\partial U$
is bounded. Then the boundary
$\partial V$
of the domain
$V$
is also bounded and the mapping
$f$
has the following property{\rm:} there exists a number
$\varepsilon > 0$
such that
$\rho_{\partial U,U}(a,b) = \rho_{\partial V,V}(f(a),f(b))$
if
$a,b \in \partial U$
and
$\rho_{\partial U,U}(a,b) < \varepsilon$.
\end{lemma}

\begin{lemma}\label{l4.3}
Under hypothesis of Lemma~{\rm\ref{l4.1}} and an additional
supposition that the boundary
$\partial U$
of the domain
$U$
is connected{\rm,} the boundary
$\partial V$
of the domain
$V$
is also connected.
\end{lemma}

The proofs of these lemmas are sufficiently simple. By this
reason, we omit them.
\vskip3mm

\textit{Proof of Theorem~{\rm\ref{t4.1}}.} {\bf Step 1.} Prove
the first part of assertion
$(i)$,
i.e., prove that if
$U$
is a bounded convex plane domain with smooth boundary then it is
uniquely determined in the class of all bounded plane domains with
smooth boundaries by the condition of local isometry of boundaries
in the relative metrics. To this end, suppose that for a bounded
convex plane domain
$U$
with smooth boundary, there exists a bounded plane domain
$V$
with smooth boundary such that its boundary
$\partial V$
is locally isometric to the boundary
$\partial U$
of
$U$
in the relative metrics of boundaries (further,
$f: \partial U \to \partial V$
is the fixed mapping realizing a such isometry). Then by
Lemmas~\ref{l4.2} and~\ref{l4.3}, the boundary of the domain
$V$
is connected (and consequently,
$V$
is a Jordan domain), and
$f$
has the property indicated in Lemma~\ref{l4.2}. Accomplishing, if
it is necessary, an additional inversion with respect to a
straight line, we can also assume that
$f: \partial U \to \partial V$
preserves the orientation of the boundary
$\partial U$
of
$U$,
induced by the canonical orientation of this domain, i.e.,
$f$
"transfers" the mentioned orientation to the orientation of the
boundary
$\partial V$
of
$V$
induced by the canonical orientation of
$V$.

Let, further,
$I = [a,b]$,
where
$a \ne b$,
be a segment such that
$I \subset \partial U$
and the image
$f(I)$
of
$I$
is no longer a segment, moreover, any another segment
$I^* = [a^*,b^*]$
($a^* \ne b^*$,
$I^* \subset \partial U$)
of
$\partial U$
having common points with
$I$
is a subset of the segment
$I$
(below, we denote the set of all such segments
$I$
by the symbol
$\Lambda$).
We assert that
$f(I)$
is a locally convex arc directed by its convexity inside the
domain
$V$.
The latter means that every point
$P \in f(I)$
has a closed neighborhood
$N = N(P)$
for which
$f(I) \cap N$
is a convex arc directed by its convexity inside
$V$,
i.e.,
$f(I) \cap N(P) = f(I_P)$,
where
$I_P = [\alpha_P,\beta_P] \subset I$,
and the closed curve
$C_P$,
composed from
$f(I_P)$
and the segment
$J_P$
joining the endpoints
$f(\alpha_P)$
and
$f(\beta_P)$
of the arc
$f(I_P)$,
either degenerates to the segment
$J_P$,
or is the boundary of a bounded convex domain with the following
property. There is found a segment
$T$
with
$\inter T \subset V$,
placed on straight line
$\tau_P$
which is perpendicular to
$J_P$
and passes through its midpoint, moreover, some endpoint of
$T$
belongs to the arc
$f(I_P)$
and its second endpoint is on the arc
$(\partial V) \setminus f(I_P)$,
both of these endpoints are on the same side from the straight
line
$j_P$
containing the segment
$J_P$,
and the endpoint belonging to the arc
$f(I_P)$
is nearer to
$j_P$
than the other endpoint. Assuming the contrary, i.e., supposing
that
$f(I)$
is no a locally convex arc directed by its convexity inside
$V$,
we (taking the smoothness of the boundary of
$V$
into account)
arrive to the existence of a segment
$I_P = [\alpha_P,\beta_P] \subset I$
such that either
$(1)$
$f(I_P)$
is a nonconvex arc, or
$(2)$
the arc
$f(I_P)$
is no a segment and is a convex arc directed by its convexity
inside the complement
$cV$
of
$V$.
In both cases, for the curve
$f(\inter I_P)$,
there exist a point
$Q \in f(\inter I_P)$
and a locally supporting segment to this curve from the side of
the complement
$cV$
of
$V$
all points of which, except the point
$Q$,
belong to the interior of
$cV$
and
$Q$
is the common point of this segment and the boundary
$\partial V$
of
$V$.
In the case of
$(2)$,
these point and supporting segment can be found on the basis of
the considerations using in the proof of the Theorem of
Leja-Wilkosz~\cite{LW} which is mentioned in~\cite{BZ}, if we
bring evident modifications corresponding to our case in it.

In the case of
$(1)$,
the curve
$f(\inter I_P)$
contains a point
$L_P$
such that if we draw the tangent in
$L_P$
to our curve then there exist points
$R_P \in f(\inter I_P)$
and
$S_P \in f(\inter I_P)$
lying on various sides of this tangent. Replace the point
$R_P$,
if it is necessary, by the point which is the nearest point to
$L_P$
on the segment
$[R_P,L_P]$
(we will remain for this point the designation
$R_P$)
and belongs to the arc
$f(\inter I_P)$.
Analogously, we replace the point
$S_P$
by the point of the arc
$f(\inter I_P)$
which is the nearest point to
$L_P$
on the segment
$[L_P,S_P]$.
And then consider two Jordan domains such that the boundary
of the first of them is the union of the segment
$[R_P,L_P]$
and of that arc from three arcs constituting the set
$f(\inter I_P) \setminus \{R_P,L_P\}$,
the endpoints of which are
$R_P$
and
$L_P$,
and the boundary of the second domain is constructed by the same
way on the basis of the points
$L_P$
and
$S_P$
and of the same arc
$f(\inter I_P)$.
By the way of constructing, one of these domains will be
contained in
$V$
and the second domain will be contained in
$cV$.
Considering the first domain of them and using the above-mentioned
argument from the proof of theorem of Leja-Wilkosz in~\cite{BZ},
it is not difficult to find the desired point
$Q$
on the part of its boundary disposed on
$f(\inter I_P)$,
and a locally supporting segment
$j$
to the curve
$f(\inter I_P)$
at that point from the side of the complement
$cV$
of the domain
$V$.
Hence, in both cases~(1) and~(2), we arrived to desired situation.
Transposing the tangent in the point
$Q$
to
$f(\inter I_P)$
in a parallel way to itself at a sufficiently short distance to
it to the side where, so to say, the domain
$V$
lies, we will easily get the following situation: there exist
three points
$R'_P$,
$L_P$
and
$S'_P$
on the boundary
$\partial V$
of
$V$
belonging to the arc
$f(\inter I_P)$
and such that
$]R'_P,S'_P[ \subset V$,
moreover, the segment
$[R'_P,S'_P]$
cuts off from
$V$
the Jordan subdomain the boundary of which contains
$L_P$.
Clear that
$f^{-1}(R'_P)$,
$f^{-1}(L_P)$
and
$f^{-1}(S'_P)$
are the successively ordered points on the interval
$\inter I_P$.
Hence, the triangle inequality holds for these points in the metric
$\rho_{\partial U,U}$,
but by their choice, for the points
$R'_P$,
$L_P$
and
$S'_P$,
the strict triangle inequality in the metric
$\rho_{\partial V,V}$
holds.
Since we could initially assume that the length of
$I_P$
is less than
$\varepsilon$,
where
$\varepsilon$
is the number from Lemma~\ref{l4.2} corresponding to the mapping
$f$
which is considered now then we arrived to the contradiction
because by this lemma, the equality in the triangle inequality
must also be fulfilled for the points
$R'_P$,
$L_P$
and
$S'_P$.
Therefore,
$f(I)$
is a locally convex arc directed by its convexity inside
$V$.

We assert that the set
$\Lambda$
is finite.
Clear that by the finiteness of the length
$l = l(\partial U)$
of the boundary
$\partial U$
of
$U$,
the finiteness of
$\Lambda$
follows from the fact that
$\Lambda$
does not contain segments the length of which does not exceed,
for example,
$\varepsilon/2$.
Assuming that a segment
$\Delta = [\alpha_{\Delta},\beta_{\Delta}] \in \Lambda$
has the length
$l(\Delta) \le \varepsilon/2$,
consider points
$Q$
and
$S$
of this segment such that
$Q \ne S$,
$Q$
is situated nearer, let us say, to the left endpoint
$\alpha_{\Delta}$
of the segment, and
$f(Q)$
and
$f(S)$
lie on the same side (and at the positive distance) from the
tangent
$\tau$
to
$\partial V$
at the point
$f(\alpha_{\Delta})$,
finally, the (least positive) angle between the tangent rays to
the arcs
$(\partial V) \setminus f(\Delta)$
and
$f([\alpha_{\Delta},S])$
in the points
$f(\alpha_{\Delta})$
and
$f(S)$,
respectively, is less than
$\pi/4$.
Further, let a point
$P \in (\partial U) \setminus \Delta$
is so close to
$\alpha_{\Delta}$
that
$\rho_{\partial U,U}(P,\alpha_{\Delta}) < \varepsilon/2$
and the points
$f(P)$
and
$f(\alpha_{\Delta})$
lie on the same side from each of the tangents to
$\partial V$
in the points
$f(Q)$
and
$f(S)$.
Under these suppositions, for the triple of the points
$P$,
$Q$
and
$S$,
the strict triangle inequality in the metric
$\rho_{\partial U,U}$
holds, and for their images
$f(P)$,
$f(Q)$
and
$f(S)$,
the triangle equality (in the metric
$\rho_{\partial V,V}$)
holds. Thereby, by virtue of the choice of the number
$\varepsilon$
(and Lemma~\ref{l4.2}), we arrive to the contradiction from which
it is follows that
$\Delta = \varnothing$
and consequently,
$\Lambda$
is finite.

Let
$\omega: [0,l] \to \partial U$
be a natural parametrization of the boundary
$\partial U$
of
$U$
corresponding to the orientation of
$\partial U$
generated by the canonical orientation of the domain
$U$,
and let
$[\alpha_1,\beta_1] \subset [0,l]$
and
$[\alpha_2,\beta_2] \subset [0,l]$
be the segments such that
$\omega([\alpha_j,\beta_j]) \in \Lambda$,
where
$j = 1,2$,
$\alpha_1  < \beta_1 < \alpha_2 <\beta_2$
and the arc
$\omega(]\beta_1,\alpha_2[)$
does not contain points of the segments from
$\Lambda$.
We assert that
$f|_{\omega([\beta_1,\alpha_2])}$
is an Euclidean isometry (i.e., there exists an Euclidean isometry
$F: \mathbb R^2 \to \mathbb R^2$
such that
$F|_{\omega([\beta_1,\alpha_2])} = f|_{\omega([\beta_1,\alpha_2])}$).
Indeed, if the arc
$\omega([\beta_1,\alpha_2])$
does not contain segments, i.e., it is strictly convex, moreover,
its convexity directed toward the interior of the complement
$cU$
of
$U$.
Hence, considering a point
$c \in \omega(]\beta_1,\alpha_2[)$
and sufficiently close to it points
$a,b \in \omega(]\beta_1,\alpha_2[)$,
where
$\beta_1 < \omega^{-1}(a) < \omega^{-1}(c) < \omega^{-1}(b) <
\alpha_2$
(the closeness of the points
$a$
and
$b$
to the point
$c$
is such that the distance between each pair of the considering
below triple of the points
$a$,
$f^{-1}(\gamma(s_0))$
and
$b$
is less than
$\varepsilon/2$;
we can easily obtain the latter using the hypothesis of theorem)
and assuming that
$[f(a),f(b)] \cap \inter (cV) \ne\varnothing$,
we arrive to a situation where for the shortest path
$\gamma: [0,s] \to \cl V$
joining the points
$f(a)$
and
$f(b)$
in the closure
$\cl V$
of the domain
$V$\footnote{The existence of such shortest path is
guaranteed, for instance, by the results of Section~3; see
also~\cite{KK}.}, there exists a point
$s_0 \in ]0,s[$
for which
$\gamma(s_0) \in \partial V$
and
$f^{-1}(\gamma(s_0)) \in f^{-1}(\Im \gamma \cap \partial V) \setminus
\{a,b\}$
($\ne\varnothing$).
But then for the triple of the points
$a$,
$f^{-1}(\gamma(s_0))$
and
$b$,
the strict triangle inequality in the metric
$\rho_{\partial U,U}$
holds,
at the same time for the points
$f(a)$,
$\gamma(s_0)$
and
$f(b)$,
takes place the equality in the triangle inequality in the metric
$\rho_{\partial V,V}$.
Therefore, by virtue of Lemma~\ref{l4.2},
$[f(a),f(b)] \subset \cl V$
from which the equality
$|f(a) - f(b)| = |a - b|$
follows. Hence, the restriction
$f|_{U_{\varepsilon} \cap \omega([\beta_1,\alpha_2])}$
of
$f$
to the intersection
$U_{\varepsilon} \cap \omega([\beta_1,\alpha_2])$
of the
$\varepsilon$-neighborhood
$U_{\varepsilon}$
($= B(P,\varepsilon)$)
of each point
$P \in \omega([\beta_1,\alpha_2])$
and arc
$\omega([\beta_1,\alpha_2])$
itself is an isometry in Euclidean metric. This circumstance
allows easily to conclude that
$f|_{\omega([\beta_1,\alpha_2])}$
is an Euclidean isometry. In the case where
$\omega([\beta_1,\alpha_2])$
contains segments (which no longer belong to the set
$\Lambda$
and, consequently, their images under the mapping
$f$
are also segments), the proof of the fact that
$f|_{\omega([\beta_1,\alpha_2])}$
is an Euclidean isometry is close to the proof of this fact in
the previous case, i.e., in the case of the strict convexity of
$\omega([\beta_1,\alpha_2])$.
The difference in the arguments consists of negligible and
easily reproducible details, and we omit them.

Now, we are able to conclude the proof of the first part of
assertion
$(i)$
of our theorem. If the boundary
$\partial U$
of
$U$
is such that
$\Lambda = \varnothing$,
then the first part of
$(i)$
is proved on the basis of the arguments from the previous item.
In the case of
$\Lambda \ne\varnothing$,
consider a segment
$\Delta \in \Lambda$
and accomplishing appropriate translation and rotation in the
plane
$\mathbb R^2$,
get the situation where the segment
$\Delta$
lies on the ordinate axis, its upper endpoint is the origin, and
the domain
$U$
is situated on the left half-plane. Let
$\gamma: [0,l] \to \partial U$
($\gamma(0) = \gamma(l) = (0,0)$)
be a natural parametrization of the boundary
$\partial U$
of
$U$
corresponding to the orientation of
$\partial U$
generated by the canonical orientation of
$U$.
If
$f|_{\gamma([0,l - l(\Delta)])}$
is an Euclidean isometry then we can assume, without loss of
generality, that
$f|_{\gamma([0,l - l(\Delta)])} = \Id_{\gamma([0,l - l(\Delta)])}$.
Taking yet into account that
$f(\gamma([l - l(\Delta),l])) = f(\Delta)$
is not a segment (because of
$\Delta \in \Lambda$),
we see that
$f(\partial U) = \partial V$
is not a closed curve
i.e.,
$f(\gamma(l)) \ne f(\gamma(0))$.
The obtained contradiction leads us to the conclusion that
$\Lambda = \varnothing$.
Thus, in this case, the first part of
$(i)$
is proved.
Further, assume that
$\Lambda$
consists of
$n$
segments
$[\gamma(\alpha_1),\gamma(\beta_1)]$,
$[\gamma(\alpha_2),\gamma(\beta_2)]$,
\ldots,
$[\gamma(\alpha_n),\gamma(\beta_n)] = \gamma([l - l(\Delta),l])
= \Delta$,
where
$0 < \alpha_1 < \beta_1 < \alpha_2 < \beta_2 < \dots < \alpha_n <
\beta_n = l$.
Since
$f|_{\gamma([0,\alpha_1])}$
is an Euclidean isometry, we can assume, with loss of generality,
that
$f|_{\gamma([0,\alpha_1])} = \Id_{\gamma([0,\alpha_1])}$.
Then, using the induction argument, it is not difficult to show
that the rotation of the vector
$\omega$,
where
$-\omega$
is the unit tangent vector to the curve
$\gamma([l - l(\Delta),l])$
(i.e., to the segment
$\Delta$) in the point
$\gamma(l)$,
is realized (under the action of the mapping
$f$)
at the angle being equal to the following quantity:
$$
V = -\sum_{k = 1}^n \biggl\{
\sup_{\alpha_k \le t_1 < t_2 < \dots < t_{\varkappa + 1} \le \beta_k}
\sum_{\nu = 1}^{\varkappa} |\theta\gamma(t_{\nu + 1}) -
\theta\gamma(t_{\nu})|\biggr\} \ne 0,
$$
where
$\theta\gamma(t)$
is the unit tangent vector to the curve
$\gamma([t,l])$
in the point
$\gamma(t)$
if
$0 < t < l$
and to the curve
$\gamma([0,l])$
in the point
$\gamma(0) = (0,0)$
when
$t = l$.
If
$|V| < 2\pi$
then
$\omega \ne \mu e_2$,
where
$\mu > 0$
and
$e_2$
is the unit base vector of the ordinate axis.
And if
$|V| \ge 2\pi$
then (since
$f$
preserves the orientation of the boundary) without
self-intersections, the curve
$f(\partial V)$
can not be close. In both cases, we got the contradiction with
the fact that the curve
$\partial V$
is closed and smooth. Hence, the first part of the assertion
$(i)$
of our theorem is completely proved.

{\bf Step 2.} Prove the second part of assertion
$(i)$.
Assuming that
$U$
is a bounded nonconvex plane domain with smooth boundary, we will
show that by the appropriate deformation, we can get another
domain
$V$
whose boundary
$\partial V$
is smooth and locally isometric to the boundary
$\partial U$
of
$U$
in the relative metrics
$\rho_{\partial U,U}$
and
$\rho_{\partial V,V}$
of boundaries, and the domains
$U$
and
$V$
themselves are not isometric each other in Euclidean metric,
i.e., there does not exist an Euclidean isometry
$J: \mathbb R^2 \to \mathbb R^2$
such that
$J(U) = V$.
In the case where the boundary of
$U$
is not connected, a construction of the above-mentioned domain
$V$
realizes by a small permutation of a connected component of the
boundary
$\partial U$
when the location of the other connected components leaves fixed.
And if the boundary of
$U$
is connected, i.e.,
$U$
is a Jordan domain, then we will argue in the following way.
By theorem of Leja-Wilkosz~\cite{LW}, there will be found a
"locally strict supporting outwards" segment
$I$
lying in
$U$
except a single interior point for
$I$,
let us say, point
$P$,
which belongs to
$\partial U$.
Consider a closed disk
$K$
centered at
$P$
and having so small radius
$r$
that the boundary circle of this disk intersects with
$I$
in two points and the interior of one of the half-disks
$K_+$
and
$K_-$
such that
$K_+ \cup K_- = K \setminus I$,
for instance,
$\inter K_-$,
does not contain points of the boundary
$\partial U$
of
$U$.
Let
$u$
and
$v$
are two straight lines which are perpendicular to
$I$,
situated on the various sides of the normal
$n$
to it at the point
$P$,
and sufficiently close to
$n$.
Let us consider the nearest points
$L$
and
$S$
of the sets
$u \cap \partial U$
and
$v \cap \partial U$
to the segment
$I$
and join them by the shortest in
$\cl U$
curve
$\mu$.
Moreover, we regard
$r$
so small that the closure of the arc representing itself lesser
of two arcs, which arise on the boundary
$\partial U$
when we remove the points
$L$
and
$S$
from it, is contained in
$(\inter K_+) \cup \{P\}$
and that (by the smoothness of
$\partial U$)
the curve
$\mu$
is convex, smooth and directed by its convexity toward
$U$.
We can get one of two cases:
$(1)$
$\mu \subset \partial U$,
and
$(2)$
$\mu$
contains segments interior of each of them is a subset of
$U$.
Further, consider (in both cases~(1) and~(2)) the points
$L^*$
and
$S^*$
belonging to
$\lambda \cap \mu$
and chosen by the following way:
$L^*$
and
$L$
lie on the same side of both the straight line containing the
segment
$I$
(moreover, the point
$L^*$
is situated nearer to this straight line than
$L$)
and the straight line
$\psi$
perpendicular to
$I$,
passing through the point
$P$,
besides,
$L^*$
is situated nearer to
$\psi$
than
$L$,
finally, the point
$S^*$
is defined in the similar way in comparison with the location of
$S$.
By symbol
$U^*$,
denote the Jordan domain with the boundary
$((\partial U) \setminus \lambda^*) \cup \mu^*$,
where
$\lambda^*$
and
$\mu^*$
are the subarcs of the arcs
$\lambda$
and
$\mu$
with the endpoints
$L^*$
and
$S^*$.

In the case of
$(1)$,
a necessary deformation of the domain
$U = U^*$
realizes in the sufficiently obvious way and is reduced to a
deformation of the curve
$\mu^*$.
The arc
$\mu^*$
replaces by a convex arc
$\widetilde {\mu^*}$
of the same length, lying in the disk
$K$
and also directed by its convexity in
$U^*$
(more exactly, in the new domain
$\widetilde {U} = \widetilde {U^*}$),
and the arc
$(\partial U^*) \setminus \lambda^* = (\partial U^*) \setminus \mu^*$
leaves fixed, moreover, the closed arc
$\widetilde {\mu^*} \cup ((\partial U^*) \setminus \mu^*)$
forms the smooth boundary of the new domain
$\widetilde {U^*}$.
It is not difficult to verify that the boundaries of
$U^*$
and
$\widetilde {U^*}$
are locally isometric in the relative metrics
$\rho_{\partial U^*,U^*}$
and
$\rho_{\partial \widetilde {U^*},\widetilde {U^*}}$
(here, as a local isometry in the relative metrics of the
boundaries of
$U^*$
and
$\widetilde {U^*}$,
we can take the mapping
$f$
of these boundaries which is leaving fixed the arc
$(\partial U^*) \setminus \mu^*$).
Clear also that in the process of the construction of our
deformation, we can get the following situation: it is impossible
to map the domain
$U^*$
onto the domain
$\widetilde {U^*}$
by an Euclidean isometry. Consequently,
$\widetilde {U^*}$
is the desired domain
$V$.

In the case of
$(2)$,
the construction of a new domain
$V$
realizes in the following way. If
$\mu^* = \lambda^*$
then
$V$
is constructed as in the case~(1). In the case where
$\mu^*$
contains segments the interior of which lie in
$U$
and their endpoints belong to
$\partial U$
(denote the set of all such segments by the symbol
$\mathcal M^*$),
we, starting from the domain
$U^*$,
realize first the construction of the domain
$\widetilde {U^*}$
circumscribed in the case~(1), but in addition, we will leave
invariant the length of every segment of
$\mathcal M^*$
under the action of the arising (in the process of the
construction) boundary mapping
$f^*: \partial U^* \to \widetilde {U^*}$.
The latter is possible by the large degree of freedom in the
construction of the curve
$\partial\widetilde {U^*}$
which is submitted by the condition to that the curve
$\mu^*$
satisfies in the case~(2)\footnote{In this connection, see
Lemma~\ref{l5.1}.}. In this case, the final mapping
$f: \partial U \to \partial \widetilde {U}$,
where
$\widetilde {U}$
($= V$)
is a desired new domain, is constructed like this: it leaves
fixed the curve
$(\partial U) \setminus \lambda$
and coincides with
$f^*$
on the set
$N = \mu^* \cap \lambda^*$.
And if the arc
$\chi$
with endpoints
$A$
and
$B$
has not common points with
$(\partial U) \setminus \lambda^*$
and is cut off from
$\partial U$
by a segment from
$\mathcal M^*$,
then we subject this curve to the action of the preserving
orientation Euclidean isometry
$J: \mathbb R^2 \to \mathbb R^2$
such that
$J(A) = f^*(A)$
and
$J(B) = f^*(B)$,
and then put
$f|_{\chi} = J|_{\chi}$.
In this case, the domain
$V$
is the Jordan domain with the boundary
$f(\partial U)$
and by the construction,
$f: \partial U \to \partial V$
is a local isometry of the boundaries of
$U$
and
$V$
in their relative metrics, moreover, the large degree of freedom
in the choice of the above-circumscribed deformation of the
domain
$U$
which still takes place, makes possible to realize this
deformation such that the domains
$U$
and
$V$
are not isometric in the Euclidean metric. So, in both cases
$(1)$
and
$(2)$,
we get the following situation: if
$U$
is nonconvex bounded plane domain with smooth boundary then it is
not uniquely determined in the class of all bounded plane domains
with smooth boundaries by the condition of local isometry of
boundaries in the relative metrics. Consequently, the assertion
$(i)$
of the theorem is completely proved.

{\bf Step 3.} Pass to prove the assertion
$(ii)$.
The fact, that an unbounded strictly convex plane domain
$U$
with smooth boundary is uniquely determined in the class of all
plane domains with smooth boundaries by the condition of local
isometry of boundaries in the relative metrics, can be proved on the
basis of the arguments from the proof of the first part of
assertion
$(i)$.
Considering one more plane domain
$V$
with smooth boundary and assuming that the boundaries of the
domains
$U$
and
$V$
are locally isometric in their relative metrics and modify
negligible the arguments from the proof of the assertion~(i), we
establish that
$\partial U$
and
$\partial V$
are isometric in the Euclidean metric from where the isometry of
the domains
$U$
and
$V$
themselves follows.

{\bf Step 4.} Proving the second part of assertion
$(ii)$,
at first, we make sure that if an unbounded plane domain
$U$
with smooth boundary is not convex then by the same method as we
used in the proof of the second part of assertion
$(i)$,
it can be deformed to a domain
$V$
with smooth boundary, moreover, to such domain that the
boundaries
$\partial U$
and
$\partial V$
are found to be locally isometric in their relative metrics, and
for these domains themselves, there does not exist an Euclidean
isometry
$J: \mathbb R^2 \to \mathbb R^2$
with property
$V = J(U)$.
In the considering case, there exists a deformation of the
boundary
$\partial U$
of the domain
$U$
which does not lead us to the desired result, but the
above-mentioned degree of freedom in a choice of a deformation
makes possible to pass easily over this difficulty.

{\bf Step 5.} Now, let
$U$
be an unbounded plane convex domain with smooth boundary which
is not strictly convex. In this case, a construction of an
above-mentioned domain
$V$
achieves by quite simple methods. Indeed, on the boundary
$\partial U$
of our domain, there exists a segment, let us say,
$I$.
We can assume that any other segment having common points with
$I$
and lying on
$\partial U$,
is a subset of
$I$.
Without loss of generality, we will also suppose that
$I$
is the segment of the abscissa axis with the endpoints
$A = (-2l,0)$
and
$B = (2l,0)$
and the domain
$U$
is found in the lower half-plane. Subject the boundary
$\partial U$
of the domain to the following transformation.
The origin divides the boundary for two curves. The curve from
those curves, which contains the segment with endpoints
$(0,0)$
and
$(0,2l)$,
leaves under this transformation fixed. The segment with
endpoints
$(-l,0)$
and
$(0,0)$
is transformed to the quarter of the circle
$\{(x,y) \in \mathbb R^2: x^2 + (y - \frac{2l}{\pi})^2 =
\frac{4l^2}{\pi^2}\}$
with endpoints
$(0,0)$
and
$P = (-\frac{2l}{\pi},\frac{2l}{\pi})$.
The remaining part of the boundary
$\partial U$
is first subjected to the translation at vector
$((1 - \frac{2}{\pi})l,\frac{2l}{\pi})$
and then to the rotation at angle
$-\frac{\pi}{2}$
with respect to
$P$.
As the final result, we get the curve
$\gamma$
dividing the plane on two unbounded domains. That domain from
them which locally adjoins from below to the segment with
endpoints
$(0,0)$
and
$(0,2l)$,
we will take for a domain
$V$.
Easily to verify that the boundary
$\partial V$
of this domain is locally isometric to the boundary
$\partial U$
of the initial domain
$U$
in the relative metrics
$\rho_{\partial U,U}$
and
$\rho_{\partial V,V}$
of boundaries, and the domains
$U$
and
$V$
themselves are not isometric to each other in Euclidean metric.
Thus, assertion
$(ii)$,
and together with it, Theorem~\ref{t4.1} are completely proved.

In connection with Theorem~\ref{t4.1}, it should be noted that
there exists a bounded plane domain
$U$
with smooth boundary which is not uniquely determined in the
class of all plane domains with smooth boundaries by the
condition of local isometry of boundaries in the relative
metrics~(see~\cite{Sl}).

In the case where the boundary of a domain
$U \subset \mathbb R^2$
is not smooth, Theorem~\ref{t4.1} ceases to be valid.
Really, the following assertion is correct.

\begin{theorem}\label{t4.2}
There exists a bounded plane domain
$U$
with Lipschitz boundary such that it is not convex but{\rm,}
at the same time{\rm,} is uniquely determined in the class of all
plane domains by the condition of local isometry of boundaries in
the relative metrics.
\end{theorem}

\textbf{Remark~4.2.}
Theorem~\ref{t4.2} due to M.~V. Korobkov (see~\cite{Ko2}).
An argument of its proof will be discussed below.

{\sc Unique determination of space domains}

Now, consider the case of space domains.
Below, we will use the following assertion which is a
generalization of Lemma~\ref{l2.3} to the case of
locally isometric mappings of the boundaries of domains.

\begin{lemma}\label{l4.4}
Let
$U,V$
be domains in space
$\mathbb R^n$
$(n \ge 2)$
such that there exists a bijective mapping
$f: \partial_H U \to \partial_H V$
local isometric in the relative metrics of their Hausdorff
boundaries. Then for every element
$w \in \partial_H U,$
there exists a number
$\varepsilon = \varepsilon_w > 0$
satisfying the following condition{\rm:} for any two elements
$a',b' \in \partial U$
such that
$]a',b'[ \subset U$
and the elements
$a,b \in \partial_H U$
generated by the path
$\gamma(t) = tb' + (1 - t)a',$
$t \in [0,1]$
{\rm(}i.e., generated by the Cauchy sequences in the intrinsic
metric
$\rho_U$
of the domain
$U$
$\{\gamma(1/n)\}_{n = 3,4,\dots}$
and
$\{\gamma(1 - 1/n)\}_{n = 3,4,\dots},$
respectively{\rm),} belong to the
$\varepsilon$-neighborhood
$Z(w) = \{z \in \partial_H U: \rho_{\partial_H U,U}(z,w) < \varepsilon\}$
of the element
$w,$
the relation
$]p_Vf(a),p_Vf(b)[ \subset V$
holds.
\end{lemma}

The proof of this lemma differs from the proof of Lemma~\ref{l2.3}
by negligible modifications, therefore, we omit it.
\vskip3mm

Now, suppose that a considering domain is strictly convex. Then
the following theorem is valid.

\begin{theorem}\label{t4.3}
Let
$n \ge 2$.
If a domain
$U$
in space
$\mathbb R^n$
is strictly convex{\rm,} then it is uniquely determined in the
class of all domains in this space by the condition of local
isometry of boundaries in the relative metrics.
\end{theorem}

\textit{Proof.}
Let
$V$
be a domain such that there exists a bijective mapping
$f: \partial_H U \to \partial_H V$
being a local isometry in the relative metrics of the Hausdorff
boundaries
$\partial_H U$
and
$\partial_H V$
of the domains
$U$
and
$V$.
Assume that
$x$
and
$y$
are points of the Euclidean boundary
$\partial U$
of the domain
$U$
(by the strict convexity of
$U$
and Remark~4.1, we can suppose that
$x$
and
$y$
are simultaneously elements of the Hausdorff boundary
$\partial_H U$
of
$U$).
By Lemma~\ref{l4.4}, each element
$w \in \partial_H U$
has an
$\varepsilon_w$-neighborhood
$Z(w) = \{z \in \partial_H U: \rho_{\partial_H U,U}(z,w) < \varepsilon\}$
with the property: for any points
$a,b \in Z(w)$
the relation
$]p_V f(a),p_V f(b)[\,\, \subset V$
holds (as for
$Z(w)$,
see Lemma~\ref{l4.4}).
From this fact, it follows that the mapping
$\bar{f}: \partial U \to \partial V$
such that
$\bar{f}(x) = p_V f(x)$
if
$x \in \partial U$
is locally isometric in Euclidean metric (i.e., if
$w \in \partial U$,
then for each point
$z \in Z(w)$,
there exist a ball
$B_x = B(x,r_x) \subset \mathbb R^n$
and an isometric mapping
$F_x: \mathbb R^n \to \mathbb R^n$
in the Euclidean metric such that
$F_x|_{(\partial U) \cup B_x} = \bar{f}|_{(\partial U) \cup B_x}$).

Let
$\bar{f}(\partial U) = T \subset \partial V$.
We assert that the closure
$\cl T$
of the set
$T$
coincides with the Euclidean boundary
$\partial V$
of
$V$.
Assuming that
$M = ((\partial V) \setminus \cl T) \ne\varnothing$,
consider a point
$z \in M$.
Since
$\cl T$
is a closed set then
$\dist \{z,T\} = \dist \{z, \cl T\} > 0$.
Taking yet into account that by Lemma~\ref{l2.1}, the set of supports
of the Hausdorff boundary of a domain is dense on its Euclidean
boundary, we can assert the existence of an element
$a$
of the Hausdorff boundary
$\partial_H V$
such that its support
$a' = p_V a$
satisfies to the condition
$\dist \{a',T\} = \dist \{a',\cl T\} > 0$.
Let
$\widetilde{a} = f^{-1}(a)$.
We have
$\bar {f}(\widetilde {a}) = p_V f (\widetilde {a}) = p_V(f(f^{-1}(a)
= p_V a = a' \in T$.
Therefore,
$\cl T = \partial V$.

Further, show that the mapping
$\bar {f}$
can be extended to an Euclidean isometry
$F: \mathbb R^n \to \mathbb R^n$
of all space
$\mathbb R^n$.
Indeed, let
$a$
and
$b$
be any two points on the Euclidean boundary
$\partial U$
of
$U$.
We will now establish that
\begin{equation}\label{eq4.1}
|\bar {f}(a) - \bar {f}(b)| = |a - b|.
\end{equation}
To this end, consider a path
$\gamma: [0,1] \to \partial U$
the endpoints of which are
$\gamma(0) = a$
and
$\gamma(1) = b$.
Since
$\bar {f}$
is a local isometry in the Euclidean metric then for each point
$t \in [0,1]$,
we can find a ball
$B_t = B(\bar {f}(\gamma(t)),r_t) \subset \mathbb R^n$
such that there exists an isometric in the Euclidean metric
mapping
$F_t: \mathbb R^n \to \mathbb R^n$
with the property
$F_t|_{(\partial U) \cap B_t} = \bar {f}|_{(\partial U) \cap B_t}$.
By the continuity of the path
$\gamma$,
the sets
$\gamma^{-1}((\partial U) \cap B_t)$,
where
$t \in [0,1]$,
form a covering of the segment
$[0,1]$
which is open with respect to
$[0,1]$.
But then we can extract a finite subcovering
$\{E_s = \gamma^{-1}((\partial U) \cap B_{t_s}), s = 1,\dots,k\}$.
If
$E_{s_1} \cap E_{s_2} \ne\varnothing$
where
$1 \le s_1,s_2 \le k$
then
$(\partial U) \cap B_{t_{s_1}} \cap B_{t_{s_2}} \ne\varnothing$.
Taking into account the strict convexity of the domain
$U$,
we easily conclude that
$F_{t_{s_1}} = F_{t_{s_2}}$.
Thereby, we can assert that there exists the single isometric
in the Euclidean metric mapping
$F: \mathbb R^n \to \mathbb R^n$
such that
$F_s = F$
for all
$s = 1,\dots,k$
and, consequently,
$\bar {f}|_{\Im \gamma} = F|_{\Im \gamma}$.
The latter implies the desired equality~(\ref{eq4.1}). And from
it, by its turn (with regard for the above-stated), the assertion
of the theorem follows.
\vspace{3mm}

\textbf{Remark~4.3.}
Theorem~\ref{t4.3} is a generalization of a theorem of A.~D.~Aleksandrov
about the unique determination of the boundary
$\partial U$
of a strictly convex domain
$U \subset \mathbb R^n$
by the relative metric
$\rho_{\partial U,U}$
(see~Theorem~\ref{t3.6}).

\textit{Proof of Theorem~{\rm\ref{t4.2}}.}
Assume that
$U$
is a bounded nonconvex domain in
$\mathbb R^2$
with Lipschitz boundary
$\partial U$,
and there exists such point
$P \in \partial U$
that on the set
$\partial U \setminus \{P\}$,
the domain
$U$
is locally strictly convex, moreover, its convexity directed to
the complement
$cU$
of this domain. We assert that
$U$
is uniquely determined by the condition of local isometry of
boundaries in the relative metrics. The proof of this assertion
realizes by the same scheme and with using the same tools as in
the proof of Theorem~\ref{t4.3}, with certain negligible
modifications. We turn our attention to them briefly.

In the considering now case let
$V$
be a one more domain in
$\mathbb R^2$
whose Hausdorff boundary is locally isometric to the Hausdorff
boundary of
$U$,
let
$f: \partial_H U \to \partial_H V$
be a bijection which is a local isometry in the relative metrics
of the Hausdorff boundaries
$\partial_H U$
and
$\partial_H V$
of
$U$
and
$V$,
finally, let
$T = \bar {f}((\partial U) \setminus \{P\})$
(since the boundary
$\partial U$
of
$U$
is Lipschitz, we, taking into account Remark~4.1, identify
$\partial_H U$
with
$\partial U$).
We assert that in this case, just as in the proof of
Theorem~\ref{t4.3},
$\cl T = \partial V$.
The latter can be proved on the basis of the arguments from the
proof of Theorem~\ref{t4.3}. Nevertheless, by Lemma~\ref{l2.1}
and the infinity of that part of the set of supports of the
Hausdorff boundary
$\partial_H V$
which is contained in the set
$M = (\partial V) \setminus \cl T$,
the indicated there point
$a'$
can be chosen so that
$\alpha = \bar {f}^{-1}(a)\,\, (\ne\varnothing) \subset (\partial U)
\setminus \{P\}$.

The further arguments iterate the arguments used in the proof of
Theorem~\ref{t4.3} almost verbatim. By this reason, we omit them.
\vspace{3mm}

As opposed to that what takes place in the case of domains in
$\mathbb R^2$
(see Theorem~\ref{t4.1}), in the case of space domains, under the
decision of problems on their unique determination by the
condition of local isometry of boundaries in the relative metrics,
the condition of convexity of a considerable domain (as in
Theorem~\ref{t4.2}) ceases to be necessary. In fact, the following
assertion holds.

\begin{theorem}\label{t4.4}
In
$\mathbb R^3,$
there exists a domain
$U$
with smooth boundary such that it is uniquely determined in the
class of all three-dimensional domains with smooth boundaries by
the condition of local isometry of boundaries in the relative
metrics but is not convex.
\end{theorem}

\textit{Proof.}
Let our domain
$U$
be made by the following way.

Consider the arc of cardioid
$$
\theta = \{(x,y,z) \in \mathbb R^3: x^2 + z^2 - \sqrt{x^2 + z^2} +
z = 0,\,\, x^2 + z^2 > 0,\,\, x \ge 0,\,\, y = 0\}.
$$
Leaving it fixed except of the part
$\theta_1$
which is cut out from it by the disk
$\{(x,y,z) \in \mathbb R^3: x^2 + z^2 \le \frac19,\,\,
y = 0\}$,
replace the arc
$\theta_1$
of the cardioid by the arc of the circle
$\{(x,y,z) \in \mathbb R^3: z = 1 -
\sqrt{\frac23 - x^2},\,\, 0 \le x \le \frac{\sqrt5}9,\,\, y = 0\}$.
It is not difficult to verify that under the rotation of the
curve obtained on this way around the axis
$Oz$
(up to the completed rotation), we obtain the closed smooth
surface being the boundary of a three-dimensional nonconvex
Jordan domain which we will accept for the desired domain
$U$,
establishing further that it is uniquely determined in the class
of all domains in
$\mathbb R^3$
with smooth boundaries by the condition of local isometry of
boundaries in the relative metrics.

So, let
$V \subset \mathbb R^3$
be another domain with smooth boundary and
$f: \partial U \to \partial V$
be a bijective mapping of the boundary
$\partial U$
of
$U$
onto the boundary
$\partial V$
of
$V$
which is a local isometry of the boundaries
$\partial U$
and
$\partial V$
in their relative metrics. Consider the curve
$\theta_0 = \theta \setminus \{(x,y,z) \in \mathbb R^3: x^2 +
z^2 \le 1/4,\,\, y = 0\}$.
Under the rotation around the axis
$Oz$,
this part of cardioid forms a locally strictly convex region
$S$
of the boundary of
$U$
directed by its convexity in the complement
$cU$
of this domain. Applying the same technique as in the proof of
Theorem~\ref{t4.3} and being based on Lemma~\ref{l4.4} in addition,
we first see for ourselves that there exists an isometry
$F: \mathbb R^3 \to \mathbb R^3$
in the Euclidean metric such that
$f|_S = F|_S$.

Without loss of generality, we can assume that
$F = \Id_{\mathbb R^3}$.
Suppose that
$S^*$
is the part of the boundary
$\partial U$
of
$U$
obtained under the rotation of the arc
$\theta^* = \cl(\theta \setminus (\theta_1 \cup \theta_0))$,
and consider the intersection of
$S^*$
with a closed half-plane for which the axis
$Oz$
is the boundary. We can also assume that this intersection is
the curve
$\theta^*$.
Now, we show that any two sufficiently close points
$a$
and
$b$
of this curve (note that the degree of closeness of these points
is determined by Lemma~\ref{l4.2} in application to the mapping
$f$)
cut out from it an arc
$ab$
the image of which under the mapping
$f$
is a plane curve. Indeed, considering the third point
$c$
of the arc
$ab$,
taking into account the local strict convexity of the arc
$\theta^*$
(with respect to the plane domain
$U_{x,z} = U \cap \{(x,y,z) \in \mathbb R^3: y = 0\}$,
moreover, by its convexity in the plane
$\tau_{x,y} = \{(x,y,z) \in \mathbb R^3: y = 0\}$,
this arc is directed to the side of the complement
$\tau_{x,y} \setminus U_{x,z}$
of
$U_{x,z}$),
and applying Lemma~\ref{l4.4} to each pair of the triple of points
$a$,
$c$
and
$b$,
we come to the conclusion that the point
$f(c)$
is on the surface
$\widetilde {S}$
formed by the rotation of the points of the arc
$f(ab) = f(a)f(b)$
around the straight line
$\zeta$
passing through the points
$f(a)$
and
$f(b)$,
and the intersection of
$\widetilde {S}$
with each half-plane, whose boundary is
$\zeta$,
has the same length as the arc
$ab$.
If we suppose that the arc
$f(ab)$
is not plane then its length will be greater than the length of the arc
$ab$.
The latter contradicts to Lemma~\ref{l4.1}. Hence, the arc
$f(ab)$
is plane. Making arguments close to those which is used in
the proof of the first part of assertion
$(ii)$
of Theorem~\ref{t4.1}, we establish further the existence of an
isometry
$F: \mathbb R^3 \to \mathbb R^3$
in the Euclidean metric such that
$F|_{ab} = f|_{ab}$.
Therefore, the arc
$f(ab)$
(together with the arc
$ab$)
is strictly convex and, consequently, if two planes contain the
arc
$f(ab)$
then they coincide. Extending our last considerations to the arc
$\theta^* \cup \theta_0$,
taking into account the above-said, and using the induction
argument, it is not difficult to establish that the curve
$f(\theta^* \cup \theta_0)$
is contained in the plane
$\tau_{x,y}$,
i.e., in the same plane that the curve
$\theta^* \cup \theta_0$.
Using again the arguments from the proof of Theorem~\ref{t4.1},
we come to the assertion that
$f|_{\theta^* \cup \theta_0} = \Id_{\theta^* \cup \theta_0}$.
Considering the rest intersections of the domain
$U$
with half-planes whose border is axis
$Oz$
and taking into account all above-stated, we obtain as the
result such relation
$$
f|_W = {\Id}_W
$$
where
$W$
is the part of the boundary
$\partial U$
of
$U$
which is obtained by the rotation of the arc
$\theta^* \cup \theta_0$
around the axis
$Oz$.

Assume that
$M = f((\partial U) \setminus W) \cap cV \cap \{(x,y,z) \in
\mathbb R^3: z \ge \frac29 \} \ne\varnothing$.
Let
$\alpha > \frac29$
and such that
$$
M_{\alpha} = M \cap \{(x,y,z) \in \mathbb R^3: z = \alpha \}
\ne\varnothing
$$
and
$$
M \cap \{(x,y,z) \in \mathbb R^3: z > \alpha \} = \varnothing.
$$
Suppose that, in
$M_{\alpha}$,
there exists a point
$(\bar {x},\bar {y},\alpha)$
such that
$\bar {x}^2 + \bar {y}^2 > 0$.
Without loss of generality, we can set that
$\bar {x}^2 + \bar {y}^2 = \max\limits_{(x,y,z) \in M_{\alpha}}
(x^2 + y^2)$.
Besides, since
$M_{\alpha} \cap f(W) = \varnothing$
then
$(\bar {x},\bar {y},\alpha) \not\in f(W)$.
Further, let
$\chi = \{\bar {x}(1 + \lambda/\sqrt{\bar {x}^2 + \bar {y}^2})e_1 +
\bar {y}(1 + \lambda/\sqrt{\bar {x}^2 + \bar {y}^2})e_2 +
(\alpha - \lambda t)e_3: \lambda \ge 0 \}$
be a ray outgoing from the point
$P_0 = (\bar {x}, \bar{y},\alpha)$,
moreover, the value of the parameter
$t\,\,(>0)$
is so small that this ray intersects
$f((\partial U) \setminus W) \setminus \{P_0\}$
and the distance between
$P_0$
and the nearest point
$P$
of the set
$(f((\partial U) \setminus W) \setminus \{P_0\}) \cap \chi$
to it is lesser than the number
$\varepsilon = \varepsilon_{P_0}$
from Lemma~\ref{l4.4} for the mapping
$f^{-1}$
(in this connection, note that the plane
$\tau_{\alpha} = \{(x,y,z) \in \mathbb R^3: z = \alpha\}$
is supporting to the surface
$f((\partial U) \setminus W)$
and, therefore, is a tangent plane to it in all points
$R \in M_{\alpha}$).
Consequently, by the lemma and the fact that the interval
$]P,P_0[$
is contained in
$V$,
the interval
$]f^{-1}(P),f^{-1}(P_0)[$
must be contained in
$U$.
But this is impossible. Therefore, it remains to consider
the case of
$\bar {x} = \bar {y} = 0$.
And in this case, we also have the contradiction, considering,
for example, the ray
$\{\lambda e_1 + (\alpha - \lambda t)e_3: \lambda \ge 0 \}$
as a desired ray, and further, repeating the arguments used in
the previous case.

We must yet discuss the case
$\alpha = \frac29$.
If
$\dist(M \cap \tau_{\frac29},W) >0$
then using the arguments from the previous two cases, we see
that this situation is also impossible.
Now, let
$\dist(M \cap \tau_{\frac29},W) = 0$.
The stated above facts and the smoothness of the boundaries
$\partial U$
and
$\partial V$
of
$U$
and
$V$
imply the following circumstance: for each point
$z^0 \in M_{\frac29}\,\, (= M \cap \tau_{\frac29})$,
there exists a number
$\varkappa_0 > 0$
such that any ray emitted from
$z^0$
and intersecting the cone
$K = \{(x,y,z) \in \mathbb R^3: z = \frac17 + \frac5{63}\sqrt{x^2 + y^2},\,
\frac17 \le z \le \frac29\}$
in a point contained between the planes
$\tau_{\frac29}$
and
$\tau_{\frac29 - \varkappa_0}$,
has common points with the surface
$(f((\partial U) \setminus W)) \setminus \{z^0\}$
(here we take into account that the generators of the cone
$K$
pass through the points of the boundary of the manifold
$\cl((\partial U) \setminus W)$,
being tangent in these points to the boundary
$\partial V$
of
$V$).
Choosing as
$z^0$
a point which is so near to
$W$
that the segment
$[z^0,\widetilde{z}]$
(where
$\widetilde{z}\in K\cap\tau_{\frac29-\frac{\varkappa_0}2}$)
of the ray
$\chi$
emitted from it and intersecting with the circle
$K \cap \tau_{\frac29 - \frac{\varkappa_0}{2}}$,
has the least of possible lengths of such segments, consider the
nearest point
$P \in (f((\partial U) \setminus W) \setminus \{z^0\}) \cap \chi$
to the point
$z^0$.
Setting in addition that
$|P - z^0| < \varepsilon_{z^0}$
(where
$\varepsilon_{z^0}$
is a number for the mapping
$f^{-1}$
from Lemma~\ref{l4.4}), we can apply the above-mentioned
arguments to make sure that this case is also impossible. At the
final result, we have the inequality
\begin{equation}\label{eq4.2}
f_3(x,y,z) < \frac29
\end{equation}
(where
$f = (f_1,f_2,f_3): \partial U \to \partial V$),
which holds for all points
$(x,y,z) \in (\partial U) \setminus W$.

Consider the bounded open set
$A \in \mathbb R^3$
whose boundary is composed from the sets
$f((\partial U) \setminus W)$
and
$\Xi = \{ (x,y,z) \in \mathbb R^3: x^2 + y^2 \le \frac5{81},\,\,
z = \frac29 \}$.
It is the three-dimensional Jordan domain contained in the
complement to
$V$.
Now, we will prove that the domain
$A$
is convex. Assume by contradiction that this is not valid. Using
the proof of theorem of Leja-Wilkosz~\cite{LW} that is set forth
in~\cite{BZ}, we can assert the existence of three points
$X \in \inter A$,
$Y \in \inter A$
and
$Z \in \inter A$
such that
$[X,Y] \subset \inter A$,
$[Y,Z] \subset \inter A$,
$[X,Z] \not\subset \inter A$,
starting from which and fixing the location of plane
$\tau$
containing these points, we can construct in this plane, for
instance, a locally supporting outwards
$A$
concave arc of ellipse
$\gamma$.
And then changing a location of the point
$Z$
in its small spherical neighborhood, we can obtain a continual
family of locally supporting outwards concave arcs of ellipses.
The plane measure of each part of the boundary
$\partial V$
of
$V$
which is found in one of the indicated plane intersections, can
not be positive, since
$\partial V$
is a smooth bounded surface and, consequently, has a finite area.
Therefore, there exist segments
$[a,b]$
of arbitrary small linear sizes such that
$]a,b[ \subset cA$
and
$a,b \in \partial A$,
moreover, we can also assume that
$a,b \not\in \Xi$.
Hence, we are again found that we have the above-discussed
situation in the process of proving of relation~(\ref{eq4.2})
from which it follows that the domain
$A$
is convex.

As the final result of our arguments for the surfaces
$\cl((\partial U) \setminus W)$
and
$f(\cl((\partial U \setminus W))$,
we are found themselves in the situation of theorem~2 of
Section 7 from Chapter 3 of monograph of
A.~V.~Pogorelov~\cite{Po}. Using it, we see that these surfaces
are equal. Taking into account the latter and also stated above
facts in the process of proving, we can assert that our theorem
is completely proved.

\section{Appendix}\label{s5}

\begin{lemma}\label{l5.1}
Let
$f_1: [0,a^*] \to \mathbb R$
$(a^* >0)$
be convex downwards strictly increasing smooth function such that
$f_1(0) = f'_1(0) = 0,$
moreover{\rm,} the graph
$\Gamma_1$
of this function contains straight line segments{\rm,} the union
of the set
$\mathcal M$
of all such segments is dense in
$\Gamma_1$
and
$(0,0)$
and
$(a^*,f_1(a^*))$
are limit points for the set of the left endpoints of the
segments from
$\mathcal M$
{\rm(}we assume that the segments
$\Delta \in \mathcal M$
are maximal in such sense that any segment
$\widetilde{\Delta} \subset \Gamma_1$
containing
$\Delta$
coincides with it{\rm)}. Then for each number
$\varepsilon > 0,$
there exists a convex downwards strictly increasing smooth function
$f_2: [0,a^*] \to \mathbb R$
differing from
$f_1$
and having the following properties{\rm:}
$||f_2 - f_1||_{C([0,a^*])} \le \varepsilon,$
$f_2(0) = f'_2(0) = 0,$
$f_2(a^*) = f_1(a^*),$
$f'_2(a^*) = f'_1(a^*)$
and the mapping
$F: \Gamma_1 \to \Gamma_2$
of the graphs of the functions
$f_1$
and
$f_2$
defined by the formula
$$
F: (x,y) \mapsto (\varphi^{-1}(x),f_2(\varphi^{-1}(f^{-1}_1(y)))) \in
\Gamma_2, \quad (x,y) \in \Gamma_1,
$$
where
$\varphi: [0,a^*] \to [0,a^*]$
is the diffeomorphic solution of the functional equation
$$
\int\limits_0^{\varphi(x)}\{1 + [f'_1(\varphi)]^2\}^{1/2} d \varphi
= \int\limits_0^x \{1 + [f'_2(t)]^2\}^{1/2} d t,
\quad 0 \le x \le a^*,
$$
is isometric in the intrinsic metrics of the curves
$\Gamma_1$
and
$\Gamma_2$
which transforms each straight line segment of
$\Gamma_1$
to a straight line segment of
$\Gamma_2$
with the same length.
\end{lemma}

\textit{Proof.}
Let
$x_1$,
$x_2$
and
$x_3$
be three points of the interval
$]0,a^*[$
such that
$x_1 < x_2 < x_3$
and these points are the left endpoints of the segments from the
set
$\mathcal M$
(the choice of the points
$x_1$,
$x_2$,
$x_3$
will be made more precise below). Assume that
$k_1$,
$k_2$,
$k_3$
and
$k_4$
are four real positive numbers. We will choose the function
$f_2$
among functions having the following form:
$$
f_2(x) =
\begin{cases}
k_1 f_1(x), & 0 \le x < x_1; \\
(k_1 - k_2)[f_1(x_1) + f'_1(x_1)(x - x_1)] + k_2 f_1(x),
& x_1 \le x < x_2;\\
\sum\limits_{s = 1}^2 (k_s - k_{s + 1})[f_1(x_s) +
f'_1 (x_s)(x - x_s)] + k_3 f_1(x), & x_2 \le x < x_3; \\
\sum\limits_{s = 1}^3 (k_s - k_{s + 1})[f_1(x_s) + f'_1(x_s)(x - x_s)]
+ k_4 f_1(x), & x_3 \le x \le a^*.
\end{cases}
$$

The equalities
$f_2(a^*) = f_1(a^*)$
and
$f'_2(a^*) = f'_1(a^*)$
leads us to the conditions
\begin{equation}\label{eq5.1}
\sum\limits_{s = 1}^3(k_s - k_{s + 1})[f_1(x_s) +
f'_1(x_s)(a^* - x_s)] + (k_4 - 1)f_1(a^*) = 0
\end{equation}
and
\begin{equation}\label{eq5.2}
\sum\limits_{s = 1}^3(k_s - k_{s + 1}) f'_1(x_s) +
(k_4 - 1) f'_1(a^*) = 0.
\end{equation}

The last condition will be result of the demand
$\varphi(a^*) = a^*$.
And since this demand is equality
$$
\int_0^{a^*}\{1 + [f'_1(t)]^2\}^{1/2} d t = \int_0^{a^*}\{1 +
[f'_2(t)]^2\}^{1/2} d t
$$
then we have
\begin{multline}\label{eq5.3}
\int_0^{a^*}\{1 + [f'_1(t)]^2\}^{1/2} -
\\
\sum\limits_{j = 0}^3 \int\limits_{x_j}^{x_{j + 1}}\biggl\{1 +
\biggl[\sum\limits_{s = 1}^j(k_s - k_{s + 1})f'_1(x_s) +
k_{j + 1}f'_1(t)\biggr]^2\biggr\}^{1/2} d t = 0,
\end{multline}
where
$x_0 = 0$,
$x_4 = a^*$
and
$\sum\limits_{s = 1}^0 \dots = 0$.

The element
$(k_1,k_2,k_3,k_4) = (1,1,1,1) \in \mathbb R^4$
is a solution to the system~(\ref{eq5.1})-(\ref{eq5.3}). At the
same time, by the construction, each straight line segment
$\Delta$
of the curve
$\Gamma_1$
is transformed to a straight line segment of the curve
$\Gamma_2$,
moreover,
$l(F(\Delta))=l(\Delta)$.
Now, it is sufficient to prove that the rank of the Jacobi matrix
of the left parts of the equalities~(\ref{eq5.1})-(\ref{eq5.3})
calculated with respect to the variables
$k_1$,
$k_2$,
$k_3$
and
$k_4$
in the point
$(1,1,1,1)$
is equal to
$3$
under the successful choice of
$x_1$,
$x_2$
and
$x_3$.

To this end, represent the mentioned matrix in the following
form:
\begin{equation}\label{eq5.4}
N = (A_{j s})_{\begin{array}
{l}
j = 1,2,3 \\
s = 1,2,3,4 \end{array}},
\end{equation}
where
$$
A_{11} = -u_1 - f'_1(x_1) =
-\int_0^{x_1} \frac{[f'_1(t)]^2 d t}{\{1 + [f'_1(t)]^2\}^{1/2}} -
f'_1(x_1) \int_{x_1}^{a^*} \frac{f'_1(t) d t}{\{1 +
[f'_1(t)]^2\}^{1/2}},
$$
$$
A_{12} = - \int_{x_1}^{x_2} \frac{f'_1(t)[f'_1(t) -
f'_1(x_1)] d t}{\{1 + [f'_1(t)]^2\}^{1/2}} -
\int_{x_2}^{a^*} \frac{f'_1(t)[f'_1(x_2) -
f'_1(x_1)] d t}{\{1 + [f'_1(t)]^2\}^{1/2}},
$$
$$
A_{13} = - \int_{x_2}^{x_3} \frac{f'_1(t)[f'_1(t) -
f'_1(x_2)] d t}{\{1 + [f'_1(t)]^2\}^{1/2}} -
\int_{x_3}^{a^*} \frac{f'_1(t)[f'_1(x_3) -
f'_1(x_2)] d t}{\{1 + [f'_1(t)]^2\}^{1/2}},
$$
$$
A_{14} = - \int_{x_3}^{a^*} \frac{f'_1(t)[f'_1(t) -
f'_1(x_3)] d t}{\{1 + [f'_1(t)]^2\}^{1/2}},
$$
$$
A_{21} = f_1(x_1) + (a^* - x_1)f'_1(x_1),
$$
$$
A_{22} = f_1(x_2) - f_1(x_1) + (a^* - x_2)f'_1(x_2) -
(a^* - x_1)f'_1(x_1),
$$
$$
A_{23} = f_1(x_3) - f_1(x_2) + (a^* - x_3)f'_1(x_3) - (a^* -x_2)f'_1(x_2),
$$
$$
A_{24} = f_1(a^*) - f_1(x_3) - (a^* - x_3)f'_1(x_3), \,\,\,
A_{31} = f'_1(x_1),
$$
$$
A_{32} = f'_1(x_2) - f'_1(x_1), \,\,\,
A_{33} = f'_1(x_3) - f'_1(x_2), \,\,\, A_{34} = f'_1(a^*) - f'_1(x_3).
$$

The rank of the matrix~(\ref{eq5.4}) coincides with the rank of
the matrix
$$
\widetilde{N} = \biggl(\sum_{\nu = 1}^s A_{j\nu}\biggr)_{
\begin{array}{l} j = 1,2,3 \\ s = 1,2,3,4 \end{array}},
$$
in which
\begin{multline*}
\sum_{\nu = 1}^2 A_{1\nu} = - u_2 - f'_1(x_2)v_2 =
\\
- \int_0^{x_2} \frac{[f'_1(t)]^2 d t}{\{1 + [f'_1(t)]^2\}^{1/2}}\,\,
-f'_1(x_2) \int_{x_2}^{a^*} \frac{f'_1(t) d t}{\{1 + [f'_1(t)]^2\}^{1/2}},
\end{multline*}
\begin{multline*}
\sum_{\nu = 1}^3 A_{1\nu} = - u_3 - f'_1(x_3)v_3 =
\\
- \int_0^{x_3} \frac{[f'_1(t)]^2 d t}{\{1 + [f'_1(t)]^2\}^{1/2}} \,\, -
f'_1(x_3) \int_{x_3}^{a^*} \frac{f'_1(t) d t}{\{1 +
[f'_1(t)]^2\}^{1/2}},
\end{multline*}
$$
\sum_{\nu =1}^4 A_{1\nu} = - u_4 = - \int_0^{a^*} \frac{[f'_1(t)]^2
d t}
{\{1 + [f'_1(t)]^2\}^{1/2}},
$$
\begin{multline*}
\sum_{\nu = 1}^2 A_{2\nu} = f_1(x_2) + (a^* - x_2)f'_1(x_2),
\,\,\,
\sum_{\nu = 1}^3 A_{2\nu} = f_1(x_3) + (a^* - x_3)f'_1(x_3),
\\
\sum_{\nu = 1}^4 A_{2\nu} = f_1(a^*), \,\,\,
\sum_{\nu = 1}^2 A_{3\nu} = f'_1(x_2), \,\,\, \sum_{\nu = 1}^3
A_{3\nu} = f'_1(x_3), \,\,\, \sum_{\nu = 1}^4 A_{3\nu} = f'_1(a^*).
\end{multline*}

Consider the determinant
\begin{multline*}
\delta_1 = \det \biggl\{\biggl(\sum_{\nu = 1}^s A_{j\nu}\biggr)_{
\begin{array}{l}
j = 2,3 \\ s = 3,4 \end{array}}\biggr\} = [f_1(x_3) +
(a^* - x_3)f'_1(x_3)]f'_1(a^*) -
\\
f_1(a^*)f'_1(x_3) =
f'_1(x_3)f'_1(a^*)\biggl\{\frac{f_1(x_3)}{f'_1(x_3)} + a^* - x_3
- \frac{f_1(a^*)}{f'_1(a^*)}\biggr\}
\end{multline*}
($0 < f'_1(x_3) < f'_1(a^*)$
by the hypothesis of the lemma and the choice of the points
$x_1$,
$x_2$
and
$x_3$).
The second factor in the right part of the last equalities
subject to the following transformations:
\begin{multline}\label{eq5.5}
\frac{\delta_1}{f'_1(a^*) f'_1(x_3)} =
\frac{f_1(x_3)}{f'_1(x_3)} -
\frac{f_1(a^*)}{f'_1(a^*)} +
a^* - x_3 =
\\
- \frac{f_1(a^*) - f_1(x_3)}{f'_1(a^*)} -
f_1(x_3)\biggl(\frac{1}{f'_1(a^*)}
- \frac{1}{f'_1(x_3)}\biggr) +
a^* - x_3
=
\\
- \frac{f'_1(\theta)(a^* - x_3)}{f'_1(a^*)} -
f_1(x_3)\frac{f'_1(x_3) - f'_1(a^*)}{f'_1(a^*)f'_1(x_3)} +
a^* - x_3 =
\\
- \biggl\{\frac{f'_1(\theta) - f'_1(a^*)}{f'_1(x_3) - f'_1(a^*)}
(a^* - x_3) + \frac{f_1(x_3)}{f'_1(x_3)}\biggr\}\frac{f'_1(x_3) -
f'_1(a^*)}{f'_1(a^*)}
\end{multline}
($x_3 < \theta < a^*$).
The convexity of
$f_1$
implies that
$$
\biggl|\frac{f'_1(\theta) - f'_1(a^*)}{f'_1(x_3) -
f'_1(a^*)}\biggr| \le 1.
$$
Therefore, if the condition
\begin{equation}\label{eq5.6}
a^* - x_3 < \frac{f_1(x_3)}{f'_1(x_3)}
\end{equation}
holds then
$\delta \ne 0$
and, consequently,
$\rank \widetilde{N} = \rank N \ge 2$.

Analogously, we can establish that if the point
$x_1$
is fixed and the point
$x_2$
($> x_1$)
is so near to
$x_1$
that the condition
\begin{equation}\label{eq5.7}
x_2 - x_1 < \frac{f_1(x_1)}{f'_1(x_1)}
\end{equation}
holds, then
\begin{multline*}
\delta_2 = \det \biggl\{\biggl(\sum_{\nu = 1}^s A_{j\nu}\biggr)_{
\begin{array}{l} j = 2,3 \\ s = 1,2 \end{array}}\biggr\} = \{f_1(x_1) +
(a^* - x_1)f'_1(x_1)\}f'_1(x_2) -
\\
f'_1(x_1)\{f_1(x_2) +
(a^* - x_2)f'_1(x_2)\} =
f'_1(x_1)f'_1(x_2)\biggl\{\frac{f_1(x_1)}{f'_1(x_1)} -
\frac{f_1(x_2)}{f'_1(x_2)} +
\\
x_2 - x_1\biggr\} \ne 0.
\end{multline*}

If we suppose that the first row of the matrix
$\widetilde{N}$
is a linear combination of two other rows of
$\widetilde{N}$
then we obtain two pairs of the relations
$$
- \frac{u_3}{f'_1(x_3)} - v_3 =
C_1\biggl\{\frac{f_1(x_3)}{f'_1(x_3)} + (a^* - x_3)\biggr\} + C_2,
$$
$$
- \frac{u_4}{f'_1(a^*)} = C_1\frac{f_1(a^*)}{f'_1(a^*)} + C_2
$$
and
$$
- \frac{u_1}{f'_1(x_1)} - v_1 =
C_1\biggl\{\frac{f_1(x_1)}{f'_1(x_1)} + (a^* - x_1)\biggr\} +C_2,
$$
$$
- \frac{u_2}{f'_1(x_2)} - v_2 =
C_1\biggl\{\frac{f_1(x_2)}{f'_1(x_2)} + (a^* - x_2)\biggr\} + C_2,
$$
from which it follows that, under the realization of the
conditions~(\ref{eq5.6}) and~(\ref{eq5.7}),
\begin{multline}\label{eq5.8}
C_1 = \biggl\{- \frac{- u_3}{f'_1(x_3)} + \frac{u_4}{f'_1(a^*)} -
v_3\biggr\}\biggl/\biggl\{\frac{f_1(x_3)}{f'_1(x_3)} -
\frac{f_1(a^*)}{f'_1(a^*)} + (a^* - x_3)\biggr\} =
\\
\biggl\{- \frac{u_2}{f'_1(x_2)} + \frac{u_1}{f'_1(x_1)} - v_2 +
v_1\biggl\}\biggl/\biggl\{\frac{f_1(x_2)}{f'_1(x_2)} -
\frac{f_1(x_1)}{f'_1(x_1)} - (x_2 -x_1)\biggr\},
\end{multline}
moreover,
$C_1$
does not depend on the location of the points
$x_1$,
$x_2$
and
$x_3$.

Further, we have
\begin{multline*}
\frac{u_4}{f'_1(a^*)} - \frac{u_3}{f'_1(x_3)} - v_3 =
\frac{u_4 - u_3}{f'_1(a^*)} + u_3\biggl(\frac{1}{f'_1(a^*)} -
\frac{1}{f'_1(x_3)}\biggr) - v_3 =
\\
\biggl\{- \frac{u_3}{f'_1(a^*)f'_1(x_3)} + O(a^* - x_3)\biggr\}
\{f'_1(a^*) - f'_1(x_3)\}.
\end{multline*}
Note that, here, we used the following estimates:
\begin{multline*}
\biggl|\frac{u_4 - u_3}{f'_1(a^*)} - v_3\biggr| =
\biggl|\frac{1}{f'_1(a^*)}
\int_{x_3}^{a^*} \frac{[f_1(t)]^2 d t}{\{1 +[f'_1(t)]^2\}^{1/2}}\,\,
- \int_{x_3}^{a^*} \frac{f'_1(t) d t}{\{1 +[f'_1(t)]^2\}^{1/2}}\biggr| =
\\
\frac{1}{f'_1(a^*)} \biggl|\int_{x_3}^{a^*}
\frac{f'_1(t)[f'_1(t) - f'_1(a^*)] d t}{\{1 + [f'_1(t)]^2\}^{1/2}}\biggr| \le
\frac{f'_1(a^*) -
f'_1(x_3)}{f'_1(a^*)}
\int_{x_3}^{a^*} \frac{f'_1(t) d t}{\{1 + [f'_1(t)]^2\}^{1/2}}
\\
\le \frac{f'_1(a^*) - f'_1(x_3)}{f'_1(a^*)}
\frac{f'_1(a^*)}{\{1 + [f'_1(a^*)]^2\}^{1/2}}.
\end{multline*}
From these calculations and~(\ref{eq5.5}), we will obtain, as a
result, the equality
\begin{multline}\label{eq5.9}
C_1 = \lim_{x_3 \to a^*}\biggl\{- \frac{u_3}{f'_1(x_3)} +
\frac{u_4}{f'_1(a^*)} - v_3\biggr\} \biggl/
\biggl\{\frac{f_1(x_3)}{f'_1(x_3)} - \frac{f_1(a^*)}{f'_1(a^*)} +
(a^* - x_3)\biggr\} =
\\
- \frac{u_4}{f_1(a^*)} \ne 0.
\end{multline}

By the analogy with~(\ref{eq5.9}) and on the basis
of~(\ref{eq5.8}), we can also to establish that
$$
C_1 = - \frac{u_1}{f_1(x_1)}.
$$
And by the convexity downwards of the function
$f_1$,
\begin{multline*}
\frac{u_1}{f_1(x_1)} = \frac{1}{f_1(x_1)} \int_0^{x_1}
\frac{[f'_1(t)]^2 d t}{\{1 + [f'_1(t)]^2\}^{1/2}} =
\\
\frac{1}{f_1(x_1)} \frac{f'_1(x_1)}{\{1 + [f'_1(x_1)]^2\}^{1/2}}
\int_{\xi}^{x_1} f'_1(t) d t =
\frac{f'_1(x_1)}{\{1 + [f'_1(x_1)]^2\}^{1/2}}
\frac{f_1(x_1) - f_1(\xi)}{f_1(x_1)}
\\
= \frac{f'_1(x_1)}{\{1 + [f'_1(x_1)]^2\}^{1/2}}
\frac{f_1(x_1) - f_1(\xi)}{f_1(x_1) - f_1(0)} \le
\frac{f'_1(x_1)}{\{1 + [f'_1(x_1)]^2\}^{1/2}} \to_{x_1 \to 0} 0
\end{multline*}
($0 < \xi < x_1$),
therefore,
$C_1 = 0$.
As the result of that, we have the contradiction
with~(\ref{eq5.9}). The latter, in turn, leads to the relations
$\rank N = \rank \widetilde{N} = 3$.
The proof of the lemma is completed.

In conclusion, note that the main results of our article were
earlier announced in~\cite{Kor1},~\cite{KK},~\cite{Ko5},~
\cite{Ko2},~\cite{Ko7}.

\vskip3mm

{\bf Acknowledgements}
\vskip3mm

The author was partially supported by the Russian Foundation for
Basic Research (Grant 11-01-00819-a), the Interdisciplinary
Project of the Siberian and Far-Eastern Divisions of the Russian
Academy of Sciences (2012-2014 no. 56), the State Maintenance
Program for the Leading Scientific Schools of the Russian
Federation (Grant NSh-921.2012.1) and the Exchange Program
between the Russian and Polish Academies of Sciences (Project
2014-2016).


\begin{thebibliography}{24}

\bibitem{Po}A.~V.~Pogorelov,
\emph{Extrinsic Geometry of Convex Surfaces}, AMS, Providence (1973).

\bibitem{Ko3}A.~P.~Kopylov,
Boundary values of mappings close to isometric mappings,
\emph{Siberian Math. J.}, \textbf{25}, no.3, 438-447 (1985).

\bibitem{Al}V.~A.~Aleksandrov,
Isometry of domains in
$\mathbb R^n$
and relative isometry of their boundaries,
\emph{Siberian Math. J.}, \textbf{25}, no.2, 339-347 (1985).

\bibitem{Al1}V.~A.~Aleksandrov,
Isometry of domains in
$\mathbb R^n$
and relative isometry of their boundaries. II,
\emph{Siberian Math. J.}, \textbf{26}, no.6, 783-787 (1986).

\bibitem{Ko1}A.~P.~Kopylov,
On the unique determination of domains in Euclidean spaces,
\emph{J. of Math. Sciences}, \textbf{153}, no.6, 869-898 (2008).

\bibitem{Kor1}M.~V.~Korobkov,
Necessary and sufficient conditions for the unique determination
of plane domains,
\emph{Dokl. Math.}, \textbf{76}, 722-723 (2007).

\bibitem{Kor2}M.~V.~Korobkov,
Necessary and sufficient conditions for unique determination of
plane domains,
\emph{Siberian Math. J.}, \textbf{49}, no.3, 436-451 (2008).

\bibitem{Kor4}M.~V.~Korobkov,
\emph{Some rigidity theorems in
Analysis and Geometry} [in Russian], Dis. Dokt. Fiz.-Mat. Nauk,
Novosibirsk (2008).

\bibitem{Kor3}M.~V.~Korobkov,
A criterion for the unique determination of domains in Euclidean
spaces by the metrics of their boundaries induced by the
intrinsic metrics of the domains, \emph{Siberian Advances in
Mathematics}, \textbf{20}, no.4, 256-284 (2010).

\bibitem{Bor}M.~K.~Borovikova,
On isometry of polygonal domains with boundaries locally
isometric in relative metrics,
\emph{Siberian Math. J.}, \textbf{33}, no.4, 571-580 (1993).

\bibitem{A}A.~D.~Aleksandrov,
\emph{Intrinsic Geometry of Convex Surfaces} [English translation],
Chapman\&Hall/CRC Taylor\&Francis Group, Boca Raton (2006).

\bibitem{Ko5}A.~P.~Kopylov,
A rigidity condition for the boundary of a submanifold in a
Riemannian manifold,
\emph{Dokl. Math.}, \textbf{77}, no.3, 340-341 (2008).

\bibitem{Ko6}A.~P.~Kopylov,
Unique determination of domains,
\emph{In{\rm:} Differential Geometry and Its Applications}.
Hackensack, NJ: World Sci. Publ., P.~157-169 (2008).

\bibitem{LW}F.~Leja, W.~Wilkosz,
Sur une propri\'{e}t\'{e} des domaines concaves,
\emph{Ann. Soc. Polon. Math.}, \textbf{2}, 222-224 (1924).

\bibitem{BZ}Yu.~D.~Burago, V.~A.~Zalgaller,
Sufficient conditions for convexity, \emph{In: Problems of Global
Geometry}, Leningrad: Nauka (1974). P.~3-53 (Zap. Nauchn. Sem.
LOMI; vol.~45).

\bibitem{KK}A.~P.~Kopylov and M.~V.~Korobkov,
Rigidity Conditions for the Boundaries of Submanifolds in
a Riemannian Manifold,
\emph{arXiv:1401.7295v2 [math.MG] 1 Oct 2014.}
http://arxiv.org/abs/1401.7295.

\bibitem{Sl}D.~A.~Slutskiy,
On Two Problems in the Theory of Unique Determination of Domains,
\emph{Vestnik, Quart. J. of Novosibirsk State Univ., Series:
Math., mech. and informatics}, \textbf{11}, no.2, 93-104 (2011).

\bibitem{Ko2}A.~P.~Kopylov,
Unique Determination of Domains by the Condition of Local
Isometry of Boundaries in the Relative Metrics,
\emph{Dokl. Math.}, \textbf{78}, no.2, 746-747 (2008).

\bibitem{Ko7}A.~P.~Kopylov,
On the Unique Determination of Domains by the Condition of the
Local Isometry of the Boundaries in the Relative Metrics
\emph{arXiv:1511.04235v6 [math.MG] 4 Oct 2016.}
http://arxiv.org/abs/1511.04235.

\end{thebibliography}
\end{document}